\documentclass[12pt]{amsart}

\usepackage{amsmath,amstext,amsthm,amssymb,amsfonts,amscd}
\usepackage{graphicx}
\usepackage{enumerate}

\usepackage[a4paper,left=3cm,right=3cm,top=2.52cm,bottom=2.52cm]{geometry}
\usepackage{fouriernc}
\usepackage{color}
\usepackage{tocvsec2}

\let\mathcal\relax
\usepackage[scr]{rsfso}
\let\mathcal\mathscr \let\mathscr\relax
\usepackage[mathscr]{euscript} 

\usepackage{tikz}
\usepackage{quiver}
\usepackage{environ}
\usepackage{etoolbox}

\usepackage[textsize=tiny]{todonotes}
\usepackage[colorlinks=true]{hyperref}          

\theoremstyle{plain}
\newtheorem{lem}{Lemme}[section]
\newtheorem{thm}[lem]{Theorem}
\newtheorem{prop}[lem]{Proposition}
\newtheorem{cor}[lem]{Corollary}

\theoremstyle{definition}
\newtheorem{defin}[lem]{Definition}
\newtheorem{exa}[lem]{Example}

\theoremstyle{remark}
\newtheorem{re}[lem]{Remark}

\DeclareMathOperator{\End}{\mathrm End}
\DeclareMathOperator{\Hom}{\mathrm Hom}

\def\temp{&} \catcode`&=\active \let&=\temp

\newcommand{\mytikzcdcontext}[2]{
  \begin{tikzpicture}[baseline=(maintikzcdnode.base)]
    \node (maintikzcdnode) [inner sep=0, outer sep=0] { 
    \begin{tikzcd}[#2]
        #1
    \end{tikzcd}
    };
  \end{tikzpicture}
}

\NewEnviron{tikzsd}[1][]{%
\def\myargs{#1}%
\edef\mydiagram{\noexpand\mytikzcdcontext{\expandonce\BODY}{\expandonce\myargs}}%
\mydiagram%
}

\makeatletter
\@namedef{subjclassname@2020}{%
  \textup{2020} Mathematics Subject Classification}
\makeatother

\numberwithin{equation}{section}
\numberwithin{figure}{section}

\begin{document}

\title[A new proof of Grauert's direct image theorem]{Grauert's direct image theorem via superconnections and desingularizations}

\author{Shu Shen}

\address{Institut de Math\'ematiques de Jussieu-Paris Rive Gauche, 
Sorbonne Universit\'e, Case Courrier 247, 4 place Jussieu, 75252 Paris Cedex 5, France.}
\email{shu.shen@imj-prg.fr}

\author{Jianqing Yu}
\address{School of Mathematical Sciences, University of Science and Technology of China, 96 Jinzhai Road, Hefei, Anhui 230026, P. R. China.}
\email{jianqing@ustc.edu.cn}

\dedicatory{}
\thanks{}

\date{\today}

\maketitle

\begin{abstract}
  We give a new differential-geometric proof of Grauert's theorem on the coherence of the higher direct image of a coherent sheaf under a proper holomorphic morphism between complex analytic spaces.
  In the smooth case, our approach is based on the antiholomorphic superconnection introduced by Block and further developed by Bismut-Shen-Wei.
  The required finiteness results follow from elliptic theory.
  In the singular case, we reduce the problem to the smooth setting using Hironaka's desingularization.
\end{abstract}

\tableofcontents

\section*{Introduction}

The aim of this paper is to give a new differential-geometric proof of Grauert's direct image theorem \cite{Grauert60}. 

Let $ \left(X,\mathscr{O}_{X} \right)$ and $ \left(Y,\mathscr{O}_{Y} \right)$ be two complex analytic spaces.
Let $ f: \left(X,\mathscr{O}_{X} \right) \to \left(Y,\mathscr{O}_{Y} \right)$ be a proper holomorphic map. 
If $ \mathcal{F} $ is a coherent $ \mathscr{O}_{X} $-module, 
let $ R f_{*} \mathcal{F} $ be the derived direct image of $ \mathcal{F} $.
If $ p\ge 0 $, denote $ R^{p} f_{*} \mathcal{F}$ the $ p$-th cohomology of $ R f_{*} \mathcal{F}$.
  
\begin{thm}\label{thm:vahu0}
  For all $ p\ge 0 $, the higher direct image $ R^{p} f_{*} \mathcal{F} $ is $ \mathscr{O}_{Y} $-coherent.
\end{thm}

Theorem \ref{thm:vahu0} is one of the fundamental results in complex geometry.
It was first established by Grauert \cite{Grauert60}, and has since been proved by numerous authors, including Forster-Knorr \cite{ForsterKnorr71}, Kiehl-Verdier \cite{KiehlVerdier71}, Levy \cite{Levy87directimage}, and recently by Clausen-Scholze \cite{ClausenScholze} (see also \cite[Chapter 10]{GrauertBook84} and \cite[Section IX.5]{DemaillyBook}).
Most of these proofs rely on deep and abstract tools from functional analysis, such as the theory of nuclear spaces and homological algebra on Fréchet complexes.

In this paper, we provide a new and more geometric proof of Theorem \ref{thm:vahu0}.
Our approach uses the antiholomorphic superconnection \cite{Block10,BismutShenWei23}, basic analysis of elliptic operators, and as a black box Hironaka’s desingularization theorem \cite{Hironaka64a,Hironaka64b}.
The proof extends the method given in \cite[Theorem 11.1.4]{BismutShenWei23}, where Theorem \ref{thm:vahu0} was proved under the assumption that $ X,Y$ are smooth and compact.

\subsection{Antiholomorphic superconnection}
Let $ X$ be a smooth compact complex manifold with holomorphic and antiholomorphic tangent bundles $ TX,\overline{TX}$.

The antiholomorphic superconnection on $ X$ is a fundamental tool introduced by Block \cite{Block10}.
If $ \left(D,v_{0 } \right)$ is a bounded holomorphic complex of vector bundles on $ X$, if $ \nabla^{D\prime \prime } $ is the holomorphic structure on $ D$, then the operator $ v_{0} + \nabla^{D\prime \prime } $ acting on the space $ \Omega^{0,\bullet} \left(X,D\right)$ of smooth antiholomorphic $ D$-valued forms is an example of antiholomorphic superconnection.
In general, the antiholomorphic superconnection does not necessarily arise from a holomorphic complex of vector bundles.
It is an operator $ A^{E\prime \prime } $ acting on $ C^{\infty} \left(X,E\right)$, where
\begin{enumerate}[\indent a)]
  \item $ E$ is a smooth $ \mathbf{Z} $-graded vector bundle of free $ \Lambda \left(\overline{T^{*} X} \right)$-module, such that there exist a smooth $ \mathbf{Z} $-graded vector bundle $ D$ and a (non-canonical) identification 
  \begin{align}\label{eq:nxco}
   E \simeq \Lambda \left(\overline{T^{*} X} \right) \widehat{\otimes} D;
  \end{align}
  \item $ A^{E\prime \prime }$ has total degree $ 1$, satisfies the integrability condition $ \left(A^{E\prime \prime } \right)^{2} = 0 $, and obeys the antiholomorphic Leibniz rules, so that under the identification \eqref{eq:nxco}, 
  \begin{align}\label{eq:nw2r}
   A^{E\prime \prime }= v_{0 } + \nabla^{D\prime \prime } + v_{\ge 2},
  \end{align}
  where $  \nabla^{D\prime \prime }$ is the antiholomorphic part of some connection and $v_{0 }, v_{\ge 2} $ are respectively sections of $\End\left(D\right) $ and $ \Lambda^{\ge 2} \left(\overline{T^{*} X} \right) \widehat{\otimes} \End\left(D\right)$ of total degree $ 1$.
\end{enumerate} 
If $ \left(E,A^{E\prime \prime } \right)$ is an antiholomorphic superconnection, the germ of smooth sections of $ E$ together with the restriction of $ A^{E\prime \prime } $ induces an $ \mathscr{O}_{X} $-complex of sheaves $ \left(\mathcal{E},A^{\mathcal{E} \prime \prime } \right) $.

Two key observations due to Block \cite{Block10} are the following.
First, locally near any point of $ X$, $ \left(E,A^{E\prime \prime } \right)$ is gauge equivalent to an antiholomorphic superconnection arising from a holomorphic complex of vector bundles, and therefore the complex $ \left(\mathcal{E},A^{\mathcal{E} \prime \prime } \right) $ has coherent cohomology.
Second, any coherent sheaf, or more generally any $ \mathscr{O}_{X} $-complex with bounded coherent cohomology, is derived isomorphic to an $ \mathscr{O}_{X} $-complex $ \left(\mathcal{E},A^{\mathcal{E} \prime \prime } \right) $ associated to some antiholomorphic superconnection.
These observations allow him to formulate a category equivalence between the homotopy category of antiholomorphic superconnections and the derived category $ {\rm D^{b}_{coh}}\left(X\right) $ of $ \mathscr{O}_{X} $-complexes with bounded coherent cohomology.

The above observations are rigorously established in \cite[Chapters 5 and 6]{BismutShenWei23}.
Several applications on compact complex manifolds are given:
\begin{enumerate}[\indent a)]
\item the Chern character for coherent sheaves in real Bott-Chern cohomology is constructed in \cite[Chapter 8]{BismutShenWei23};
\item a new proof of Grauert’s direct image theorem is given in \cite[Section 11.1]{BismutShenWei23};
\item a Riemann-Roch-Grothendieck formula for coherent sheaves in real Bott-Chern cohomology is established in \cite[Chapters 9-16]{BismutShenWei23}.
\end{enumerate}

Note that a) and b) are also obtained in \cite{Qiang1,BondalRoslyi23,Papayanov24} by methods similar to those of \cite{BismutShenWei23}.
See \cite{MTTW} and \cite{X25}, for other related works on antiholomorphic superconnections.
In rational Bott-Chern cohomology, a), as well as c), with respect to a more restrictive holomorphic map, are addressed in \cite{Wu23} through different methods inspired by \cite{Grivaux10}.

\subsection{The case where $ X$ and $ Y$ are compact manifolds}\label{s:XYcompact}
Assume $ X,Y$ are compact complex manifolds.
Recall briefly the proof of Theorem \ref{thm:vahu0} in this case given in \cite[Section 11.1]{BismutShenWei23}.

Let $ i_{f}  : x\in  X\to \left(x,f\left(x\right)\right)\in X\times Y$ be the graph of $ f$.
Let $ \pi : X \times Y \to Y$ be the obvious projection.
Then, $ f = \pi \circ i_{f} $.
It is easy to see that the conclusion of Theorem \ref{thm:vahu0} holds for $ i_{f} $, and the proof of Theorem \ref{thm:vahu0} can thus be reduced to establishing the analogous result for $ \pi $.

If $ \mathcal{F} $ is a coherent sheaf on $ X\times Y$, using Block's theorem \cite{Block10} and \cite[Theorem 6.5.1]{BismutShenWei23}, we can replace $ \mathcal{F} $ by the complex of sheaves $ \left(\mathcal{E},A^{\mathcal{E} \prime \prime } \right) $ associated to some antiholomorphic superconnection $ \left(E,A^{E\prime \prime } \right)$ on $ X \times Y$.
The derived direct image $ R\pi_{*} \mathcal{F} $ is then isomorphic to the direct image $ \left(\pi_{*} \mathcal{E},A^{\pi_{*} \mathcal{E}\prime \prime }  \right) $.

If $ V \subset Y$ is an open subset of $ Y$, we have
\begin{align}\label{eq:rshb}
 \pi_{*} \mathcal{E}\left(V\right)= C^{\infty} \left( X \times V,E_{|X \times V} \right).
\end{align}
Analysis tools such as differential operators can be used to study the sheaf $\pi_{*}\mathcal{E}$.  
Fixing a non-canonical identification \eqref{eq:nxco}, and given metrics on $D$ and $TX$, one constructs a Laplace-type operator on $\pi_{*}\mathcal{E}$ \cite[(11.1.2)]{BismutShenWei23}.  
This is a second-order differential operator acting on \eqref{eq:rshb}, which commutes with $A^{\pi_{*}\mathcal{E}\prime\prime}$.  
It differentiates only along $X$ and is elliptic in this direction.

By ellipticity and the spectral theory of this Laplacian, it is shown in \cite[Theorem 11.1.3]{BismutShenWei23} that with respect to a suitable spectral parameter $a > 0$,  locally on $Y$, there is a decomposition of complexes of sheaves,
\begin{align}\label{eq:b2xl}
  \pi_{*}\mathcal{E} = \left(\pi_{*}\mathcal{E}\right)_{a,-} \oplus \left(\pi_{*}\mathcal{E}\right)_{a,+}.
\end{align}
Moreover, $\left(\pi_{*}\mathcal{E}\right)_{a,-}$ is the complex of sheaves associated with some antiholomorphic superconnection defined locally on $Y$, while $\left(\pi_{*}\mathcal{E}\right)_{a,+}$ is acyclic.  
This yields the coherence of $\pi_{*}\mathcal{E}$, as stated in \cite[Theorem 11.1.4]{BismutShenWei23}.  

We emphasize the similarity of the above argument with the heat kernel proof of the local family index theorem due to Bismut \cite{B86}, as well as with the analytic theory of determinant line bundles initiated by Quillen \cite{Quillendeter} and further developed by Bismut–Gillet–Soulé \cite{BGS3}.

\subsection{The case where $ X$ and $ Y$ are non-compact manifolds}\label{s:XYnoncompact}
Assume now that $X$ and $Y$ are non-compact complex manifolds.  
To generalize the proof described in Section \ref{s:XYcompact} to this setting, we need to extend Block's theory to non-compact manifolds and develop a spectral truncation argument for the fibration $\pi : X \times Y \to Y$ with non-compact fibre $X$.  

Observe that since $f: X \to Y$ is proper, if $V \subset Y$ is a relatively compact open subset, then $f^{-1}(V)$ is relatively compact in $X$.  
If $\mathcal{F}$ is a coherent sheaf on $X$, the support of $\left(i_{f*}\mathcal{F}\right)_{|X \times V}$ is contained in $f^{-1}(V) \times V$, which is relatively compact in $X \times Y$.  
In this situation, Block’s theory can be developed (see, e.g., \cite{ChuangHolsteinLazarev21}): namely, $\left(i_{f*}\mathcal{F}\right)_{|X \times V}$ is derived isomorphic to the sheaf associated with some antiholomorphic superconnection on $X \times V$, which is acyclic on $\left(X \setminus f^{-1}(V)\right) \times V$.  
However, due to the lack of control on $\left(X \setminus f^{-1}(V)\right) \times V$, the associated Laplace-type operator does not appear to have good analytic properties.  

For this reason, in Section \ref{S:AntiTame}, we introduce a new object, the so-called \emph{tame antiholomorphic superconnection}.  
Near infinity, a tame object is induced from a constant complex of vector bundles.  
We extend Block’s theory to this tame framework (Theorem \ref{thm:poev}).  
In particular, $\left(i_{f*}\mathcal{F}\right)_{|X \times V}$ is derived isomorphic to the sheaf associated with a tame antiholomorphic superconnection.
 
Thanks to tameness, if suitable metrics are chosen, the Laplace-type operator on $\pi_{*}\mathcal{E}$ acquires a very simple structure on $ \left(X \backslash f^{-1} \left(V\right)\right)\times V$ : it becomes independent of the base variable in $ V$ and is given by the classical Kodaira Laplacian on $ X \backslash f^{-1} \left(V\right)$ shifted by a constant positive definite matrix.  

To apply the spectral truncation argument, we introduce a subcomplex $\left(\pi_{*}\mathcal{E}\right)_{(2)}$ of $\pi_{*}\mathcal{E}$, which is quasi-isomorphic to $\pi_{*}\mathcal{E}$ and fibrewise is the operator core of the Laplacian.
Since the Laplacian near infinity is independent of the base variable, its fibrewise resolvent and spectral projections extend naturally to the sheaf $\left(\pi_{*}\mathcal{E}\right)_{(2)}$.  
In this way, we obtain an analogue of the decomposition \eqref{eq:b2xl} for $\left(\pi_{*}\mathcal{E}\right)_{(2)}$ with similar properties (Theorem \ref{thm:mitj} and Proposition \ref{prop:sn2o}).  
This allows us to extend \cite[Theorem 11.1.4]{BismutShenWei23} to the non-compact setting.

\subsection{The case where $ X$ and $ Y$ have singularities}
Assume now that $X$ and $Y$ are not smooth.  
Since coherence is a local property, the singularities of $Y$ can be handled by applying an extension principle for coherent sheaves \cite[Section 1.2.7]{GrauertBook84}.

This local property can also be used to reduce the proof of Theorem \ref{thm:vahu0} to the case where $X$ is second countable and of finite dimension.  

In Section \ref{s:more}, using Rückert’s Nullstellensatz with respect to a  closed reduced analytic subspace $Z \subset X$, we show that the coherence of higher direct images for coherent sheaves defined on $ X$ supported on $Z$ follows from Grauert’s theorem for the restriction map $f_{|Z} : Z \to Y$ (Proposition \ref{prop:izht}).  

We apply this result in two cases: when $Z$ is the reduction $X_{\mathrm{red}}$, and when $Z$ is the singular locus $X_{\mathrm{sing}}$.  
The former allows us to reduce the proof of Theorem \ref{thm:vahu0} for $f: X \to Y$ to the restriction map $f_{|X_{\mathrm{red}}} : X_{\mathrm{red}} \to Y$ (Corollary \ref{cor:tnlv}).
The latter, combined with Hironaka’s resolution of singularities \cite{Hironaka64a,Hironaka64b} and the smooth case of Theorem~\ref{thm:vahu0}, allows us to further reduce the proof of Theorem \ref{thm:vahu0} to $f_{|X_{\mathrm{sing}}} : X_{\mathrm{sing}} \to Y$ (see Section \ref{s:endp}).  

Ultimately, our main result can be established by induction on $\dim X$.

\subsection{Some remarks}
Since our goal is to provide a new proof of a classical result, we try to make the argument as self-contained as possible.
Nevertheless, the following standard results, mostly consequences of the Weierstrass preparation theorem, are used, sometimes implicitly.
\begin{enumerate}[\indent a)]
  \item Oka's coherence theorem : $ \mathscr{O}_{X} $ is coherent \cite[Section 2.5.3]{GrauertBook84}.
  \item The nilradical $ \mathscr{N}_{X} \subset \mathscr{O}_{X} $ is coherent, so that the reduction $ X_{{\rm red}} $ is a closed complex analytic subspace of $ X$ \cite[Section 4.2.5]{GrauertBook84};
  \item The singular locus $ X_{{\rm sing}} $ is a closed complex analytic subspace of $ X$ \cite[Section 6.2.2]{GrauertBook84};
  \item The inverse image preserves the coherence \cite[Section 1.2.6]{GrauertBook84};
  \item The direct image of a closed embedding preserves the coherence \cite[Section 1.2.7]{GrauertBook84}
\end{enumerate}
We note that a) has already been used in \cite{Block10,BismutShenWei23}.

With potential applications in mind, in Section \ref{S:mainresult}, we reformulate our main result in derived category, giving a slightly more general version as Theorem \ref{thm:main}.
Basic results on the derived category are therefore used throughout the paper.

In Section \ref{S:AntiTame}, we develop a complete theory of the tame antiholomorphic superconnection.
For readers interested only in a proof of Grauert’s theorem, only part of the results in this section are required.

All of the above makes the paper somewhat long, but we believe the added depth justifies it.

\subsection{Organisation of the paper}
The paper is organised as follows.
In Section \ref{S:mainresult}, we state the main result of the paper.

In Section \ref{S:AntiTame}, we introduce a tameness condition for antiholomorphic superconnections on non-compact manifolds. 
We then generalise Block's theory \cite{Block10} and \cite[Theorem 6.5.1]{BismutShenWei23} to the non-compact setting.

In Section \ref{S:product}, we study the direct image of an antiholomorphic superconnection on a product manifold and show Theorem \ref{thm:vahu0} in the smooth case.

Finally, in Section \ref{S:Gg}, we establish Theorem \ref{thm:vahu0} in fully generality.

In the whole paper, we use the $ \mathbf{Z} $-graded formalism.
We denote the graded commutator by $[,]$.
The symbol $\widehat{\otimes}$ stands for the $\mathbf{Z}$-graded tensor product, which should not be confused with the tensor product of Fréchet spaces.
We also adopt the convention that $ \mathbf{N} = \{0,1,2, \ldots \}$ and $ \mathbf{N}^{*}  = \{1,2, \ldots \}$.

\subsection{Acknowledgement}
This work was initiated following a suggestion by Xiaojun Huang during the Complex Geometry Conference held in Taipei in December 2023. 
We are deeply grateful to him for his valuable insights and encouragement. 
We also wish to thank Frédéric Paugam for his help in clarifying several points on derived categories.
We acknowledge the partial financial support from NSFC under Grant No. 12371054. S.S. is also partially supported by ANR Grant ANR-20-CE40-0017.

\settocdepth{subsection}
 
\section{Statement of the main result}\label{S:mainresult}
The purpose of this section is to review some basic results on complex analytic spaces and to state Grauert's direct image theorem.

This section is organised as follows.
In Section \ref{s:cas}, we recall some basic definition related to complex analytic spaces.
 
In Section \ref{s:dc}, we introduce the derived categories and the derived direct image functors.

In Section \ref{s:gdi}, we state Grauert's direct image theorem.

Finally, in Section \ref{s:imr}, we give some immediate reductions for the proof of Grauert's direct image theorem.

\subsection{Complex analytic spaces}\label{s:cas}
Let $ \left(X,\mathscr{O}_{X} \right)$ be a complex analytic space\footnote{It is called complex analytic schemes in \cite[Section II.9]{DemaillyBook}.} \cite[Section 1.1.5]{GrauertBook84}.

Then, $ X$ is a Hausdorff (not necessarily second countable) topological space, and $ \mathscr{O}_{X} $ is a sheaf of $ \mathbf{C} $-algebras on $ X$, so that $ \left(X,\mathscr{O}_{X} \right)$ is a $ \mathbf{C} $-ringed space.
Moreover, for each $ x_{0 } \in X$, there exist 
\begin{enumerate}[\indent a)]
  \item an open neighbourhood $ U \subset X$ of $ x_{0 } $,
  \item an open subset $ \Omega $ of $ \mathbf{C}^{n} $ with a sheaf of holomorphic functions $ \mathscr{O}_{\Omega }$, 
  \item a coherent ideal sheaf $ \mathscr{I} \subset \mathscr{O}_{\Omega } $ with the zero set $ N\left(\mathscr{I} \right) \subset \Omega $,
\end{enumerate}
such that we have an identification of $ \mathbf{C} $-ringed spaces, 
\begin{align}\label{eq:d1oi}
 \left(U,\mathscr{O}_{X|U} \right)= \left(N\left(\mathscr{I} \right),\left(\mathscr{O}_{\Omega } / \mathscr{I}\right)_{|N\left(\mathscr{I} \right)}  \right).
\end{align}

If $ x\in X$, we denote $ \dim_{x} X\in \mathbf{N} $ the  dimension \cite[Section 5.1.1]{GrauertBook84} of $ X$ at $ x$.
Put
\begin{align}\label{eq:ptlf}
 \dim X= \sup_{x\in X} \dim_{x} X\in \mathbf{N} \cup \left\{\infty\right\}.
\end{align}
By \cite[p.~94]{GrauertBook84}, the function $x\in X \to \dim_{x} X\in \mathbf{N} $ is upper semi-continuous. Thus, if $ 
X$ is compact or if $ X$ is a relatively compact open subspace of some complex analytic space, then  $ \dim X$ is finite.

 

  
Let $ \left(Y,\mathscr{O}_{Y} \right)$ be another complex analytic space.
Let $ \left(f,\widetilde{f} \right): \left(X,\mathscr{O}_{X} \right)\to \left(Y,\mathscr{O}_{Y} \right)$ be a holomorphic map, i.e., $ \left(f,\widetilde{f} \right)$  is a morphism of $ \mathbf{C} $-ringed spaces. 
Then, $ f : X \to Y $ is a continuous map between the underlying topological spaces, and $ \widetilde{f} : \mathscr{O}_{Y} \to f_{*} \mathscr{O}_{X} $ is a morphism of sheaves of $ \mathbf{C} $-algebras.
If $ x_{0 } \in X$ and if $ y_{0 } = f\left(x_{0 } \right)\in Y$, on the stalks, $ \widetilde{f} $ induces a morphism of $ \mathbf{C} $-algebras,\footnote{It is defined by the composition of $\widetilde{f}_{y_{0 } } : \mathscr{O}_{Y,y_{0 } } \to \left(f_{*} \mathscr{O}_{X}\right)_{y_{0 } }   $ and the resection $ \left(f_{*} \mathscr{O}_{X}\right)_{y_{0 } } \to \mathscr{O}_{X,x_{0 } } $.}
\begin{align}\label{eq:hdvt}
 \widetilde{f}_{x_{0 }, y_{0 } } : \mathscr{O}_{Y,y_{0 } } \to \mathscr{O}_{X,x_{0 } }.
\end{align}
Recall that the stalks $ \mathscr{O}_{X,x_{0 } }$,  $ \mathscr{O}_{Y,y_{0 } }$ are local $ \mathbf{C} $-algebras.
If $ \mathfrak{m} \left(\mathscr{O}_{X,x_{0 } }\right)$,  $ \mathfrak{m} \left(\mathscr{O}_{Y,y_{0 } }\right)$ are the maximal ideals of $ \mathscr{O}_{X,x_{0 } }$, $ \mathscr{O}_{Y,y_{0 } }$, we have the canonical splittings,  
\begin{align}\label{eq:dkmn}
& \mathscr{O}_{X,x_{0 } }= \mathbf{C} \oplus \mathfrak{m} \left(\mathscr{O}_{X,x_{0 } }\right), &\mathscr{O}_{Y,y_{0 } }= \mathbf{C} \oplus \mathfrak{m} \left(\mathscr{O}_{Y,y_{0 } }\right).
\end{align}
By construction, \eqref{eq:hdvt} preserves the splittings \eqref{eq:dkmn}, so that $ \widetilde{f}_{x_{0 }, y_{0 } } $ is automatically a morphism of local $ \mathbf{C} $-algebras.
Therefore, $ \left(f,\widetilde{f} \right)$ is a morphism of local $ \mathbf{C} $-ringed spaces.

In the sequel, for ease of the notation, we will write $ X, Y$ for the complex analytic spaces, and $ f$ for the holomorphic map.

\subsection{The derived categories}\label{s:dc}
Let $ {\rm M}\left(X\right) $ be the category of $ \mathscr{O}_{X} $-modules and $ \mathscr{O}_{X} $-morphisms.
Then, $ {\rm M}\left(X\right)$ is an abelian category.
Let $ {\rm D} \left(X \right)$ be the derived category of $ \mathscr{O}_{X} $-complexes.
If $ p \in \mathbf{Z} $, let $ \mathcal{H}^{p} $ be the $ p$-th cohomology functor from $ {\rm D} \left(X\right) $ to $  {\rm M}\left(X\right)$.

Let $ {\rm D}^{+}  \left(X \right)$ be the full subcategory of $ {\rm D} \left(X \right)$ consisting of objects $ \mathcal{F} $ whose cohomology is bounded from below, i.e., if $ p \ll -1$, $ \mathcal{H}^{p} \mathcal{F} = 0 $.
Let ${\rm D_{coh}^{b}}\left(X\right)$ be the full subcategory of $ {\rm D} \left(X\right) $ consisting of objects with bounded coherent cohomology, i.e., all the $ \mathcal{H}^{p} \mathcal{F} $ are $ \mathscr{O}_{X} $-coherent and for $ |p|\gg 1$, $ \mathcal{H}^{p} \mathcal{F} = 0 $.
Let $   {\rm D_{coh}^{lb,+}}\left(X\right) $ be the full subcategory of $ {\rm D}^{+}  \left(X \right)$ consisting of objects $ \mathcal{F} $, whose restriction to any relatively compact open subset $ U$ of $ X$ is an object of $ {\rm D^{b}_{coh}}\left(U\right)$.

Let $ f: X\to Y$ be a holomorphic map.
The direct image $ f_{*} $ defines a covariant functor from $ {\rm M}\left(X\right)$ to $ {\rm M}\left(Y\right)$.
The functor $ f_{*} $ is left exact, and in general not right exact.
Let $ R f_{*} :{\rm D}^{+} \left(X\right) \to  {\rm D}^{+} \left(Y\right) $ be the associated right derived functor.
If $ \mathcal{F} $ is an object of $ {\rm D}^{+} \left(X\right)$, for $ p \in \mathbf{Z} $, denote  $ R^{p} f_{*} \mathcal{F} $ the $ p$-th cohomology of $ R f_{*} \mathcal{F}$. 
Then, $ R^{p} f_{*} \mathcal{F} $ is an $ \mathscr{O}_{Y} $-module.

Recall briefly the construction of $ Rf_{*} \mathcal{F} $.
Since $ {\rm M}\left(X\right)$ has enough injectives, if $ \mathcal{F} $ is an object of $ {\rm D}^{+}\left(X\right)$, by \cite[\href{https://stacks.math.columbia.edu/tag/013I}{Tag 013I}, \href{https://stacks.math.columbia.edu/tag/013K}{Tag 013K}]{stacks-project}, there exist an object $ \mathcal{I} $ in $ {\rm D}^{+}\left(X\right)$, which is bounded from below and injective degree by degree, and a quasi-isomorphism of $ \mathscr{O}_{X} $-complexes,
\begin{align}\label{eq:evfv}
 \mathcal{F} \to \mathcal{I}. 
\end{align}
By \cite[\href{https://stacks.math.columbia.edu/tag/013P}{Tag 013P}, \href{https://stacks.math.columbia.edu/tag/013S}{Tag 013S}]{stacks-project}, in the homotopy category of $ \mathscr{O} _{X} $-complexes, up to canonical isomorphism, the above $ \mathcal{I} $ and quasi-isomorphism \eqref{eq:evfv} are unique.
By \cite[\href{https://stacks.math.columbia.edu/tag/05TH}{Tag 05TH}]{stacks-project}, we have a canonical isomorphism in $ {\rm D}^{+} \left(X\right) $,  
\begin{align}\label{eq:simg}
 Rf_{*} \mathcal{F} \simeq  f_{*} \mathcal{I}. 
\end{align}
Moreover, given $ c_{0 } \in \mathbf{Z} $, if $ \mathcal{H} \mathcal{F} $ is concentrated in degrees $ [c_{0 } ,+\infty)$, by \cite[\href{https://stacks.math.columbia.edu/tag/05T6}{Tag 05T6}]{stacks-project}, the above $ \mathcal{I} $ can be chosen to be concentrated in the same range, so that the cohomology of $ Rf_{*} \mathcal{F} $ remains concentrated in degrees $ [c_{0 } ,+\infty)$.

If $ f$ is a closed embedding, by \cite[Section 1.2.7, Properties of holomorphic embeddings 1)]{GrauertBook84}, $ f_{*} $ is exact, so that
\begin{align}\label{eq:gpfo}
 Rf_{*} = f_{*}. 
\end{align}
By the above exactness, for each $ p\in \mathbf{Z} $, we have 
\begin{align}\label{eq:bzic}
  R^{p} f_{*}\mathcal{F} = f_{*} \left( \mathcal{H}^{p} \mathcal{F} \right). 
\end{align}
By \eqref{eq:bzic}, $ R^{p} f_{*}\mathcal{F}= 0 $ if and only if $ \mathcal{H}^{p} \mathcal{F} = 0 $.
Moreover, by \cite[Section 1.2.7, Properties of holomorphic embeddings 2)]{GrauertBook84}, $ R^{p} f_{*}\mathcal{F}$ is coherent if and only if $ \mathcal{H}^{p} \mathcal{F} $ is coherent.
In particular, $ Rf_{*} $ sends $ {\rm D^{b}_{coh}} \left(X\right)$ to $ {\rm D^{b}_{coh}} \left(Y\right)$.



\subsection{Grauert's direct image theorem}\label{s:gdi}
Recall that $ f:X\to Y$ is called proper, if the preimage of a compact subset in $ Y$ is compact.
In this sequel, we will always assume that $ f$ is a proper holomorphic map.

The purpose of this paper is to give a new proof for the following theorem of Grauert \cite{Grauert60}.

\begin{thm}\label{thm:main}
 The following statements hold.
 \begin{enumerate}[\indent 1)]
   \item \label{main1} The functor $ Rf_{*} $ restricts from $ {\rm D^{lb,+}_{coh}}\left(X\right) $  to $ {\rm D^{lb,+}_{coh}}\left(Y\right)$.
   In particular, Theorem \ref{thm:vahu0} holds.
   \item \label{main2} Suppose $ \dim X < + \infty $.
   \begin{enumerate}[a)]
     \item \label{main2a} There exists $ c\left(f\right)\in \mathbf{N} $ such that 
     if $ c_{0} ,c_{1} \in \mathbf{Z} $ with $ c_{0 } \le c_{1} $, if $ \mathcal{F} $ is an object in $ {\rm D^{b}_{coh}}\left(X\right) $ whose cohomology is concentrated in degree $ \left[c_{0} ,c_{1} \right]$, then the cohomology of $ R  f_{*} \mathcal{F} $ is coherent and concentrated in degree $ \left[c_{0} ,c_{1}+c\left(f\right)\right]$.
     \item \label{main2b} The functor $ Rf_{*} $ restricts from $ {\rm D^{b}_{coh}}\left(X\right) $  to $ {\rm D^{b}_{coh}}\left(Y\right)$.
   \end{enumerate}
 \end{enumerate}
\end{thm}
\begin{proof}
 We will show in Corollary \ref{cor:mq1y} that Statement \ref{main1}) is a consequence of \ref{main2b}).
 
 Statement \ref{main2a}) will be established in several steps.
 In Section \ref{s:imr}, we reduce the proof to the case where $ X$ has second countable topology and $ Y$ is a smooth manifold.
 In Section \ref{S:Gs}, we consider the case when $ X$ is also smooth.
 In Section \ref{S:Gg}, we complete the proof in full generality. 

 Statement \ref{main2b}) is an immediate consequence of \ref{main2a}).
\end{proof}


\begin{re}\label{re:pnjk}
   If $ f$ is an embedding, by the discussion following \eqref{eq:bzic},  our theorem holds with  
  \begin{align}\label{eq:itxp}
    c\left(f\right)= 0. 
  \end{align}
\end{re}

\begin{re}\label{re:d2jg}
  If $ \dim X=+\infty$, $ R f_{*} \mathscr{O}_{X}  $ may have unbounded cohomology.
  In fact, if $ g\in \mathbf{N}^{*} $, let $ \Sigma_{g} $ be a compact Riemann surface of genus $ g$.
  Then, $\dim H^{1} \left(\Sigma_{g}, \mathscr{O}_{\Sigma_{g} } \right) = g\ge 1$, so that  $ H^{1} \left(\Sigma_{g}, \mathscr{O}_{\Sigma_{g} } \right) \neq 0 $.
  If $ n\in \mathbf{N}^{*}  $, denote $ \Sigma_{g}^{n} $ the $ n$-th product of $ \Sigma_{g} $. 
  Put 
  \begin{align}\label{eq:d1ui}
 &  X= \bigsqcup_{n\in \mathbf{N}^{*}  }  \Sigma_{g}^{n},& Y= \mathbf{N}^{*}. 
  \end{align}
  Let $ f : X\to Y$ be the map, which sends $ \Sigma_{g}^{n}  $ to $ n\in  Y$.
  If $ p\in \mathbf{N} $, then $ R^{p} f_{*}  \mathscr{O}_{X} $ is a vector bundle on $ Y$.
  Over $ p\in Y$,  we have  
  \begin{align}\label{eq:yrbt}
   \left(R^{p} f_{*}  \mathscr{O}_{X}\right)_{|p}  = H^{p} \left(\Sigma_{g}^{p} ,\mathscr{O}_{\Sigma_{g}^{p}  }  \right) =  \left(H^{1} \left(\Sigma_{g} ,\mathscr{O}_{\Sigma_{g}  }  \right)\right)^{p} \neq0. 
  \end{align}
  Therefore, for all $ p\in \mathbf{N}^{*} $, we have $ R^{p}  f_{*}  \mathscr{O}_{X}\neq0  $. 
 \end{re}

\begin{re}\label{re:qbql}
  When $ X$ is smooth, in Section \ref{S:Gs}, we show 
  \begin{align}\label{eq:wghp}
    c\left(f\right) \le \dim X.
  \end{align}
  When $ X$ is non-smooth, our proof given in Section \ref{s:endp} does not provide an optimal control on $ c\left(f\right)$. 
  In Section \ref{s:cont}, which is included for expository purposes, using Andreotti-Grauert’s theorem \cite{AndreottiGrauert62} \cite[IX (4.15) Corollary]{DemaillyBook}, we show that \eqref{eq:wghp} still holds for arbitrary complex analytic space.
\end{re}
 

 
\subsection{Some immediate reductions : coherence as a local property}\label{s:imr}
Let $f: X \to Y$ be a proper holomorphic map, and let $ \mathcal{F} $ be an object in $ {\rm D^{+}} \left(X\right)$.

Let $ V $ be an open subset of $ Y$. 
Set 
\begin{align}\label{eq:peew}
 U= f^{-1} \left(V\right).
\end{align}
Then, $ U$ is an open subset of $ X$.
Let  
\begin{align}\label{eq:rboa}
 f_{U,V} : U \to V
\end{align}
be the restriction of $ f$ to $ U$.
Then, $ f_{U,V} $ is still a proper holomorphic map.

\begin{prop}\label{prop:dpzs}
  We have the canonical identification in $ {\rm D^{+}} \left(V\right)$, 
  \begin{align}\label{eq:ebfi}
    \left(Rf_{*} \mathcal{F}\right)_{|V} =    Rf_{U,V*}\left( \mathcal{F}_{|U}\right). 
   \end{align}
\end{prop}
\begin{proof}   
Let $ \mathcal{F} \to \mathcal{I} $ be a right injective resolution for $ \mathcal{F} $ as in \eqref{eq:evfv}.
By \cite[Proposition 2.4.1 (i)]{KashiwaraSchpira}, $ \mathcal{F} _{|U} \to \mathcal{I} _{|U} $ is a right injective resolution for $ \mathcal{F}_{|U} $.
By \eqref{eq:simg}, we have 
\begin{align}\label{eq:a1ty}
\left(Rf_{*} \mathcal{F}\right)_{|V}= \left(f_{*} \mathcal{I}\right)_{|V}=  f_{U,V*} \left(\mathcal{I} _{|U} \right)=  Rf_{U,V*}\left( \mathcal{F}_{|U}\right).
\end{align}
The proof of our proposition is complete.
\end{proof}

\begin{cor}\label{cor:mq1y}
  If for any relatively compact open subset $ V\subset Y$, the functor $ Rf_{U,V*} $ restricts from $ {\rm D^{b}_{coh}}\left(U\right) $  to $ {\rm D^{b}_{coh}}\left(V\right)$, then $ Rf_{*} $ restricts from $ {\rm D^{lb,+}_{coh}}\left(X\right) $  to $ {\rm D^{lb,+}_{coh}}\left(Y\right)$.
  In particular, Theorem \ref{thm:main} \ref{main1}) is a consequence of Theorem \ref{thm:main} \ref{main2b}).
\end{cor}
\begin{proof}
 Since $ f$ is proper, if $ V$ is relatively compact in $ Y$,  then $ U= f^{-1} \left(V\right)$ is relatively compact in $ X$.
 By this observation and by Proposition \ref{prop:dpzs}, we get the first statement of our corollary.
 Using the fact that $ U$ has finite dimension, we obtain the second statement, and finish the proof of our corollary.
\end{proof}




Assume now $ V$ is small enough.
By \eqref{eq:d1oi}, $ V$ is a closed analytic subspace of an open domain $ V_{1} $ in $ \mathbf{C}^{n} $.
Let 
\begin{align}\label{eq:zgvp}
 i: V \to V_{1}
\end{align}
be the natural closed (proper) embedding.
Set 
\begin{align}\label{eq:dine}
  f_{U,V_{1} } = i \circ f_{U,V}.
\end{align}
Then,  $ f_{U,V_{1} }$ is a proper holomorphic map.
 
\begin{prop}\label{prop:vimt}
 For $ p \in \mathbf{Z} $,  we have the canonical isomorphism of $ \mathscr{O}_{V_{1} } $-modules, 
  \begin{align}\label{eq:j2fs}
   i_{*} R^{p} f_{U,V*} \left(\mathcal{F}_{|U}  \right)= R^{p} f_{U,V_{1}* } \left(\mathcal{F} _{|U} \right).
   \end{align}
 In particular, $ R^{p} f_{U,V*} \left( \mathcal{F}_{|U} \right)$ is coherent (resp. vanishes) if and only if $ R^{p} f_{U ,V_{1}*} \left( \mathcal{F}_{|U } \right)$ is coherent (resp. vanishes).
\end{prop}
\begin{proof} 
  Since $ i$ is a closed embedding, by 
  \eqref{eq:gpfo}, \eqref{eq:dine}, and by functoriality \cite[\href{https://stacks.math.columbia.edu/tag/0D5T}{Tag 0D5T}]{stacks-project}, we have 
  \begin{align}\label{eq:gubv}
   Rf_{U,V_{1} *} = Ri_{*} R f_{U,V*}= i_{*} R f_{U,V*}.
  \end{align}
  By \eqref{eq:bzic} and \eqref{eq:gubv}, we get \eqref{eq:j2fs}.
The last statements follow immediately from the discussion after \eqref{eq:bzic}.
\end{proof}

\begin{re}\label{re:fq3s}
  If $ \mathcal{V} = \left\{V\right\}$ is a cover of $ Y$ by above small open subsets, 
  by Propositions \ref{prop:dpzs} and \ref{prop:vimt}, we see that the proof of Theorem \ref{thm:main} \ref{main2a}) can be reduced to show the corresponding results for the holomorphic maps $ f_{U,V_{1} } :U \to V_{1} $ together with
  \begin{align}\label{eq:xofp}
    \sup_{V\in \mathcal{V} } c\left(f_{U,V_{1} }\right)<+\infty.
  \end{align}
  Note that $ U$ is relatively compact in $ X$, so it has a second countable topology, and that $ V_{1} $ is a smooth manifold.
\end{re}

\section{Antiholomorphic superconnection and tameness}\label{S:AntiTame}
The purpose of this section is to generalise \cite[Theorem 4.3]{Block10} and \cite[Theorem 6.5.1]{BismutShenWei23} to coherent sheaves on non-compact complex manifolds with certain tame structure at infinity.
Essentially, our approach consists of gluing,  outside a compact set, a tame model to the antiholomorphic superconnection, a fundamental object introduced in \cite{Block10, BismutShenWei23}.
In these two references, results are formulated in derived category.
However, it is well-known that gluing 
in derived category is a delicate issue.
Here, we instead formulate our results in a certain homotopy category.


This section is organised as follows.
In Section \ref{s:BX}, we recall the definition and some basic constructions on antiholomorphic superconnections.

In Section \ref{s:BtX}, we introduce a tameness condition for the antiholomorphic superconnection on a non-compact complex manifold $ X$.

In Section \ref{s:Cbt}, we introduce a tameness condition for $ \mathscr{O}_{X} $-complexes.

In Section \ref{s:Dbt}, we introduce another tameness condition for $ \mathscr{O}_{X} $-complexes which are also modules over the Dolbeault complex.
We state Theorem \ref{thm:poev} the main result of this section.

Finally, in Sections \ref{sEssSuj} and \ref{s:thmb}, we show Theorem \ref{thm:poev}.

\subsection{Antiholomorphic superconnection}\label{s:BX}
Let $ X$ be a complex manifold of dimension $ n$.\footnote{To be consistent with \eqref{eq:ptlf}, this means that $ X$ has one $ n$-dimensional connected component, and the others have dimension $ \le n$.}
Let $ TX,\overline{TX} $ be the holomorphic and antiholomorphic tangent bundles of $ X$.
Let $ \left(\Omega^{0,\bullet } \left(X,\mathbf{C} \right),\overline{\partial}^{X} \right)$ be the Dolbeault complex of smooth antiholomorphic forms on $ X$.

Let $ \left(E,A^{E \prime\prime } \right)$ be an antiholomorphic superconnection \cite[Sections 4.1 and 5.1]{BismutShenWei23} on $ X$.
Let us recall its definition and some related constructions.

By definition, there exist $ r,r^{\prime } \in \mathbf{Z} $ with $ r \le r^{\prime } $ such that $ E= \bigoplus_{i= r}^{r'}  E^{i}$ is a $ \mathbf{Z} $-graded finite dimensional smooth vector bundle on $ X$.
Moreover, $ E$ is assumed to have a $ \mathbf{Z} $-graded free $ \Lambda \left(\overline {T^{*} X}\right)$-action.

Then, $ E$ is equipped with a canonical filtration
\begin{align}\label{eq:uucw}
  E \supset \overline{T^{*} X} \cdot E \supset \Lambda^{2} \left(\overline{T^{*} X} \right)\cdot E \supset \cdots \supset \Lambda^{n} \left(\overline{T^{*} X} \right)\cdot E \supset 0. 
\end{align}
Put  
\begin{align}\label{eq:qxmz}
  D= E/ \overline {T^{*} X} \cdot E.
 \end{align}
The multiplication induces a point-wise linear map, 
\begin{align}\label{eq:my2s}
  \Lambda^{\bullet } \left( \overline{T^{*} X}\right) \otimes  D \to \Lambda^{\bullet } \left(\overline{T^{*} X}\right) \cdot E/\Lambda^{\bullet+1} \left(\overline{T^{*} X} \right) \cdot E.
\end{align}
The freeness of the $ \Lambda \left(\overline{T^{*} X} \right)$-action on $ E$ means that all the objects or morphisms in \eqref{eq:uucw}-\eqref{eq:my2s}  are smooth vector bundles or smooth bundle maps, and that \eqref{eq:my2s} is an isomorphism of smooth vector bundles.
Therefore, we have a non-canonical isomorphism of smooth vector bundles of $ \mathbf{Z} $-graded $ \Lambda \left(\overline{T^{*} X} \right)$-modules, 
\begin{align}\label{eq:aayz}
 E\simeq \Lambda \left(\overline {T^{*} X}\right)\widehat{\otimes}  D.
\end{align}

Moreover, $ A^{E \prime\prime }: C^{\infty} \left(X,E\right)\to C^{\infty} \left(X,E\right) $ is a first order differential operator such that 
\begin{enumerate}[\indent a)]
  \item $ A^{E\prime\prime } $ increases the degree of sections of $ E$ by $ 1$; 
  \item $ A^{E\prime\prime } $ is flat, i.e., 
\begin{align}\label{eq:gpgm}
  \left(A^{E \prime\prime }\right)^{2}  = 0;
\end{align}
  \item $ A^{E\prime\prime } $ satisfies the antiholomorphic Leibniz rule,  i.e., for $ \alpha \in \Omega^{0,\bullet} \left(X,\mathbf{C} \right)$ and $ s\in C^{\infty}\left(X,E\right)$,  
\begin{align}\label{eq:vytt}
  A^{E \prime\prime } \left(\alpha \cdot s\right) = \overline{\partial}^{X} \alpha \cdot s + \left(-1\right)^{{\rm deg}\alpha } \alpha \cdot A^{E \prime\prime } s.
\end{align}
\end{enumerate}


By \eqref{eq:qxmz} and \eqref{eq:vytt}, $ A^{E\prime \prime } $ induces a 
smooth section $ v_{0 } $ of $ \End \left(D\right)$ of degree $ 1$, so that $ \left(D,v_{0 } \right)$ is a bounded smooth complex of vector bundles.

\begin{defin}\label{def:ms2r}
 The pair $ \left(D,v_{0 } \right)$ is called the diagonal of $ \left(E,A^{E\prime\prime } \right)$.
\end{defin}

\begin{defin}\label{def:aer3}
  The support of $ \left(E,A^{E\prime\prime } \right)$ is a closed subset of $ X$ consisting of points $ x\in X$ such that $ \left(D,v_{0 } \right)_{|x} $ is not exact.
\end{defin}

Non-canonically, using \eqref{eq:aayz}, as in \cite[(5.1.6)]{BismutShenWei23}, we can write 
\begin{align}\label{eq:ktja}
  A^{E \prime\prime } = v_{0 } +\nabla^{D\prime\prime } +v_{2} + \ldots + v_{n},
\end{align}
where $ \nabla^{D\prime\prime } $ is the antiholomorphic part of some connection on $ D$ and for $ i\ge 2$,  $ v_{i }\in  \Omega^{0,i} \left(X,\End^{1-i} \left(D\right)\right)$.

Let $\mathcal{D},  \mathcal{E} $ be complexes of sheaves, defined on an open subset $U \subset X$ by  
\begin{align}\label{eq:agac}
& \mathcal{D} \left(U\right)= \left(C^{\infty} \left(U,D_{|U} \right),v_{0 |U} \right),&\mathcal{E} \left(U\right)= \left(C^{\infty} \left(U,E_{|U} \right),A^{E\prime \prime } _{|U} \right).
\end{align}
We will call $ \mathcal{D} $ the diagonal of $ \mathcal{E} $.
By \eqref{eq:vytt}, $ A^{E\prime \prime }$ is an $ \mathscr{O}_{X} $-morphism, so that $ \mathcal{E} $ is an $ \mathscr{O}_{X} $-complex.
By a fundamental result in \cite[Lemma 4.5]{Block10} and \cite[Theorem 5.2.1]{BismutShenWei23},  $ \mathcal{E} $ has $ \mathscr{O}_{X} $-coherent cohomology.

\begin{exa}\label{ex:fhyj}
  Put
    \begin{align}\label{eq:paqs}
    C_{X} =   \left(\Lambda \left(\overline{T^{*} X} \right), \overline{\partial}^{X} \right).
    \end{align}
Then, $ C_{X} $  is an antiholomorphic superconnection on $ X$.
The associated $ \mathscr{O}_{X}$-complex is the sheaf of Dolbeault complex, which will be denoted by $ \mathscr{D}_{X}$.
\end{exa}


\begin{exa}\label{exa:k2ge}
  Let $ D$ be a bounded holomorphic\footnote{This means each degree $ D^{\bullet} $ as well as the differential $ v_{0} : D^{\bullet } \to D^{\bullet +1} $ are holomorphic.} complex of vector bundles on $ X$.
  If  $ v_{0 } $ and $ \nabla^{D\prime\prime }  $ are respectively the differential and holomorphic structure on $ D$, then
  \begin{align}\label{eq:alyf}
   \left(\Lambda \left(\overline{T^{*} X} \right)\widehat{\otimes}  D, v_{0 } + \nabla^{D\prime\prime }\right)
  \end{align}
 is an antiholomorphic superconnection on $ X$. 
\end{exa}
  
Let $ \left(\underline E,A^{\underline E\prime\prime } \right)$ be another antiholomorphic superconnection with diagonal $ \left(\underline D,\underline v_{0 } \right)$.
As in \cite[Section 4.1]{BismutShenWei23}, if $ p\in \mathbf{Z} $, let  $ \Hom_{}^{p}  \left(E,\underline E\right)$ be the set of $ \Lambda \left(\overline {T^{*} X}\right)$-linear\footnote{The linearity is in graded sense \cite[(4.1.9)]{BismutShenWei23}.
 } maps from $ E$ to $ \underline E$ of degree $ p$.
If $ \phi $ is a section of $ \Hom^{p}  \left(E,\underline E\right)$, put 
\begin{align}\label{eq:vgot}
  A^{\Hom\left(E,\underline E\right)\prime\prime }\phi= A^{\underline E\prime\prime }\phi -\left(-1\right)^{p} \phi  A^{ E\prime\prime }.
\end{align}
Then, $ \left(\Hom \left(E,\underline E\right),A^{\Hom\left(E,\underline E\right)\prime\prime } \right)$ forms an antiholomorphic superconnection.

A morphism of antiholomorphic superconnections \cite[Section 6.2]{BismutShenWei23}, 
\begin{align}\label{eq:23cb}
 \phi: \left(E,A^{E\prime\prime } \right) \to \left(\underline E,A^{\underline E\prime\prime } \right)
\end{align}
is a smooth section of $ \Hom^{0}  \left(E,\underline E\right)$ on $ X$, which is $ A^{\Hom\left(E,\underline E\right)\prime\prime } $-closed, or equivalently, 
\begin{align}\label{eq:2xri}
  A^{\underline E \prime\prime } \phi = \phi A^{E \prime\prime } .
\end{align}

By \eqref{eq:qxmz} and \eqref{eq:2xri}, $ \phi$ induces a morphism of smooth complexes of vector bundles, 
\begin{align}\label{eq:jwot}
 \phi_{0 } :  \left(D,v_{0 } \right) \to \left(\underline D,\underline v_{0 } \right).
\end{align}
As in \eqref{eq:ktja}, using \eqref{eq:aayz}, we can write non-canonically, 
\begin{align}\label{eq:bjwn}
 \phi= \phi_{0 } + \phi_{1} + \ldots + \phi_{n},
\end{align}
where for $ i \ge 1$, $ \phi_{i}\in \Omega^{0,i} \left(X,\Hom^{-i} \left(D,\underline{D} \right)\right)$.

\begin{defin}\label{def:pkpt}
 Let $ {\rm B^{b} } \left(X\right)$ be the category of antiholomorphic superconnections and morphisms on $ X$.\footnote{The subscript $ {\rm b}$ insists that $ E$ is bounded on $ X$.}
\end{defin}

Let us introduce a (non full) subcategory of $ {\rm B^{b}}\left(X\right)$, which is the model at infinity for objects we will introduce later.

If $ \nabla^{D} $ is a grading preserving flat connection on $ D$ with trivial holonomy and if $ v_{0 } $ is flat with respect to $ \nabla^{D} $, then $ \left(D,v_{0 },\nabla^{D} \right)$ will be called a constant complex.
In this case, using the parallel transport with respect to $ \nabla^{D} $,  we see that each degree $ D^{\bullet }\simeq  \mathbf{C}^{{\rm rk}D^{\bullet } } $ is trivial and $ v_{0 } $ is a constant matrix.

If $ \left(D,v_{0 }, \nabla^{D} \right)$ is constant, if $ \nabla^{D\prime \prime }$ is the antiholomorphic part of $ \nabla^{D} $, then 
the antiholomorphic superconnection defined as in \eqref{eq:alyf} is called constant.
In the sequel, we will often omit $ \nabla^{D} $ and say $ \left(D,v_{0 } \right)$ or $\left(\Lambda  \left(\overline{T^{*} X} \right)\widehat{\otimes}  D, v_{0 } +\nabla^{D\prime\prime }\right)$ is constant.

If $ \left(E,A^{E \prime\prime } \right)$ is constant with constant diagonal $ \left(D,v_{0 } \right)$, a section of $ E$ is called constant if it takes values in $ D$ and is flat with respect to $ \nabla^{D} $.
If $ \left(\underline E,A^{\underline E \prime\prime } \right)$ is another constant antiholomorphic superconnection, then $ \left(\Hom\left(E,\underline{E} \right),A^{{\rm Hom}\left(E,\underline{E} \right)\prime \prime } \right)$ is still constant.
A morphism  $\phi :  \left(E,A^{E\prime\prime } \right) \to  \left(\underline E,A^{\underline E\prime\prime } \right)$ is called  constant if $ \phi$ is a constant section of $ {\rm Hom}\left(E,\underline E\right)$.

Let $ Z \subset X$ be a subset of $ X$. 
Two antiholomorphic superconnections defined on some neighbourhoods of $ Z$ are called equivalent if they restrict to the same antiholomorphic superconnection on some smaller open neighbourhood of $ Z$.
The associated equivalence class 
will be called an antiholomorphic superconnection on $ Z$.
Its support is now a closed subset in $ Z$.

A section of an antiholomorphic superconnection on $ Z$,  a morphism between two antiholomorphic superconnections on $ Z$, as well as the similar notations for constant ones are defined in an obvious way.






\subsection{Tame antiholomorphic superconnection}\label{s:BtX}
In the rest part of this section, we assume that  $ X  $ is an open submanifold  of some ambient manifold $ \widetilde{X} $.
Let $ \overline X$ be the closure of $ X$ in $ \widetilde{X} $.
The objects introduced later depend only on the germ of the closed subset $ \overline X$.

\begin{defin}\label{def:ubgi}
An antiholomorphic superconnection $ \left( E,A^{ E \prime\prime } \right)$ on $ \overline{X} $ is called tame, if it is constant outside a compact subset of $ \overline{X} $.


A morphism of two tame antiholomorphic superconnections 
\begin{align}\label{eq:pbhl}
  \phi :\left( E,A^{ E \prime\prime } \right) \to \left(\underline E,A^{\underline E \prime\prime } \right)
 \end{align}
on $ \overline{X} $ is called tame, if $ \phi $ is constant outside a compact subset of $ \overline{X} $.

Let $ {\rm B^{b,t}}  \left(\overline X\right)$ be the category of tame antiholomorphic superconnections and  tame morphisms on $ \overline X$.

  Let $ {\rm B_{c}^{b,t} }\left(\overline X\right)$ be the full subcategory of $ {\rm B^{b,t} }\left(\overline X\right)$ consisting of objects with compact support in $ \overline X$. 
\end{defin}

\begin{exa}\label{ex:dfqz}
 The antiholomorphic superconnection $ C_ {\widetilde{X} } $ on $ \widetilde{X} $ defined in Example \ref{ex:fhyj}  induces a tame antiholomorphic superconnection $ C_ {\overline{X} } $ on $ \overline{X} $.
\end{exa}

An open cover $ \mathcal{U} $ of $ \overline{X} $ is called adapted if except one $ U_{0 }\in \mathcal{U} $, all the other open sets in $ \mathcal{U} $ are relatively compact in $ \overline{X} $.

\begin{thm}\label{thm:jltn}
  If $ \left(E,A^{E\prime\prime } \right)$ is an object of $ {\rm B^ {b,t}} \left(\overline{X} \right)$, then there exists an adapted finite open cover  $\mathcal{U} $ of $ \overline{X}  $  such that 
 \begin{enumerate}[\indent a)]
     \item on $ U_{0 } $, $ \left(E,A^{E\prime\prime } \right)_{|U_{0 } }$ is constant;
   \item on the other $ U\in \mathcal{U} \backslash \{U_{0 } \}$, $ \left(E,A^{E\prime\prime } \right)_{|U}  $ is (non-canonically) isomorphic to an antiholomorphic superconnection associated to a bounded holomorphic complex of vector bundles on $ U$ (see Example \ref{exa:k2ge}). 
 \end{enumerate}
 If $ \left(E,A^{E\prime\prime } \right)$ has compact support in $ \overline{X} $, we can require further that 
 \begin{enumerate}[\indent c)]
   \item the constant diagonal $ \left(D,v_{0 } \right)_{|U_{0 } } $ is exact.
 \end{enumerate}
\end{thm}
\begin{proof}
  By \cite[Theorem 5.2.1]{BismutShenWei23}, for any $ x\in \overline X$, there exists a small open neighbourhood $ \widetilde{U} _{x}  $ of $ x$ in the ambient manifold $ \widetilde{X} $ such that  $ \left(E,A^{E\prime\prime } \right)_{|{\widetilde{U}} _{x} }$ is well-defined and isomorphic to the one defined by a bound holomorphic complex of vector bundles on $ \widetilde{U}_{x}$.
  Our theorem follows easily from this consideration and the tameness assumption.
\end{proof}

If $ \left( E,A^{ E \prime\prime } \right) $  and $  \left(\underline E,A^{\underline E \prime\prime } \right)$ are objects in $ {\rm B^{b,t}} \left(\overline{X} \right)$, 
then $ \left({\rm Hom}\left(E,\underline{E} \right),A^{{\rm Hom}\left(E,\underline{E} \right)\prime \prime } \right)$  is still an object in $ {\rm B^{b,t}} \left(\overline{X} \right)$.

If $ \phi : \left( E,A^{ E \prime\prime } \right) \to \left(\underline E,A^{\underline E \prime\prime } \right)$ is a  morphism in $ {\rm B^{b,t}}\left(\overline{X} \right) $, as in \cite[Section 5.6]{BismutShenWei23}, put
\begin{align}\label{eq:yrmd}
  &{\rm cone}^{\bullet} \left(\phi\right)=E^{\bullet+1} \oplus \underline E^{\bullet} ,&A^{{\rm cone} \left(\phi\right)\prime\prime } = \begin{pmatrix}
    A^{E\prime\prime }  & 0 \\
    \phi\left(-1\right)^{{\rm deg}}  & A^{\underline E \prime\prime }
  \end{pmatrix}.
\end{align}
Then, $  \left({\rm cone} \left(\phi\right),A^{{\rm cone} \left(\phi\right)\prime\prime } \right)$ is still an object in $ {\rm B^{b,t}}\left(\overline{X} \right) $.
If $ \left(E,A^{E^{\prime\prime } } \right)$ and $ \left(\underline{E} ,A^{\underline{E} ^{\prime\prime } } \right)$ have compact support, then $  \left({\rm cone} \left(\phi\right),A^{{\rm cone} \left(\phi\right)\prime\prime } \right)$ has also compact support.

A tame morphism $ \phi :\left( E,A^{ E \prime\prime } \right) \to  \left(\underline E,A^{\underline E \prime\prime } \right)$ is called tame null-homotopic if there is a tame section $ h$  of $ {\rm Hom} \left(E,\underline{E} \right)$ of degree $ -1$ on $ \overline{X} $ such that 
\begin{align}\label{eq:fxpd}
 \phi = A^{{\rm Hom}\left(E,\underline{E} \right)\prime \prime }h,
\end{align}
or equivalently,
\begin{align}\label{eq:qee2}
 \phi = h A^{E\prime\prime } + A^{\underline E\prime\prime } h.
\end{align}
Two tame morphisms $ \phi_{1}, \phi_{2} :\left( E,A^{ E \prime\prime } \right) \to  \left(\underline E,A^{\underline E \prime\prime } \right)$ are called tame homotopic if $ \phi_{1}-\phi_{2}$ is a tame null-homotopy.

\begin{defin}\label{def:q2gd}
 Let $ \underline {\rm B}_{}^{{\rm b,t}} \left(\overline{X} \right)$ be the homotopy category of $ {\rm B^{b,t}} \left(\overline{X} \right)$. 
 Its objects coincide with the objects of $ {\rm B^{b,t}} \left(\overline{X} \right)$, and the morphisms are given by the tame homotopy equivalence classes of morphisms of $ {\rm B^{b,t}} \left(\overline{X} \right)$.

 Let $ \underline{{\rm B}}^{{\rm b,t}}_{{\rm c}}  \left(\overline{X} \right)$ be the full subcategory of $ \underline {\rm B}_{}^{{\rm b,t}} \left(\overline{X} \right)$ consisting of objects with compact support in $ \overline{X} $.
\end{defin}

A triangle in $ \underline {\rm B}^{{\rm b,t}}  \left(\overline{X} \right)$ or $ \underline{{\rm B}}_{{\rm c}} ^{{\rm b,t}} \left(\overline{X} \right) $  is called distinguished if it is isomorphic to a triangle associated to some morphism $ \phi :\left( E,A^{ E \prime\prime } \right) \to  \left(\underline E,A^{\underline E \prime\prime } \right)$, 
\begin{equation}\label{dia:b3rs}
\begin{tikzsd}
    E  &	\underline{E}  & {\rm cone}\left(\phi \right) & E^{\bullet +1}.
     \arrow["", from=1-1, to=1-2]
     \arrow["", from=1-2, to=1-3]
     \arrow["", from=1-3, to=1-4]
  \end{tikzsd} 
\end{equation}

\begin{prop}\label{prop:xiqu}
  The categories $ \underline {\rm B}^{{\rm b,t}}  \left(\overline{X} \right)$ and $ \underline{{\rm B}}_{{\rm c}} ^{{\rm b,t}} \left(\overline{X} \right) $  together with their distinguished triangles form triangulated categories.
\end{prop}
\begin{proof}
 We need to verify (TR0)-(TR5) of \cite[Proposition 1.4.4]{KashiwaraSchpira}.
 This is immediate since in our context all the morphisms or sections constructed in the proof of \cite[Lemma 1.4.2, Proposition 1.4.4]{KashiwaraSchpira} are tame.
\end{proof}

\subsection{The tame $ \mathscr{O} $-complexes on $ \overline X$}\label{s:Cbt}
We call interchangeably an $ \mathscr{O}_{X} $-module or an $ \mathscr{O} $-module on $ X$.

If $ Z $ is a subset of $ X$, we say two $ \mathscr{O} $-modules defined on some open neighbourhoods of $ Z$ are equivalent if they restrict to the same $ \mathscr{O} $-module on a smaller open neighbourhood of $ Z$.
The associated equivalence class will be called an $ \mathscr{O} $-module on $ {Z} $.
We define morphisms of $ \mathscr{O} $-modules on $Z $ in an obvious way.
This gives an Abelian category $ {\rm M}\left(Z\right)$ of $ \mathscr{O} $-modules on  $ {Z} $.

An $ \mathscr{O} $-module on $ Z$ is called coherent if it has a representative which is coherent on some open neighbourhood of $ Z$.
We define the category of coherent $ \mathscr{O} $-modules on $ Z$ as a full subcategory of $ {\rm M}\left(Z\right)$. 
Clearly, this is still an Abelian category.

An $ \mathscr{O} $-complex on $ Z$ is defined in a similar way.
A holomorphic complex of vector bundles on $ Z$ induces naturally an $ \mathscr{O}$-complex on $ Z$.
 
An $ \mathscr{O} $-complex on $ Z$ is called constant, if it is induced by a constant complex of vector bundles on $ Z$.
We define the constant morphisms in the same way.

\begin{defin}\label{def:32oh}
 An $ \mathscr{O} $-complex $ \mathcal{F} $ on $ \overline{X} $ is called tame if it is constant outside a compact subset of $ \overline{X} $.

 An $ \mathscr{O} $-morphism $ \phi : \mathcal{F} \to \mathcal{F} '$ of two tame $ \mathscr{O} $-complexes on $ \overline{X} $ is called tame, if $ \phi$ is constant outside a compact subset of $ \overline X$.
 
 Let  $ {\rm C^{b,t}_{coh}} \left(\overline X\right)$ be the category of  tame $ \mathscr{O} $-complexes with bounded coherent cohomology and tame $ \mathscr{O} $-morphisms on $ \overline X$.
\end{defin}



\begin{exa}\label{exa:nqs3}
  An $ \mathscr{O} $-complex on $ \overline{X} $ with compact support is tame.
  Morphisms of two such $ \mathscr{O} $-complexes on $ \overline{X} $ are also tame.
\end{exa}

\begin{exa}\label{exa:nqs4}
  If $ \mathcal{F}  $ is an object in $ {\rm C^{b,t}_{coh}} \left(\overline{X} \right)$ with differential $ v^{\mathcal{F} } $, if $ p\in \mathbf{Z} $, then the canonical truncations \cite[(1.3.10), (1.3.11)]{KashiwaraSchpira}, 
  \begin{equation}\label{dia:vtj1}
    \begin{aligned}
    \tau_{\ge p} \mathcal{F} &:   
\begin{tikzsd}
      0 & \mathcal{F}^{p}/{\rm im}\left(v^{\mathcal{F} } \right)\cap \mathcal{F}^{p}  &	\mathcal{F}^{p+1}  & \mathcal{F}^{p+2}   & \cdots,
       \arrow["", from=1-1, to=1-2]
       \arrow["", from=1-2, to=1-3]
       \arrow["", from=1-3, to=1-4]
       \arrow["", from=1-4, to=1-5]
    \end{tikzsd} \\
    \tau_{\le p} \mathcal{F} &:
\begin{tikzsd}
    \cdots & \mathcal{F}^{p-2}  &	\mathcal{F}^{p-1}  &  \ker v^{\mathcal{F} } \cap \mathcal{F}^{p} & 0,
     \arrow["", from=1-1, to=1-2]
     \arrow["", from=1-2, to=1-3]
     \arrow["", from=1-3, to=1-4]
     \arrow["", from=1-4, to=1-5]
  \end{tikzsd} 
  \end{aligned}
  \end{equation}
are still in $ {\rm C^{b,t}_{coh}} \left(\overline{X} \right)$.
We have canonical morphisms in $ {\rm C^{b,t}_{coh}} \left(\overline{X} \right)$, 
\begin{align}\label{eq:zh1x}
&\mathcal{F} \to \tau_{\ge p} \mathcal{F}, & \tau_{\le p}  \mathcal{F} \to \mathcal{F}. 
\end{align}
\end{exa}

\begin{thm}\label{thm:yndp}
  If $ \mathcal{F} $ is an object in $ {\rm C^{b,t}_{coh}}\left(\overline{X} \right)$, then there exist an adapted finite open cover $\mathcal{U} $  of $ \overline{X} $ and $ \ell\in \mathbf{N} $, such that
  \begin{enumerate}[\indent a)]
    \item on $ U_{0 } \in \mathcal{U} $, $ \mathcal{F}_{|U_{0 } } $ is constant;
    \item on the other $ U\in \mathcal{U} \backslash \{U_{0 } \}$, for each $ p\in \mathbf{Z} $, there exist a holomorphic complex of trivial vector bundles $ R_{U} $ on $ U $ concentrated in degrees $ \left[-\ell ,0\right]$, an $ \mathscr{O} $-morphism $ \alpha_{U} :\mathscr{O} \left(R^{0 }_{U}  \right)\to \left(\mathcal{H}^{p} \mathcal{F}\right)_{|U} $ on $ U$ such that we have an exact $ \mathscr{O} $-complex on $ U$,  
  \begin{equation}\label{dia:a1eh}
\begin{tikzsd}
      0 & 	 \mathscr{O} \left({R} ^{-\ell }_{U} \right) & \cdots  & \mathscr{O}  \left({R}  ^{0 }_{U} \right)  & \left(\mathcal{H}^{p} \mathcal{F}\right)_{|U } &0 .
       \arrow["", from=1-1, to=1-2]
       \arrow["", from=1-2, to=1-3]
       \arrow["", from=1-3, to=1-4]
       \arrow["\alpha_{U} ", from=1-4, to=1-5]
      \arrow["", from=1-5, to=1-6]
    \end{tikzsd} 
  \end{equation}
\end{enumerate}
\end{thm}
\begin{proof}
 Our theorem follows immediately from the definition of coherence.
 Thanks to our tameness condition on $ \mathcal{F} $ and the boundedness condition on $ \mathcal{H} \mathcal{F} $, $ \ell$ can be chosen to be independent of $ p$ and $ U\in \mathcal{U} \backslash \left\{U_{0 } \right\}$.
\end{proof}

If $ \phi :  \mathcal{F} \to \mathcal{F} '$ is a morphism in $ {\rm C^{b,t}_{coh}}\left(\overline{X} \right) $, then $ {\rm cone}\left(\phi \right)$ is still an object in $ {\rm C^{b,t}_{coh}}\left(\overline{X} \right) $.  
If $ \mathcal{F}, \mathcal{F} '$ have compact support in $ \overline{X} $, then $ {\rm cone}\left(\phi \right)$ also has compact support in $ \overline{X}$.

Let $ \mathscr{O}^{\infty} $ be the sheaf of smooth functions, and let $ \mathscr{D}  $ be the sheaf of Dolbeault complex (see Example \ref{ex:fhyj}).
Put 
\begin{align}\label{eq:d2vy}
& \mathscr{D}^{i}  \mathcal{F} = \mathscr{D}^{i}  \otimes_{\mathscr{O} } \mathcal{F},&\mathscr{D} \mathcal{F} = \mathscr{D} \widehat{\otimes }_{\mathscr{O} } \mathcal{F}.  
\end{align}
If $ \overline{\partial}$ and $ v^{\mathcal{F} }$ are respectively the differentials of $\mathscr{D}$ and $ \mathcal{F}$, the differential on $ \mathscr{D} \mathcal{F} $ is given by 
\begin{align}\label{eq:ute1}
 A^{\mathscr{D} \mathcal{F} \prime \prime } = v^{\mathcal{F} }+\overline{\partial}.
\end{align}
Then, $ \mathscr{D} \mathcal{F} $ is an $ \mathscr{O} $-complex on $ \overline{X} $, which however is not tame.

By the holomorphic Poincaré lemma and the fact that $\mathscr{O}^{\infty},  \mathscr{D} $ are flat over $ \mathscr{O} $ \cite{Malgrange67}, the canonical morphism of $ \mathscr{O} $-complexes on $ \overline{X} $, 
\begin{align}\label{eq:gvcd}
  \mathcal{F} \to \mathscr{D} \mathcal{F}
 \end{align}
is a quasi-isomorphism.
Moreover, 
\begin{align}\label{eq:1qmc}
\mathscr{O}^{\infty} \otimes_{\mathscr{O} } \left(\mathcal{H} \mathcal{F} \right) = \mathcal{H} \left(\mathscr{D}^{0 }  \mathcal{F}\right).
\end{align}
Since $ \mathscr{O}^{\infty} \otimes_{\mathscr{O} }  \cdot $ is faithful,\footnote{\label{Footnote:FF}If $ A$ is an $ \mathscr{O} $-module, using the above argument, we have an exact sequence $ 0 \to A \to \mathscr{O}^{\infty} \otimes_{\mathscr{O} } A \to \mathscr{D}^{1} A\to \cdots $. Therefore, $ \mathscr{O}^{\infty} \otimes_{\mathscr{O} } A= 0$ implies $ A=0$.
This shows $ \mathscr{O}^{\infty} $ is faithfully flat over $ \mathscr{O} $ \cite[\href{https://stacks.math.columbia.edu/tag/00HP}{Tag 00HP}]{stacks-project}.
} 
by \eqref{eq:1qmc}, we have 
\begin{align}\label{eq:3sjd}
 \mathcal{H} \mathcal{F}= 0 \iff \mathcal{H} \left(\mathscr{D}^{0 }  \mathcal{F}\right)= 0. 
\end{align}

\subsection{The $ \mathscr{D} $-tame $ \mathscr{O} $-complexes on $ \overline{X} $} \label{s:Dbt}
The $ \mathscr{O} $-complex $ \mathcal{E} $ associated to an object of $ {\rm B^{b,t}} \left(\overline{X} \right)$ and the $ \mathscr{O} $-complex $ \mathscr{D} \mathcal{F} $ in \eqref{eq:d2vy} associated to an object in $ {\rm C^{b,t}_{coh}} \left(\overline{X} \right)$
are not tame.
Let us introduce another tameness condition.


\begin{defin}\label{def:vojk}
  An $ \mathscr{O} $-complex on $ \overline{X} $ is called $ \mathscr{D} $-tame, if it is also a $ \mathbf{Z} $-graded $ \mathscr{D} $-module, its differential satisfies the antiholomorphic Leibniz rules (cf. \eqref{eq:vytt}), and it is defined by a constant antiholomorphic superconnection outside a compact subset of $ \overline{X}$.

  A morphism of two above $ \mathscr{O} $-complexes is called $ \mathscr{D} $-tame, if it is also a morphism of $ \mathscr{D} $-modules, and if it is defined by a constant morphism of constant antiholomorphic superconnections outside a compact subset of $ \overline{X} $.

  Two above $ \mathscr{D} $-tame morphisms are called $ \mathscr{D} $-tame homotopy equivalent if their difference is $ \mathscr{D} $-tame null-homotopic defined in an obvious way.

  Let $ \widetilde{{\rm B}}^{{\rm t}} \left(\overline{X} \right)$ be the category of $ \mathscr{D} $-tame $ \mathscr{O} $-complexes and $ \mathscr{D} $-tame morphisms on $ \overline{X} $.

  Let $ \widetilde{\underline{{\rm B}} }^{{\rm t}} \left(\overline{X} \right)$ be the homotopy category of $ \widetilde{{\rm B}}^{{\rm t}} \left(\overline{X} \right)$.
  Its objects coincide with the objects of $ \widetilde{{\rm B}}^{{\rm t}} \left(\overline{X} \right)$, and the morphisms are given by the $ \mathscr{D} $-tame homotopy equivalence classes of morphisms of $ \widetilde{{\rm B}}^{{\rm t}} \left(\overline{X} \right)$.
\end{defin}
 
If $ \psi $ is a morphism in $ \widetilde{{{\rm B}} }^{{\rm t}} \left(\overline{X} \right)$, then $ {\rm cone}\left(\psi \right)$ is still an object in $ \widetilde{{{\rm B}} }^{{\rm t}} \left(\overline{X} \right)$.
We define the distinguished triangle in $ \widetilde{\underline{{\rm B}} }^{{\rm t}} \left(\overline{X} \right)$ in an obvious way.

\begin{prop}\label{prop:vkig}
 The category $ \widetilde{\underline{{\rm B}} }^{{\rm t}} \left(\overline{X} \right)$ together with its distinguished triangles forms a triangulated category.
\end{prop}
\begin{proof}
 The proof is the same as the one of Proposition \ref{prop:xiqu}.
\end{proof}


The assignment $ \left(E,A^{E\prime \prime } \right) \to \mathcal{E} $ 
induces fully faithful functors $ {\rm B^{b,t}}\left(\overline{X} \right) \to \widetilde{{\rm B}}^{{\rm t}} \left(\overline{X} \right)$ and $ {\rm \underline{B} ^{b,t}}\left(\overline{X} \right) \to \widetilde{{\rm \underline{B} }}^{{\rm t}} \left(\overline{X} \right)$.
Moreover, $ \mathcal{F} \to \mathscr{D} \mathcal{F} $  induces a functor $ {\rm C^{b,t}_{coh}}\left(\overline{X} \right)  \to \widetilde{{\rm B}}^{{\rm t}} \left(\overline{X} \right)$.

\begin{prop}\label{prop:ceuw}
If $ \left(E,A^{{E}\prime \prime  } \right)$ is an object in $ {\rm B^{b,t}}\left(\overline{X} \right)$, if $ \mathcal{F} $ is an object in $ {\rm C^{b,t}_{coh}}\left(\overline X\right)$, and if  $ \phi :  \mathcal{E} \to  \mathscr{D} \mathcal{F}$ is a morphism in $ \widetilde{{\rm B}}^{{\rm t}} \left(\overline{X} \right)$, then $ \phi$ is a quasi-isomorphism if and only if it induces a quasi-isomorphism on the diagonals, 
\begin{align}\label{eq:bunb}
 \phi_{0 } : \mathcal{D} \to \mathscr{D}^{0} \mathcal{F}. 
\end{align}
\end{prop}
\begin{proof}
 Since our statement is local, by \cite[Theorem 5.2.1]{BismutShenWei23} or Theorem \ref{thm:jltn}, we can assume that $ \left(E,A^{E\prime \prime } \right)$ is defined by a bounded holomorphic complex of vector bundles $ \left(D,v_{0 } \right)$.
 Let $\left( \mathcal{E}_{r}\right)_{r\ge 0 }  $ and $ \left(\underline{\mathcal{E}} _{r}\right)_{r\ge 0 } $ be the spectral sequences of $ \mathcal{E} $ and $ \mathscr{D} \mathcal{F} $ \cite[Definition 5.3.1]{BismutShenWei23} associated to the filtration induced by \eqref{eq:uucw}.

 By \cite[(5.3.1), Theorem 5.3.4]{BismutShenWei23} and by \eqref{eq:1qmc}, we know that $\left( \mathcal{E}_{r}\right)_{r\ge 0 }  $ and $ \left(\underline{\mathcal{E}} _{r}\right)_{r\ge 0 } $ degenerate at $ r= 2$, and\footnote{Here, $ \mathscr{O}^{\infty}  \left(D\right)$ or $ \mathscr{O} \left(D\right)$ are respectively the sheaves of smooth or holomorphic sections of $ D $.}
 \begin{align}\label{eq:laec}
  &\mathcal{E}_{0} = \left(\mathscr{D} \otimes_{\mathscr{O}^{\infty}  }   \mathscr{O}^{\infty} \left(D\right),v_{0 } \right),&\mathcal{E}_{1} = \left(\mathscr{D} \otimes_{\mathscr{O}}  \mathcal{H} \mathscr{O} \left(D\right), \overline{\partial}  \right),&&\mathcal{E}_{2} = \mathcal{H} \mathscr{O} \left(D\right),\\
  &\underline{\mathcal{E}} _{0} = \left(\mathscr{D} \otimes_{\mathscr{O}^{\infty}  }  \mathcal{D}^{0 } \mathcal{F},v^{\mathcal{F} } \right),&\underline{\mathcal{E}} _{1} = \left(\mathscr{D} \otimes_{\mathscr{O}}  \mathcal{H} \mathcal{F}, \overline{\partial}  \right),&&\underline{\mathcal{E}} _{2} = \mathcal{H} \mathcal{F}.\notag
 \end{align}


For each $ r \ge 0 $, $ \phi$ induces morphisms of complexes $ \mathcal{E}_{r} \to \underline{\mathcal{E}}_{r} $.
When $ r= 0 $, this morphism is just $ \phi_{0} $.
The degree $ 1$ part of the equation $ \left(v^{\mathcal{F} } +\overline{\partial} \right)\phi = \phi\left(v_{0 } + \nabla^{D\prime \prime } \right) $ gives 
\begin{align}\label{eq:q11a}
 \overline{\partial} \phi_{0 } + v^{\mathcal{F} } \phi_{1} = \phi_{1} v_{0 } + \phi_{0 } \nabla^{D\prime \prime }. 
\end{align}
By \eqref{eq:q11a}, $ \phi_{0} $ induces a morphism 
\begin{align}\label{eq:bpqh}
 \mathcal{H} \phi_{0} : \mathcal{H} \mathscr{O} \left(D\right)\to \mathcal{H} \mathcal{F}, 
\end{align}
which gives the morphisms $ \mathcal{E}_{1} \to \underline{\mathcal{E} }_{1} $ and $ \mathcal{E}_{2} \to \underline{\mathcal{E} }_{2} $.
By \cite[Theorem 5.3.4]{BismutShenWei23}, the filtration on the cohomologies $ \mathcal{H} \mathcal{E} $ and $ \mathcal{H} \mathscr{D} \mathcal{F} $ are trivial, so the morphism $ \mathcal{E}_{2} \to \underline{\mathcal{E} }_{2} $ coincides with the natural morphism $ \mathcal{H} \mathcal{E} \to \mathcal{H} \mathscr{D} \mathcal{F} $ induced by $ \phi $.

By above, $ \phi$ is a quasi-isomorphism if and only if \eqref{eq:bpqh} is an isomorphism, and $ \phi_{0 } $ is a quasi-isomorphism if and only if \eqref{eq:bpqh} induces an isomorphism $ \mathscr{O}^{\infty} \otimes_{\mathscr{O} } \mathcal{H} \mathscr{O} \left(D\right)\to \mathscr{O}^{\infty} \otimes_{\mathscr{O} } \mathcal{H} \mathcal{F}$.
Our proposition follows now from the fact that $ \mathscr{O}^{\infty} $ is faithfully flat over $ \mathscr{O} $ (see Footnote \ref{Footnote:FF}).
\end{proof}

The main result of this section is the following.

\begin{thm}\label{thm:poev}
  The following statements hold.
  \begin{enumerate}[\indent a)]
    \item\label{thm:poeva}  If $ \mathcal{F} $ is an object in $ {\rm C^{b,t}_{coh}}\left(\overline X\right)$, then there exist an object $ \left(E,A^{{E}\prime \prime  } \right)$ in $ {\rm B^{b,t}}\left(\overline{X} \right)$ with the associated $ \mathscr{O} $-complex $ \mathcal{E} $ and a quasi-isomorphism in $ \widetilde{{\rm B}}^{{\rm t}} \left(\overline{X} \right)$,  
    \begin{align}\label{eq:fdcq}
   \phi :  \mathcal{E} \to  \mathscr{D} \mathcal{F}.
    \end{align}
    \item\label{thm:poevb} If $ \left(\left(\underline{E} ,A^{\underline{E} \prime \prime } \right),\underline{\mathcal{E}} , \underline{\mathcal{F}} , \underline{\phi} \right)$ is another data as $ \left(\left(E,A^{E\prime \prime } \right),\mathcal{E}, \mathcal{F}, \phi\right)$ in a),  if $ \psi: \mathcal{F} \to \underline{{\mathcal{F} }}  $ is a morphism in $ {\rm C^{b,t}_{coh}} \left(\overline{X} \right)$, there exists a unique morphism $ \widetilde{\psi } $ in $ {\rm \underline{B} ^{b,t}}\left(\overline{X} \right) $ such that the  diagram
    \begin{equation}\label{dia:1qca}
\begin{tikzsd}
          \mathcal{E}  & \underline{{\mathcal{E} }}   \\
          \mathscr{D} \mathcal{F}  & \mathscr{D} \underline{{\mathcal{F} }}
          \arrow["\phi ", from=1-1, to=2-1]
          \arrow["\widetilde{\psi }", from=1-1, to=1-2,dotted]
          \arrow["\underline{\phi} ", from=1-2, to=2-2]
          \arrow["\psi ", from=2-1, to=2-2]
        \end{tikzsd}
    \end{equation}
    in $ \widetilde{\underline{{\rm B}}}^{{\rm t}}  \left(\overline{X} \right)$ commutes.

    
  \end{enumerate}
\end{thm}
\begin{proof}
 The proof of our theorem will be given in Sections \ref{sEssSuj} and \ref{s:thmb}.
\end{proof}

\subsection{Proof of Theorem \ref{thm:poev} \ref{thm:poeva}) }\label{sEssSuj}
The tameness condition for $ \mathscr{O}^{\infty} $-modules is defined in a similar way.

Let $ \mathcal{G} $ be an $ \mathscr{O}^{\infty} $-complex on $ \overline{X} $ with the properties: there exist an adapted finite open cover $ \mathcal{U} $ of $ \overline{X} $  and $ \ell\in \mathbf{N} $ such that  
 \begin{enumerate}[\indent a)]
   \item on $ U_{0 } $, $ \mathcal{G} _{|U_{0 } } $ is defined by a constant complex of vector bundles $ G_{U_{0 } } $, i.e.,  
   \begin{align}\label{eq:in1y}
    \mathcal{G}_{|U_{0 } } =  \mathscr{O}^{\infty} \left(G_{U_{0 } } \right);
   \end{align}
   \item on $ U\in \mathcal{U} \backslash \{U_{0 } \}$, for each $ p\in \mathbf{Z} $, there exist a smooth complex $ R_{U} $ of trivial vector bundles on $ U$ concentrated in degrees $ \left[-\ell, 0 \right]$ and an $ \mathscr{O}^{\infty} $-morphism $ \alpha_{U} : \mathscr{O}^{\infty} \left(R^{0 }_{U}\right) \to \left(\mathcal{H}^{p} \mathcal{G} \right)_{|U} $ such that we have an exact $ \mathscr{O}^{\infty} $-complex on $ U$,
\begin{equation}\label{dia:l3ut}
\begin{tikzsd}
    0 & \mathscr{O}^{\infty} \left(R_{U}^{- \ell } \right) &	\cdots & \mathscr{O}^{\infty} \left(R_{U}^{0}\right)  & \left(\mathcal{H}^{p} \mathcal{G}\right) _{|U} & 0 .
     \arrow["", from=1-1, to=1-2]
     \arrow["", from=1-2, to=1-3]
     \arrow["", from=1-3, to=1-4]
     \arrow["\alpha_{U} ", from=1-4, to=1-5]
     \arrow["", from=1-5, to=1-6]
  \end{tikzsd} 
\end{equation}
\end{enumerate}

Note that for $ p\in \mathbf{Z} $, the canonical truncations $ \tau_{\ge p} \mathcal{G} $ and $ \tau_{\le p} \mathcal{G} $ (see Example \ref{exa:nqs4}) are still $ \mathscr{O}^{\infty} $-complexes with the same properties. 

If $ U\in \mathcal{U} $, there exists an open subset $ \widetilde{U} $ of $ \widetilde{X} $ such that $ U= \widetilde{U} \cap \overline{X} $.
Let $ \left( \phi_{U}\right)_{U\in \mathcal{U} } $ be a smooth partition of unity of $ \overline{X} $ subordinate to $ \mathcal{U} $.
This means that $\phi_{U}$ is a restriction of $ \phi_{\widetilde{U} }\in C^{\infty} \left(\widetilde{X}, [0, 1]\right)$\footnote{Observe that $ \left\{\widetilde{U} \right\}_{U\in \mathcal{U} } \cup \left\{\overline{X}^{c} \right\}$ forms an open cover for  $ \widetilde{X} $.
A partition of unity subordinate to this open cover in classical sense gives our family $ \left\{\phi_{\widetilde{U}} \right\}_{U\in \mathcal{U} } $.} such that 
\begin{align}\label{eq:klag}
 {\rm Supp}\left(\phi_{\widetilde{U} } \right)
 \subset \widetilde{U}.
\end{align}
Moreover, on some open neighbourhood of $ \overline{X} $, we have 
\begin{align}\label{eq:ityw}
 \sum_{U\in \mathcal{U} } \phi_{\widetilde{U} } = 1.
\end{align}

\begin{prop}\label{prop:cfiv}
  If $ p,q\in \mathbf{Z} $ with $ p\le q$ and if $ \mathcal{G} $ is as above such that the cohomology $ \mathcal{H} \mathcal{G}  $ is concentrated in degrees $ [p,q]$, then 
  there exist a tame smooth complex of vector bundles $ \left(D,v_{0 } \right)$ on $ \overline{X} $ concentrated in degrees $  [p-\ell,q]$, and a tame $ \mathscr{O}^{\infty}$-quasi-isomorphism on $ \overline X$,   
   \begin{align}\label{eq:2kcx}
    \mathscr{O}^{\infty} \left(D\right)\to \mathcal{G} .
   \end{align}
  \end{prop}
  \begin{proof}
  We follow the proof of \cite[Proposition 6.3.2]{BismutShenWei23}.
  Up to a shift, we may and we will assume that $ p= 0 $, i.e., $ \mathcal{H} \mathcal{G} $ is concentrated in degrees $ [0, q]$.
  
  \underline{Step 1 : a reduction.}
  We claim that the proof of our proposition can be reduced to the case where $ \mathcal{G} $ is concentrated in positive degrees, i.e., for all $ k<0 $, 
  \begin{align}\label{eq:wf1p}
    \mathcal{G}^{k} = 0.
  \end{align}
  Indeed, the truncation $ \tau_{\ge 0 } \mathcal{G} $ has the same property as $ \mathcal{G} $ and is concentrated in the positive degrees.
  Moreover, the canonical tame $ \mathscr{O}^{\infty} $-morphism on $ \overline{X} $, 
  \begin{align}\label{eq:gnfk}
  \pi : \mathcal{G} \to \tau_{\ge 0 } \mathcal{G}
  \end{align}
  is a quasi-isomorphism.
  If our proposition holds for $ \tau_{\ge 0 } \mathcal{G} $, there exist a tame smooth complex of vector bundles $ \left(D,v_{0 }\right)$ on $ \overline{X} $ concentrated in degrees $ [-\ell, q]$ and a tame $ \mathscr{O}^{\infty} $-quasi-isomorphism on $ \overline{X} $, 
  \begin{align}\label{eq:ornt}
   \phi : \mathscr{O}^{\infty} \left(D\right) \to \tau_{\ge 0 } \mathcal{G}. 
  \end{align}
  By \cite[\href{https://stacks.math.columbia.edu/tag/0649}{Tag 0649}]{stacks-project}, there exists a morphism of $ \mathscr{O}^{\infty} $-complexes $ \phi ^{\prime } $ such that the diagram 
  \begin{equation}\label{dia:zbro}
\begin{tikzsd}
         & \mathcal{G} \\
        \mathscr{O}^{\infty} \left({D} \right) & \tau_{\ge 0 } \mathcal{G}
        \arrow["\phi ^{\prime } ", from=2-1, to=1-2,dotted]
        \arrow["\pi ", from=1-2, to=2-2]
        \arrow["\phi", from=2-1, to=2-2]
      \end{tikzsd}
  \end{equation}
  is commute up to homotopy.  
  In general,  $ \phi ^{\prime } $ is not necessarily tame.

  Since $ D,\tau_{\ge 0 } \mathcal{G}, \mathcal{G} $ are constant on $ U_{0 } $, 
  applying \cite[\href{https://stacks.math.columbia.edu/tag/0649}{Tag 0649}]{stacks-project} again, on $ U_{0 } $ there exist a constant lift such that \eqref{dia:zbro} commutes up to constant homotopy.
  If $ v^{\mathcal{G} } $ denotes the differential of $ \mathcal{G} $, by the uniqueness of the lift \cite[\href{https://stacks.math.columbia.edu/tag/064A}{Tag 064A}]{stacks-project},  there is $ h : \mathscr{O}^{\infty} \left(D_{|U_{0 } } \right) \to \mathcal{G}_{|U_{0 } }  $ of degree $ -1$  such that 
  \begin{align}\label{eq:ckgx}
   \phi^{\prime }   + v^{\mathcal{G} } h + h v_{0 } 
  \end{align}
  is well-defined and constant on $ U_{0 } $.

  By \eqref{eq:klag}, $h_{U_{0 } } =  \phi_{U_{0 } } h $ is well-defined\footnote{\label{foot:Us}We can take $ \widetilde{U}_{0} $ sufficiently small so that $ h$ extends to $ \widetilde{U}_{0} $.} on $ \overline{X} $.    
  Put 
  \begin{align}\label{eq:bbwx}
   \phi ^{\prime \prime } = \phi ^{\prime }  + v^{\mathcal{G} }  h_{U_{0 } }  +   h_{U_{0 } }  v_{0 }.
  \end{align}
  Since $ h_{U_{0 } } = h$ outside a sufficiently large compact subset of $ \overline{X} $, we see that $ \phi^{\prime \prime }  $ is tame.
  Moreover, when replacing $ \phi ^{\prime } $ by $ \phi ^{\prime \prime } $, the diagram \eqref{dia:zbro} still commutes up to homotopy.
  Therefore, $ \phi ^{\prime \prime } $ is a quasi-isomorphism.
  Thus, $ D$ and $ \phi ^{\prime \prime } $ are respectively the desired tame complex and the desired tame quasi-isomorphism.

  \underline{Step 2 : induction on $ q\in \mathbf{N} $}.
  By Step 1, we assume that $ \mathcal{G} $ is concentrated in positive degrees.
  Let us show our proposition for such $ \mathcal{G} $ by induction on $ q \in \mathbf{N} $.

  By \eqref{eq:wf1p}, we have an injective tame $ \mathscr{O}^{\infty} $-morphism on $ \overline{X} $, 
  \begin{align}\label{eq:vwdv}
   \mathcal{H}^{0} \mathcal{G} \to \mathcal{G}^{0 }.
  \end{align}
  Consider $ \mathcal{H}^{0}\mathcal{G}  $ as a tame $ \mathscr{O} ^{\infty} $-complex concentrated at degree $ 0 $. 
  We have a tame morphism of $ \mathscr{O}^{\infty} $-complexes on $ \overline{X} $,
  \begin{align}\label{eq:sy3g}
   \mathcal{H}^{0} \mathcal{G}  \to \mathcal{G}. 
  \end{align}
  
  \underline{Step 2.1 : the case $ q= 0 $.} 
  Assume now $ q= 0 $. 
  By \eqref{eq:in1y} and \eqref{dia:l3ut}, 
  \begin{enumerate}[\indent a)]
    \item on $ U_{0 } $, if $ HG_{U_{0 } }  $ is the cohomology of $ G_{U_{0 } } $, we have a canonical isomorphism of constant sheaves,  
    \begin{align}\label{eq:stra}
     \alpha_{U_{0 } } : \mathscr{O}^{\infty} \left(H^{0 } G_{U_{0 } } \right) \simeq \left(\mathcal{H}^{0} \mathcal{G}\right)_{|U_{0 } };
    \end{align}
    \item on $ U \in \mathcal{U} \backslash \left\{U_{0 } \right\}$, we have an exact $ \mathscr{O}^{\infty} $-complex on $ U$, 
    \begin{equation}\label{dia:ff3z}
\begin{tikzsd}
         0 & \mathscr{O}_{}^{\infty}  \left(R^{-\ell }_{U } \right)  &	\cdots  & \mathscr{O}_{}^{\infty} \left(R^{0 }_{U}\right)   & \left(\mathcal{H}^{0} \mathcal{G} \right)_{|U }  & 0.
          \arrow["", from=1-1, to=1-2]
          \arrow["", from=1-2, to=1-3]
          \arrow["", from=1-3, to=1-4]
          \arrow["\alpha _{U} ", from=1-4, to=1-5]
          \arrow["", from=1-5, to=1-6]
       \end{tikzsd} 
     \end{equation} 
  \end{enumerate} 
  
  We will establish Step 2.1 by induction on $ \ell \in \mathbf{N} $.
   
  If $ \ell = 0 $, by \eqref{eq:stra} and \eqref{dia:ff3z}, $ \mathcal{H}^{0} \mathcal{G} $ is locally free on $ \overline{X} $.
  There exists a smooth vector bundle $ D $ on $ \overline{X} $ such that  
  \begin{align}\label{eq:n1tt}
   \mathcal{H}^{0} \mathcal{G}   =  \mathscr{O}^{\infty} \left(D\right). 
  \end{align}
  Moreover, $ D$  inherits from \eqref{eq:stra} a tame structure.
  We consider $ D$ as a tame complex of vector bundles concentrated at degree $ 0$.
  Since $ q= 0 $, \eqref{eq:sy3g} is a quasi-isomorphism.
  Thus, $ D $ and \eqref{eq:sy3g} are the desired tame $ \mathscr{O}^{\infty} $-complex and the desired tame $ \mathscr{O}^{\infty} $-quasi-isomorphism.
  This finishes the proof of Step 2.1 when $ \ell = 0 $.
  
  Assume now $ \ell \ge 1$ and that our proposition holds for all the $ \mathcal{G}^{\prime } $ with similar properties and with the obvious associated constants $ p^{\prime},q^{\prime }, \ell^{\prime }  $ satisfying $ p^{\prime } = q^{\prime }  $ and $ \ell'\le \ell-1$.
  As in \cite[(6.3.7)]{BismutShenWei23}, we extend the constant vector bundles $ H^{0 } G_{U_{0 } } $ and $ R_{U}^{0 } \big|_{U\in  \mathcal{U} \backslash \left\{U_{0 } \right\}}  $ trivially to $ \overline{X} $.
  Put
  \begin{align}\label{eq:cwdu}
   R= H^{0 } G_{U_{0 } } \oplus \bigoplus_{U\in \mathcal{U} \backslash \left\{U_{0 } \right\}} R^{0 } _{U}. 
  \end{align}
  Then, $ R$ is a trivial vector bundle on $ \overline{X} $.
  
  By \eqref{eq:klag}, we see that $ \phi_{U} \alpha_{U} $ is well-defined
  on $ \overline{X} $ (cf. Footnote \ref{foot:Us}).
  Put
  \begin{align}\label{eq:3gtu}
    r=  \sum_{U\in \mathcal{U} }^{} \phi_{U} \alpha _{U}.
  \end{align}
  Then, 
  \begin{align}\label{eq:ii2l}
   r: \mathscr{O}_{}^{\infty} \left(R\right)\to \mathcal{H}^{0} \mathcal{G}
  \end{align}
  is a surjective $ \mathscr{O}_{}^{\infty} $-morphism.
  Outside a sufficiently large compact subset of $ \overline{X} $, all the $ \phi_{U} $ vanish except $ \phi_{U_{0 } } $ is equal to $ 1$.
  Thus, $ r$  is tame on $ \overline X$.
  
  Composing \eqref{eq:vwdv} and \eqref{eq:ii2l}, we get a new tame $ \mathscr{O}_{}^{\infty}  $-complex on $ \overline{X} $, 
  \begin{equation}\label{dia:fbgs}
\begin{tikzsd}
      0 & \mathscr{O}_{}^{\infty} \left(R\right) &	\mathcal{G}^{0}  & \mathcal{G}^{1}  & \cdots
       \arrow["", from=1-1, to=1-2]
       \arrow["", from=1-2, to=1-3]
       \arrow["", from=1-3, to=1-4]
       \arrow["", from=1-4, to=1-5]
    \end{tikzsd} 
  \end{equation}
  
  We claim that we can apply our induction assumption to \eqref{dia:fbgs}.
  Indeed, the cohomology of \eqref{dia:fbgs} is concentrated at degree $ -1$.
  We need only to show that if $ U \in \mathcal{U} \backslash \left\{U_{0 } \right\}$, the lengths of the analogue objects in \eqref{dia:ff3z} are $ \le \ell-1$.
  This is exactly what is shown in \cite[(6.3.12), (6.3.13)]{BismutShenWei23}.

  By induction assumptions, there exist a tame $ \mathscr{O}^{\infty} $-complex $ D'$ on $ \overline{X} $ concentrated in degrees $[-\ell,-1]$  and a tame $ \mathscr{O}^{\infty} $-quasi-isomorphism given by the vertical arrows of the diagram, 
  \begin{equation}
    \label{dia:vqpw2}
\begin{tikzsd}
         \cdots & \mathscr{O}^{\infty}\left( D^{\prime, -2} \right) & \mathscr{O}^{\infty}\left( D^{\prime, -1} \right)  & 0  &0  &\cdots  \\
         \cdots & 0 & \mathscr{O}^{\infty} \left(R \right)& \mathcal{G}^{0 }   & \mathcal{G}^{1}   & \cdots
          \arrow["", from=1-1, to=1-2]
          \arrow["", from=1-2, to=1-3]
          \arrow["", from=1-3, to=1-4]
          \arrow["", from=1-4, to=1-5]
          \arrow["", from=1-5, to=1-6]
          \arrow["", from=1-2, to=2-2]
          \arrow["", from=2-1, to=2-2]
          \arrow["", from=2-2, to=2-3]
          \arrow["", from=2-3, to=2-4]
          \arrow["", from=2-4, to=2-5]
          \arrow["", from=2-5, to= 2-6]
          \arrow["", from=1-3, to=2-3]
          \arrow["", from=1-4, to=2-4]
          \arrow["", from=1-5, to=2-5]
        \end{tikzsd}
    \end{equation}
    
  We can rewrite the diagram \eqref{dia:vqpw2} in an obvious way as
  \begin{equation}
    \label{dia:vqpw3}
\begin{tikzsd}
         \cdots & \mathscr{O}^{\infty} \left(D^{\prime, -2} \right) & \mathscr{O}^{\infty}\left( D^{\prime, -1}  \right) & \mathscr{O}^{\infty} \left(R \right) &0  &\cdots  \\
         \cdots & 0 & 0  & \mathcal{G}^{0}_{}  & \mathcal{G}^{1}_{}   & \cdots
          \arrow["", from=1-1, to=1-2]
          \arrow["", from=1-2, to=1-3]
          \arrow["", from=1-3, to=1-4]
          \arrow["", from=1-4, to=1-5]
          \arrow["", from=1-5, to=1-6]
          \arrow["", from=1-2, to=2-2]
          \arrow["", from=2-1, to=2-2]
          \arrow["", from=2-2, to=2-3]
          \arrow["", from=2-3, to=2-4]
          \arrow["", from=2-4, to=2-5]
          \arrow["", from=2-5, to= 2-6]
          \arrow["", from=1-3, to=2-3]
          \arrow["", from=1-4, to=2-4]
          \arrow["", from=1-5, to=2-5]
        \end{tikzsd}
    \end{equation}
  By a direct verification, or as detailed in \cite[(6.3.15)-(6.3.19)]{BismutShenWei23}, in \eqref{dia:vqpw3}, the first line is our desired tame $ \mathscr{O}^{\infty} $-complex, the vertical arrows give the desired tame $ \mathscr{O}^{\infty} $-quasi-isomorphism.
  This finishes the proof of the induction argument on $ \ell\in \mathbf{N} $ and completes the proof for Step 2.1.
  
  \underline{Step 2.2 : the case $ q\ge 1$.}
  Assume now $ q\ge 1$ and our result holds for all $ \mathcal{G}^{\prime }  $ with similar properties and with $ q^{\prime }-p^{\prime } \le q-1$.
  The proof is essentially the same as the one given in \cite[(6.3.20)-(6.3.27)]{BismutShenWei23}.
  Special care is required to ensure all the involved objects and morphisms can be chosen to be tame.
  
  Apply the induction assumption to $ \tau_{\le q-1} \mathcal{G} $.
  There exist a tame smooth complex of vector bundles  $ D_{1} $ on $ \overline{X} $ concentrated in degrees $ [-\ell,q-1]$ and a tame $ \mathscr{O}^{\infty} $-quasi-isomorphism on $ \overline{X} $, 
  \begin{align}\label{eq:wady}
   \mathscr{O}^{\infty} \left(D_{1}\right) \to \tau_{\le q-1} \mathcal{G}. 
  \end{align}
  
  Composing \eqref{eq:wady} with the canonical morphism $ \tau_{\le q-1} \mathcal{G} \to \mathcal{G} $ and applying the cone construction, we get a triangle, 
  \begin{equation}\label{dia:1jsw}
\begin{tikzsd}
      \mathscr{O}^{\infty} \left(D_{1} \right) & \mathcal{G}  & {\rm cone}\left(\mathscr{O}^{\infty} \left(D_{1}\right), \mathcal{G} \right) & \mathscr{O}^{\infty} \left(D_{1}^{\bullet +1} \right). 
       \arrow["", from=1-1, to=1-2]
       \arrow["", from=1-2, to=1-3]
       \arrow["", from=1-3, to=1-4]
    \end{tikzsd} 
  \end{equation}
  We claim that there is a tame homotopy equivalence, 
  \begin{align}\label{eq:ttto}
  {\rm cone}\left( {\rm cone}^{\bullet -1} \left(\mathscr{O}^{\infty} \left(D_{1}\right) , \mathcal{G} \right), \mathscr{O}^{\infty}\left(D_{1}\right) \right)  \simeq \mathcal{G}.
  \end{align}
  Indeed, the morphism \eqref{eq:ttto}, its homotopy inverse, and all the other involved maps are explicitly given in the proof of \cite[Lemma 1.4.2]{KashiwaraSchpira}, which are all tame.
  
  Since the cohomology of $ {\rm cone}\left(\mathscr{O}^{\infty} \left(D_{1}\right), \mathcal{G} \right)$ is concentrated at degree $q$ and is given by $ \mathcal{H} ^{q}\mathcal{G}$, by our induction assumptions, there exist a tame smooth complex of vector bundles $ D_{2} $ concentrated in degrees $ [q-\ell,q]$ and a tame $ \mathscr{O}^{\infty} $-quasi-isomorphism on $ \overline{X} $, 
  \begin{align}\label{eq:3jky}
   \mathscr{O}^{\infty} \left(D_{2}\right) \to {\rm cone}\left(\mathscr{O}^{\infty} \left(D_{1}\right), \mathcal{G} \right).
  \end{align}
  
  Using the last arrow in \eqref{dia:1jsw} and \eqref{eq:3jky}, we get to a tame sequence, 
  \begin{equation}\label{dia:pfom}
\begin{tikzsd}
      \mathscr{O}^{\infty} \left(D_{2}\right) & {\rm cone}\left(\mathscr{O}^{\infty} \left(D_{1}\right), \mathcal{G} \right) &	\mathscr{O}^{\infty} \left(D_{1}^{\bullet +1}\right).
       \arrow["", from=1-1, to=1-2]
       \arrow["", from=1-2, to=1-3]
    \end{tikzsd} 
  \end{equation}
Applying \cite[Proposition 1.4.1 (TR5)]{KashiwaraSchpira} to \eqref{dia:pfom}, as before all the involved map in the proof of \cite[Proposition 1.4.1 (TR5)]{KashiwaraSchpira} are tame, we get a tame $ \mathscr{O}^{\infty} $-quasi-isomorphism on $ \overline{X} $,
  \begin{align}\label{eq:gnco}
    {\rm cone}\left( \mathscr{O}^{\infty}\left( D_{2} ^{\bullet -1} \right), \mathscr{O}^{\infty}\left(D_{1} \right)\right)  \simeq {\rm cone}\left( {\rm cone}^{\bullet -1} \left(\mathscr{O}^{\infty} \left(D_{1}\right) , \mathcal{G} \right), \mathscr{O}^{\infty}\left(D_{1}\right) \right) .
  \end{align}

  By \eqref{eq:ttto} and \eqref{eq:gnco}, we get a tame $ \mathscr{O}^{\infty} $-quasi-isomorphism on $ \overline{X} $,
  \begin{align}\label{eq:vgv3}
   {\rm cone}\left( \mathscr{O}^{\infty}\left( D_{2} ^{\bullet -1} \right), \mathscr{O}^{\infty}\left(D_{1} \right)\right)  \simeq \mathcal{G}.
  \end{align}
  Thus, $ {\rm cone}\left( D_{2} ^{\bullet -1} , D_{1} \right)$ is the desired tame complex of vector bundles which is concentrated in degrees $ [-\ell,q]$ and \eqref{eq:vgv3} is the desired tame quasi-isomorphism.

  This finishes the proof of our induction argument on $ q\in \mathbf{N} $, and completes the proof for Step 2.2. 
  The proof of our proposition is complete.
  \end{proof}

Let $ \mathcal{F}  $ be an object in $ {\rm C^{b,t}_{coh}} \left(\overline{X} \right) $.
We use the notation in \eqref{eq:d2vy}.
By Theorem \ref{thm:yndp} and \eqref{eq:1qmc}, $ \mathscr{D}^{0 }  \mathcal{F}$ is an $ \mathscr{O}^{\infty}  $-complex with same properties as $ \mathcal{G} $. 
By Proposition \ref{prop:cfiv}, there exist a tame bounded smooth complex of vector bundles $ \left(D,v_{0 } \right)$ and tame $ \mathscr{O}^{\infty} $-quasi-isomorphism,  
\begin{align}\label{eq:dknd}
 \phi_{0 } : \mathscr{O}^{\infty} \left(D\right)\to \mathscr{D}^{0 } \mathcal{F}.
\end{align}

Put 
\begin{align}\label{eq:bsce}
 E= \Lambda \left(\overline{T^{*} X } \right) \widehat{\otimes }  D.
\end{align}
Clearly, $ E$ extends to an open neighbourhood of $ \overline{X} $, so that $ E$ defines a smooth $ \mathbf{Z} $-graded vector bundle on $ \overline{X} $. 
Let $ \mathcal{E} $ be the $ \mathbf{Z} $-graded sheaf of smooth sections of $ E$.

\begin{prop}\label{prop:axew}
There exist a tame antiholomorphic superconnection $\left(E, A^{E\prime \prime } \right)$ on $ \overline{X} $ and a $ \mathscr{D} $-tame $ \mathscr{O} $-quasi-isomorphism $ \phi : \mathcal{E} \to \mathscr{D} \mathcal{F} $ on $ \overline{X} $  of the following forms, 
\begin{align}\label{eq:gmbo}
  &A^{E\prime \prime } = v_{0 } + \nabla^{D\prime \prime } + \sum_{i\ge 2} v_{i},& \phi = \sum_{i\ge 0 } \phi_{i},
 \end{align}
where  $ \nabla^{D\prime \prime } $ is an antiholomorphic part of some connection on $ D$, for each $ i \ge 2$, $ v_{i} $ is a smooth $ \left(0, i\right)$-form with values in $ \End^{1-i} \left(D\right)$, and for each $ i \ge 1$, $ \phi_{i} : \mathscr{O}^{\infty} \left(D^{\bullet } \right) \to \mathscr{D}^{i} \otimes_{\mathscr{O} } \mathcal{F}^{\bullet -i}     $ is a degree preserving morphism of $ \mathbf{Z} $-graded $ \mathscr{O}^{\infty} $-modules on $ \overline{X}$.
\end{prop}
\begin{proof}
Let us first construct a tame antiholomorphic superconnection $ A^{E\prime \prime } $ and a $ \mathscr{D} $-tame $ \mathscr{O} $-morphism $ \phi$, i.e., 
\begin{align}\label{eq:yld3}
& \left(A^{E\prime \prime } \right)^{2} = 0, & \phi A^{E\prime \prime } = A^{\mathscr{D} \mathcal{F} \prime \prime } \phi,
\end{align}
or equivalently,
 \begin{align}\label{eq:zyko}
  \begin{pmatrix}
  A^{E\prime \prime }  & 0 \\
  \phi \left(-1\right)^{{\rm deg}}  & A^{\mathscr{D} \mathcal{F} \prime \prime }  
  \end{pmatrix}^{2} = 0. 
\end{align}
 
We take an arbitrary $ \nabla^{D\prime \prime }_{0 }  $ which restricts to  the antiholomorphic part of the trivial connection on $ U_{0 } $.
Let us construct a tame pair $ \left(v_{1}, \phi_{1}\right) $ such that 
 \begin{align}\label{eq:bcal}
  \begin{pmatrix}
    v_{0 } +\nabla^{D\prime \prime }_{0} +v_{1}   & 0 \\
    \left(\phi_{0 } +\phi_{1} \right) \left(-1\right)^{{\rm deg}}  & A^{\mathscr{D} \mathcal{F} \prime \prime } 
    \end{pmatrix}^{2}_{\le 1}  = 0.
 \end{align}
 Here, the subscript $ {\le 1} $ indicates that we only consider terms involving $ \Lambda^{\le 1} \left(\overline{T^{*} X} \right)$, and similar notation will appear later.
 Since $ \phi_{0} $ in \eqref{eq:dknd} is a morphism of $ \mathscr{O} ^{\infty} $-complexes, we know that \eqref{eq:bcal} is equivalent to
 \begin{align}\label{eq:gomp}
  \left[ \begin{pmatrix}
    v_{0 }    & 0 \\
   \phi_{0}  \left(-1\right)^{{\rm deg}}  & v^{\mathcal{F} }
    \end{pmatrix},\begin{pmatrix}
      \nabla^{D\prime \prime }_{0}    & 0 \\
     0   & \overline{\partial}
      \end{pmatrix}\right] + \left[ \begin{pmatrix}
        v_{0 }    & 0 \\
       \phi_{0}  \left(-1\right)^{{\rm deg}}  & v^{\mathcal{F} }
        \end{pmatrix},\begin{pmatrix}
         v_{1 }    & 0 \\
         \phi_{1}  \left(-1\right)^{{\rm deg}}  & 0 
          \end{pmatrix}\right]= 0. 
 \end{align}
 
 In \cite[(6.3.55)-(6.3.62)]{BismutShenWei23}, the authors construct a (not necessarily tame) solution\footnote{The proof given in the above reference does not rely on the compactness assumption.
 } $\left( v_{1}^{\prime }, \phi_{1}^{\prime }\right) $ to \eqref{eq:gomp}.
 Note that on $ U_{0 } $, $ v_{1} ^{\prime \prime } = 0$ and $ \phi_{1}^{\prime \prime }  = 0$ is another obvious solution.
 Put 
 \begin{align}\label{eq:ept1}
  &v_{1} = \left(1-\phi_{U_{0 } }\right)v_{1} ^{\prime },  & \phi _{1} = \left(1-\phi_{U_{0 } }\right)\phi _{1} ^{\prime }.
 \end{align}
 Then, $\left( v_{1}, \phi_{1}\right) $ is a solution to \eqref{eq:gomp}.
 By \eqref{eq:klag} and \eqref{eq:ityw}, $\left( v_{1}, \phi_{1}\right) $ is tame on $ \overline{X}$.

 Assume now $ i\ge 2$ and we have found tame $v_{1}, v_{2}, \ldots,v_{i-1}  $ and $\phi_{1},  \phi_{2}, \ldots, \phi _{i-1}  $, such that if
 \begin{align}\label{eq:gmbo1}
  &A^{E\prime \prime }_{\le i-1}  = v_{0 } + \nabla^{D\prime \prime }_{0 }  + \sum_{j=  1}^{i-1}  v_{j},& \phi_{\le i-1}  = \sum_{j=  0 }^{i-1}  \phi_{j},
 \end{align}
 then 
 \begin{align}\label{eq:lm3n}
  \begin{pmatrix}
    A^{E\prime \prime }_{\le i-1}    & 0 \\
    \phi_{\le i-1}  \left(-1\right)^{{\rm deg}}  & A^{\mathscr{D} \mathcal{F} \prime \prime } 
    \end{pmatrix}^{2}_{\le i-1}  = 0.
 \end{align}

 Let us construct a tame solution $ \left(v_{i}, \phi_{i}\right) $ to the equation  
 \begin{align}\label{eq:zcg1}
  \begin{pmatrix}
    A^{E\prime \prime }_{\le i-1} + v_{i}    & 0 \\
    \left(\phi_{\le i-1} + \phi_{i}\right)  \left(-1\right)^{{\rm deg}}  & A^{\mathscr{D} \mathcal{F} \prime \prime }  
    \end{pmatrix}^{2}_{\le i}  = 0.
 \end{align}
 By \eqref{eq:lm3n},  we know that \eqref{eq:zcg1} is equivalent to 
 \begin{align}\label{eq:wbzq}
  \begin{pmatrix}
  A^{E\prime \prime }_{\le i-1}    & 0 \\
    \phi_{\le i-1}  \left(-1\right)^{{\rm deg}}  & A^{\mathscr{D} \mathcal{F} \prime \prime } 
    \end{pmatrix}^{2}_{=i}
    + \left[ \begin{pmatrix}
      v_{0 }    & 0 \\
     \phi_{0}  \left(-1\right)^{{\rm deg}}  & v^{\mathcal{F} }
      \end{pmatrix},\begin{pmatrix}
       v_{i}    & 0 \\
       \phi_{i}  \left(-1\right)^{{\rm deg}}  & 0 
        \end{pmatrix}\right]= 0. 
 \end{align}
 A (not necessarily tame) solution $ \left(v_{i}^{\prime } , \phi^{\prime } _{i} \right)$ to \eqref{eq:wbzq} is constructed in \cite[(6.3.63)-(6.3.69)]{BismutShenWei23}.
  Also, outside a compact subset of $ \overline{X} $, another obvious solution is $ v_{i}^{\prime \prime }  = 0$ and $ \phi_{i}^{\prime \prime } = 0$.
 By the same construction as \eqref{eq:ept1}, we get a tame solution $ \left(v_{i}, \phi_{i} \right)$ to \eqref{eq:wbzq}.
 
 This procedure stops at $ i= n $. Ultimately,  we get a tame solution $ \left(A^{E\prime \prime },  \phi\right) = \left( A^{E\prime \prime }_{\le n}, \phi _{\le n} \right)$ to \eqref{eq:zyko}.

 Since $ \phi_{0 } $ is a quasi-isomorphism in \eqref{eq:dknd}, by Proposition \ref{prop:ceuw}, we know that $ \phi : \mathcal{E} \to \mathscr{D} \mathcal{F} $ is a quasi-isomorphism.
 The proof of our proposition is complete.
\end{proof}

\begin{proof}[Proof of Theorem \ref{thm:poev} \ref{thm:poeva})]
 The part \ref{thm:poeva}) of Theorem \ref{thm:poev} is an immediate consequence of Proposition \ref{prop:axew}.
\end{proof}

\begin{re}\label{re:dqvn}
 If $ c_{1} \in \mathbf{Z} $ and if $ \mathcal{H} \mathcal{F} $ is concentrated in degrees $ \le c_{1} $, by Proposition \ref{prop:cfiv},  the diagonal $ \left(D,v_{0 } \right)$ of $ \left(E,A^{E\prime \prime } \right)$ in Theorem \ref{thm:poev} \ref{thm:poeva}) can be chosen to be concentrated in degrees $ \le c_{1} $, so that $ E$ is concentrated in degrees $ \le c_{1}  + n$.
\end{re}

\subsection{Proof of Theorem \ref{thm:poev} \ref{thm:poevb})}\label{s:thmb}
Let $ \mathcal{E} $ be an $ \mathscr{O} $-complex associated to an object in $ {\rm B^{b,t}} \left(\overline{X} \right)$.
Let $ \left(\underline{\mathcal{E}} , \underline{\mathcal{F}} , \underline{\phi} \right)$ be as in Theorem \ref{thm:poev} \ref{thm:poevb}). 

\begin{prop}\label{prop:2gpa}
   If $ \psi : \mathcal{E} \to \mathscr{D} \underline{\mathcal{F}}  $ is a morphism in $ {\rm \underline{\widetilde{B}}  }^{{\rm t}} \left(\overline{X} \right)$, 
  then there exists a morphism $ \widetilde{\psi} $ in  ${\underline{{\rm B}} }^{{\rm b,t}} \left(\overline{X} \right)$ such that the diagram 
  \begin{equation}\label{dia:chiq}
\begin{tikzsd}
         & \underline{{\mathcal{E}}}     \\
        \mathcal{E}   & \mathscr{D} \underline{{\mathcal{F} } }
        \arrow["\widetilde{\psi}  ", from=2-1, to=1-2, dotted]
        \arrow["\underline{\phi } ", from=1-2, to=2-2]
        \arrow["\psi", from=2-1, to=2-2]
      \end{tikzsd}
  \end{equation}  
 in $ \widetilde{\underline{{\rm B}} }^{{\rm t}} \left(\overline{X} \right)$ commutes. 
\end{prop}
\begin{proof}
 Equivalently, we need construct a morphism in $ \widetilde{{\rm B}}^{{\rm t}} \left(\overline{X} \right)$ of the following type, 
 \begin{align}\label{eq:vgoo2}
 \begin{pmatrix}
 \widetilde{\psi}   & 0 \\
 h\left(-1\right)^{{\rm deg}}  & 1
 \end{pmatrix}: {\rm cone}\left(\mathcal{E},\mathscr{D} \underline{\mathcal{F}   } \right)\to {\rm cone}\left(\underline{\mathcal{E} }, \mathscr{D} \underline{\mathcal{F} } \right),
 \end{align}
 where $ h$ is the one making the diagram \eqref{dia:chiq} commute in homotopy categories.
 More precisely, we need find a tame solution $ \left(\widetilde{\psi} , h\right)$ to the equation,
 \begin{align}\label{eq:suye}
  A^{{\rm cone}\left(\underline{\mathcal{E} }, \mathscr{D} \underline{\mathcal{F} }  \right)\prime \prime}    \begin{pmatrix}
    \widetilde{\psi }   & 0 \\
    h\left(-1\right)^{{\rm deg}}   & 1
    \end{pmatrix}-
    \begin{pmatrix}
     \widetilde{\psi}   & 0 \\
     h \left(-1\right)^{{\rm deg}}   & 1
   \end{pmatrix}
     A^{{\rm cone}\left({\mathcal{E} }, \mathscr{D} \underline{\mathcal{F} }  \right)\prime \prime} = 0. 
 \end{align}

We fix splittings for $ \mathcal{E} $ and $ \underline{\mathcal{E} } $, which are induced by \eqref{eq:aayz}.
We have identifications of $ \mathscr{D} $-modules on $ \overline{X} $, 
 \begin{align}\begin{aligned}\label{eq:qvjg}
  &{\rm cone}\left(\mathcal{E},\mathscr{D} \underline{\mathcal{F}   } \right)= \mathscr{D} \otimes_{\mathscr{O}^{\infty} }  {\rm cone}\left(\mathscr{O}^{\infty} \left(D\right),\mathscr{D}^{0} \underline{\mathcal{F}}  \right),\\
  &{\rm cone}\left(\underline{\mathcal{E}} ,\mathscr{D} \underline{\mathcal{F}   } \right)= \mathscr{D} \otimes_{\mathscr{O}^{\infty} }  {\rm cone}\left(\mathscr{O}^{\infty} \left(\underline{D} \right),\mathscr{D}^{0} \underline{\mathcal{F}}  \right).
 \end{aligned}\end{align}
Put 
\begin{align}\label{eq:lg2r}
 \mathcal{P} =  \Hom\left({\rm cone}\left(\mathscr{O}^{\infty} \left( {D} \right),\mathscr{D}^{0 }  \underline{\mathcal{F} } \right),{\rm cone}\left(\mathscr{O}^{\infty} \left(\underline{D}\right), \mathscr{D}^{0 }  \underline{\mathcal{F} } \right)\right).
\end{align}
By \eqref{eq:qvjg} and \eqref{eq:lg2r}, we have\footnote{\label{Footenote:Hom}Here, $ {\rm Hom}$ is the one defined before \eqref{eq:vgot}.} 
\begin{align}\label{eq:xolz}
 \mathscr{D} \mathcal{P} = \Hom\left({\rm cone}\left(\mathcal{E} ,\mathscr{D} \underline{\mathcal{F} } \right),{\rm cone}\left(\underline{\mathcal{E} }, \mathscr{D} \underline{\mathcal{F} } \right)\right).
\end{align}
Then, $ \mathcal{P} $ is a complex equipped a differential $ v^{\mathcal{P} } $ induced by $ \psi_{0 } $ and $ \underline{\phi }_{0} $.
The differential of $ \mathscr{D} \mathcal{P} $ has an expansion analogue to \eqref{eq:ktja}, 
\begin{align}\label{eq:pteq}
 A^{\mathscr{D} \mathcal{P} \prime \prime } = v^{\mathcal{P} } +\overline{\partial}+ \ldots 
\end{align}

We can rewrite \eqref{eq:suye} as 
\begin{align}\label{eq:fsib}
 A^{\mathscr{D}\mathcal{P} \prime \prime  }  \begin{pmatrix}
  \widetilde{\psi }   & 0 \\
  h\left(-1\right)^{{\rm deg}}   & 1
  \end{pmatrix}= 0. 
\end{align}
Let us construct, by induction on the antiholomorphic degrees, a tame solution $\left( \widetilde{\psi }, h\right)=\left( \sum_{i= 0 }^{n} \widetilde{\psi }_{i} , \sum_{i= 0}^{n} h_{i} \right) $ to \eqref{eq:fsib}.

By \cite[\href{https://stacks.math.columbia.edu/tag/0649}{Tag 0649}]{stacks-project}, there exists a morphism of $ \mathscr{O}^{\infty} $-complexes on $ \overline{X} $ of the following type
\begin{align}\label{eq:vgoo1}
  \begin{pmatrix}
  \widetilde{\psi}_{0 }^{\prime }    & 0 \\
  h_{0 }^{\prime }\left(-1\right)^{{\rm deg}}     & 1
  \end{pmatrix}: {\rm cone}\left(\mathscr{O}^{\infty} \left(D\right),\mathscr{D}^{0 }  \underline{\mathcal{F}   } \right)\to {\rm cone}\left(\mathscr{O}^{\infty} \left(\underline{D }\right), \mathscr{D}^{0 }  \underline{\mathcal{F} } \right). 
  \end{align}
 Using the above reference again on $ U_{0 } $, we have a constant morphism of constant complexes, 
 \begin{align}\label{eq:vgoo}
  \begin{pmatrix}
  \widetilde{\psi}_{0 }^{\prime \prime }    & 0 \\
  h_{0 }^{\prime \prime }\left(-1\right)^{{\rm deg}}     & 1
  \end{pmatrix}: {\rm cone}\left(\mathscr{O}^{\infty} \left(D\right),\mathscr{D}^{0 }  \underline{\mathcal{F}   } \right)_{|U_{0 } } \to {\rm cone}\left(\mathscr{O}^{\infty} {\left(\underline{D }\right)} , \mathscr{D}^{0 }  \underline{\mathcal{F} } \right)_{|U_{0 } } . 
  \end{align} 
  As before, \eqref{eq:vgoo1} and \eqref{eq:vgoo} are respectively equivalent to
  \begin{align}\label{eq:3fg3}
&   v^{\mathcal{P} } \begin{pmatrix}
  \widetilde{\psi}_{0 }^{\prime }    & 0 \\
  h_{0 }^{\prime }\left(-1\right)^{{\rm deg}}     & 1
  \end{pmatrix}= 0, &v^{\mathcal{P} } _{|U_{0 } } \begin{pmatrix}
  \widetilde{\psi}_{0 }^{\prime \prime }    & 0 \\
  h_{0 }^{\prime \prime }\left(-1\right)^{{\rm deg}}     & 1
  \end{pmatrix}= 0.
  \end{align}
  Put 
  \begin{align}\label{eq:xhjm}
   &\widetilde{\psi}_{0 } = \left(1-\phi_{U_{0 } }\right) \widetilde{\psi}_{0 }^{\prime } + \phi_{U_{0 } } \widetilde{\psi}_{0 }^{\prime \prime }, & h_{0 } = \left(1-\phi_{U_{0 } }\right) h_{0 }^{\prime } + \phi_{U_{0 } } h_{0 }^{\prime \prime }.
  \end{align}
  As before, $ \widetilde{\psi}_{0 }$ and $ h_{0 } $ are tame.
  By \eqref{eq:3fg3} and \eqref{eq:xhjm}, we have 
  \begin{align}\label{eq:ubla}
    \left\{A^{\mathscr{D}\mathcal{P} \prime \prime  }  \begin{pmatrix}
      {\widetilde{\psi } }_{0 }   & 0 \\
      h_{0 } \left(-1\right)^{{\rm deg}}   & 1
      \end{pmatrix}\right\}_{= 0 } = v^{\mathcal{P} } \begin{pmatrix}
      {\widetilde{\psi } }_{0 }   & 0 \\
      h_{0 } \left(-1\right)^{{\rm deg}}   & 1
      \end{pmatrix}= 0. 
  \end{align} 
 
  Assume now $ i\ge 1$ and that we have constructed tame $ \widetilde{\psi} _{0},\ldots, \widetilde{\psi} _{i-1} $ and  $ h_{0 },\ldots, h_{i-1} $ such that if 
  \begin{align}\label{eq:21wy}
   &\widetilde{\psi} _{\le i-1} = \sum_{j=0}^{i-1} {\widetilde{\psi}  }_{j}, & h_{\le i-1} = \sum_{j=0}^{i-1} h_{j},
  \end{align}
  then
  \begin{align}\label{eq:ltc1}
    \left\{A^{\mathscr{D}\mathcal{P} \prime \prime  }  \begin{pmatrix}
      {\widetilde{\psi}  }_{\le i-1}   & 0 \\
      h_{\le i-1 } \left(-1\right)^{{\rm deg}}   & 1
      \end{pmatrix}\right\}_{\le i-1 } = 0. 
  \end{align}
  
  We need find a tame solution $ \left({\widetilde{\psi}  }_{i}, h_{i} \right)$ to the equation,  
  \begin{align}\label{eq:ltc2}
    \left\{A^{\mathscr{D}\mathcal{P} \prime \prime  }  \begin{pmatrix}
      {\widetilde{\psi}  }_{\le i-1} + {\widetilde{\psi}  }_{i}   & 0 \\
      \left(h_{\le i-1 }+h_{i} \right) \left(-1\right)^{{\rm deg}}   & 1
      \end{pmatrix}\right\}_{\le i } = 0,
  \end{align}
By \eqref{eq:ltc1}, we see that \eqref{eq:ltc2} is equivalent to 
\begin{align}\label{eq:ltc3}
  v^{\mathcal{P} }  \begin{pmatrix}
  {\widetilde{\psi}  }_{i}  & 0 \\
  h_{i}\left(-1\right)^{{\rm deg}}   & 0 
  \end{pmatrix}     = -\left\{A^{\mathscr{D}\mathcal{P} \prime \prime  }  \begin{pmatrix}
      {\widetilde{\psi}  }_{\le i-1}   & 0 \\
      h_{\le i-1 }\left(-1\right)^{{\rm deg}}   & 1
      \end{pmatrix}\right\}_{=  i } . 
\end{align}

Note also that the right-hand side of \eqref{eq:ltc3} has the following form,
\begin{align}\label{eq:ytor}
-  \left\{A^{\mathscr{D} \mathcal{P} \prime \prime } \begin{pmatrix}
    {\widetilde{\psi}  }_{\le i-1}  & 0 \\
    h_{\le i-1}\left(-1\right)^{{\rm deg}}   & 1
    \end{pmatrix}\right\}_{= i}= \begin{pmatrix}
    u_{i}  & 0 \\
    v_{i}  & 0 
    \end{pmatrix}.
\end{align}
Put 
\begin{align}\label{eq:1txs}
 \mathcal{R} = \Hom\left(\mathscr{O}^{\infty} \left( {D} \right),{\rm cone}\left(\mathscr{O}^{\infty} \left(\underline{D}\right), \mathscr{D}^{0}  \underline{\mathcal{F} } \right)\right).
\end{align}
With obvious notation, we can rewrite \eqref{eq:ltc3} as an equation in $ \mathscr{D}^{i} \mathcal{R} $,
\begin{align}\label{eq:uo31}
  v^{\mathcal{R} }  \begin{pmatrix}
  \widetilde{\psi}_{i}  \\
  h_{i}\left(-1\right)^{{\rm deg}}   
  \end{pmatrix}     = \begin{pmatrix}
  u_{i} \\ v_{i} 
  \end{pmatrix}.
\end{align}

We claim that $ \begin{pmatrix}
u_{i} \\
v_{i} 
\end{pmatrix}$ is $ v^{\mathcal{R} } $-closed.
Indeed, we have the obvious identity
\begin{align}\label{eq:xksw}
  A^{\mathscr{D} \mathcal{P} \prime \prime } A^{\mathscr{D} \mathcal{P} \prime \prime } \begin{pmatrix}
   {\widetilde{\psi}  }_{\le i-1}  & 0 \\
   h_{\le i-1}\left(-1\right)^{{\rm deg}}   & 1
   \end{pmatrix}= 0. 
 \end{align}
 By \eqref{eq:ltc1}, the degree $ i$ part of \eqref{eq:xksw} is given by 
 \begin{align}\label{eq:c2lg}
  v^{\mathcal{P} }  \left\{A^{\mathscr{D} \mathcal{P} \prime \prime } \begin{pmatrix}
   {\widetilde{\psi}  }_{\le i-1}  & 0 \\
   h_{\le i-1}\left(-1\right)^{{\rm deg}}   & 1
   \end{pmatrix}\right\}_{= i} = 0. 
 \end{align}
 By \eqref{eq:ytor} and \eqref{eq:c2lg}, we see that $ \begin{pmatrix}
  u_{i} \\ v_{i} 
  \end{pmatrix}$ is $ v^{\mathcal{R} } $-closed.

Since $ \underline{\phi } $ is a quasi-isomorphism, by Proposition \ref{prop:ceuw}, $ \underline{\phi}_{0 } $ is a quasi-isomorphism, so that $ {\rm cone}\left(\mathscr{O}^{\infty} \left(\underline{D}\right), \mathscr{D}^{0}  \underline{\mathcal{F} } \right)$ is acyclic.
As \cite[Proposition 6.3.4]{BismutShenWei23}, the complex $ \left(\mathscr{D}^{i}  \mathcal{R},v^{\mathcal{R} } \right) $ as well as the complex of its global sections on $ \overline{X} $ is acyclic.
Therefore, the equation \eqref{eq:uo31} has a solution, which is not necessarily tame.
However, outside a sufficiently large compact subset of $ \overline{X} $, $ \begin{pmatrix}
0 \\
0 
\end{pmatrix} $ is an obvious solution to \eqref{eq:uo31}.
As \eqref{eq:ept1}, we get a tame solution for \eqref{eq:uo31}, hence also a tame solution to \eqref{eq:ltc2}.

The procedure stops at $ i= n$.
We obtain a tame solution $ \left(\widetilde{\psi}, h\right)= \left(\widetilde{\psi}_{\le n} , h_{\le n} \right)$ to \eqref{eq:fsib}, which finishes the proof of our proposition.
\end{proof}

\begin{prop}\label{prop:y1hu}
 If $ \widetilde{\psi } : \mathcal{E} \to \underline{\mathcal{E} } $ is a morphism in $ {\rm \underline{B} ^{b,t}}\left(\overline{X} \right)$ such that the diagram 
 \begin{equation}\label{dia:logl}
\begin{tikzsd}
        & \underline{\mathcal{E}}    \\
       \mathcal{E}   & \mathscr{D} \underline{\mathcal{F} }
       \arrow["\widetilde{\psi } ", from=2-1, to=1-2]
       \arrow["\underline{\phi} ", from=1-2, to=2-2]
       \arrow["0 ", from=2-1, to=2-2]
     \end{tikzsd}
 \end{equation}
 in $ \widetilde{\underline{\rm B} }^{{\rm t}} \left(\overline{X} \right) $ commutes, then $ \widetilde{\psi } = 0 $ in $ {\rm \underline{B} _{}^{b,t}}\left(\overline{X} \right)$.
\end{prop}
\begin{proof}
 Clearly, \eqref{dia:logl}  can be rewritten as the commutative diagram in $ \widetilde{\underline{\rm B} }^{{\rm t}} \left(\overline{X} \right) $, 
 \begin{equation}\label{dia:ymmc}
\begin{tikzsd}
      \mathcal{E}  & \underline{\mathcal{E} }  \\
      0  & \mathscr{D}  \underline{\mathcal{F} } .
      \arrow["", from=1-1, to=2-1]
      \arrow["\widetilde{\psi } ", from=1-1, to=1-2]
      \arrow["\underline{\phi } ", from=1-2, to=2-2]
      \arrow["", from=2-1, to=2-2]
    \end{tikzsd}
 \end{equation}
 By \cite[Proposition 1.4.4 (TR1) (TR3) (TR4)]{KashiwaraSchpira} and Proposition \ref{prop:vkig}, the diagram \eqref{dia:ymmc} can be extended to the commutative diagram in $ \widetilde{\underline{{\rm B}} }^{{\rm t}}  \left(\overline{X} \right)$, 
 \begin{equation}\label{dia:ymmc1}
\begin{tikzsd}
    \mathcal{E}  & {\rm cone}^{\bullet -1} \left(\underline{\phi } \right)  \\
    \mathcal{E}  & \underline{\mathcal{E} }  \\
    0  & \mathscr{D} \underline{\mathcal{F} } .
    \arrow["{\rm id}", from=1-1, to=2-1]
    \arrow["\widetilde{\widetilde{\psi } } ", from=1-1, to=1-2,dotted]
    \arrow["", from=1-2, to=2-2]
    \arrow["\widetilde{\psi } ", from=2-1, to=2-2]
    \arrow["", from=2-1, to=3-1]
    \arrow["\underline{\phi } ", from=2-2, to=3-2]
    \arrow["", from=3-1, to=3-2]
  \end{tikzsd}
\end{equation}

Since $ \underline{\phi} $ is a quasi-isomorphism, $ {\rm cone}\left(\underline{\phi } \right)$ is acyclic.
The vertical arrow in the following diagram \eqref{dia:fvrw} is therefore a quasi-isomorphism.
By Proposition \ref{prop:ceuw}, $ {\rm cone}\left(\underline{\phi }_{0 }\right)$ is also acyclic.
Proceeding as in the proof of Proposition \ref{prop:2gpa},\footnote{\label{footnoteDFA}We can not apply direct Proposition \ref{prop:2gpa}, since $ {\rm cone}\left({\underline{\phi}  } \right)$, whose differential contains more than two terms, is not of type $ \mathscr{D} \mathcal{F} $.
However, this does not affect the validity of the proof.} there is a lift, represented by the dotted arrow, so that
\begin{equation}\label{dia:fvrw}
\begin{tikzsd}
    & 0  \\
    \mathcal{E}  & {\rm cone}^{\bullet -1} \left(\underline{\phi }  \right) 
    \arrow[from=1-2, to=2-2]
    \arrow[dotted, from=2-1, to=1-2]
    \arrow["\widetilde{\widetilde{\psi } }"',from=2-1, to=2-2]
  \end{tikzsd}
\end{equation}
commutes in $ \widetilde{\underline{{\rm B}} }^{{\rm t}} \left(\overline{X} \right) $.
Therefore, $ \widetilde{\widetilde{\psi} } = 0 $ in $ \widetilde{\underline{{\rm B}} }^{{\rm t}}  \left(\overline{X} \right)$.
 
By \eqref{dia:ymmc1}, we see that $ \widetilde{\psi } = 0 $ in $ \underline{{\rm B}}^{{\rm b,t}} \left(\overline{X} \right)$ and finish the proof of our proposition. 
\end{proof}

\begin{proof}[Proof of Theorem \ref{thm:poev} \ref{thm:poevb})]
  Theorem \ref{thm:poev} \ref{thm:poevb}) is an immediate consequence of Propositions \ref{prop:2gpa} and \ref{prop:y1hu}.
  The former establishes the existence of a lift, while the latter ensures its uniqueness.
\end{proof}

\section{Direct image and tame antiholomorphic superconnection}\label{S:product}
The purpose of this section is to study the direct image of a tame antiholomorphic superconnection on a product of complex manifolds.
More precisely, if $ S $ is a relatively compact open subset of a complex manifold $ \widetilde{S} $, if $ X$ is a complex manifold, we study the tame antiholomorphic superconnection $ \left(E,A^{E\prime \prime } \right)$ with associated sheaf $ \mathcal{E} $ on the product $ \overline{S} \times X$. If $ \pi : \overline{S} \times X \to \overline{S} $ is the obvious projection, when $ \left(E,A^{E\prime \prime } \right)$ has compact support in $ \overline{S} \times X$, we show that the direct image $ \pi_{*} \mathcal{E} $ has coherent cohomology.

This section is organised as follows. 
In Section \ref{s:superx}, we specialise the construction of Section \ref{S:AntiTame} to the product $ \overline{S} \times X$, and state the main results of this section.


In Section \ref{s:2sub}, we introduce subcomplexes $\left( \pi_{*} \mathcal{E}\right)_{c}, \left(\pi_{*} \mathcal{E}\right)_{\left(2\right)} $ of $ \pi_{*} \mathcal{E} $, and the corresponding diagonals.
The former $ \left( \pi_{*} \mathcal{E}\right)_{c}$ consists of sections with fibrewise compact support, whereas $ \left(\pi_{*} \mathcal{E}\right)_{\left(2\right)}$ consists of sections with $ L^{2} $-higher derivatives, once some metrics are fixed.
If $ \left(E,A^{E\prime \prime } \right)$ has compact support, we show that these three complexes $\left( \pi_{*} \mathcal{E}\right)_{c}, \left(\pi_{*} \mathcal{E}\right)_{\left(2\right)} $, $ \pi_{*} \mathcal{E} $ are homotopy equivalent.

In Section \ref{s:split}, with respect to a non-canonical splitting on $ E$, we construct the associated Laplacians, acting on the fibre of $ \left(\pi_{*} \mathcal{E}\right)_{\left(2\right)}$ as well as its diagonal.
We extend these actions also to the sheaves themselves.

In Section \ref{s:triv}, assuming $ S$ is small enough, we trivialise various constructions on $ S$.

In Section \ref{s:trunc}, we study the resolvent of the various Laplacians.

In Section \ref{s:trunc1}, given a small spectral parameter $ a>0$, on some suitable open set subset $ \widetilde{S} _{a}  \subset \widetilde{S} $, we introduce a spectral decomposition $\left( \pi_{*} \mathcal{E}\right)_{\left(2\right),a} = \left( \pi_{*} \mathcal{E}\right)_{\left(2\right),a,-} \oplus \left( \pi_{*} \mathcal{E}\right)_{\left(2\right),a,+} $ as well as a similar decomposition on their diagonal.
We show that $ \left( \pi_{*} \mathcal{E}\right)_{\left(2\right),a,-}$ is the sheaf associated to some antiholomorphic superconnection, while $ \left( \pi_{*} \mathcal{E}\right)_{\left(2\right),a,+} $ is acyclic.

In Sections \ref{s:p2}-\ref{s:p1}, we show the results stated in Section \ref{s:superx}.

In Section \ref{s:proj}, we reinterpret our main results in the context of the derived category.

Finally, in Section \ref{S:Gs}, we show Theorem \ref{thm:main} \ref{main2a}) in the case when $ X$ is smooth.  
\subsection{Tame antiholomorphic superconnection on a product}\label{s:superx}
Let $ S$ be a relatively compact open subset of some ambient complex manifold $ \widetilde{S} $ of dimension $ m$.
Let $ \overline{S} \subset \widetilde{S} $ be the closure of $ S$. Then, $ \overline{S} $ is compact.
Let $ X$ be a complex manifold of dimension $ n$.

Put
\begin{align}\label{eq:3yiz}
 &M= S \times X, &\overline{M} = \overline{S} \times X, &&\widetilde{M} =  \widetilde{S} \times X.
\end{align}
Then,  $ M$ is an open subset of $ \widetilde{M} $, and $ \overline{M} $ is the closure of $ M$ in $ \widetilde{M} $.
By abuse of notation, the obvious projections $ M \to S$, $ \overline{M} \to \overline{S}$, $\widetilde{M} \to \widetilde{S} $ are denoted by $ \pi$, and $ M \to X$, $ \overline{M} \to X$, $\widetilde{M} \to X $ are denoted by $p$.

If $ s\in \widetilde{S}$, put
\begin{align}\label{eq:jio2}
 X_{s} = \pi^{-1} \left(s\right).
\end{align}
Then, 
\begin{align}\label{eq:2dzb}
 X_{s} = \left\{s\right\}\times X.
\end{align} 

We apply the construction of Section \ref{S:AntiTame} to $ \overline{M} $.
Let $ \left(E,A^{E\prime\prime } \right)$ be an object in $ {\rm B^{b,t}} \left(\overline{M} \right)$, with diagonal $ \left(D,v_{0 } \right)$ and associated $\mathscr{O}  $-complex $ \mathcal{E} $ on $ \overline{M} $.

\begin{prop}\label{prop:njjy}
 There exist a small open neighbourhood $ \widetilde{\widetilde{S} } $ of $ \overline{S} $ and an open subset $ U_{0 } $ of $ X$ with $ X \backslash  U_{0 } $ compact such that $ \left(E,A^{E\prime\prime } \right)$ extends to an antiholomorphic superconnection on $ \widetilde{\widetilde{S}}  \times X$ which is constant on $ \widetilde{{\widetilde{S}}}  \times U_{0 } $.
\end{prop}
\begin{proof} 
 By definition, $ \left(E,A^{E\prime \prime } \right)$ extends to some open neighbourhood of $ \overline{S} \times X$. 
 Since $ \overline{S} $ is compact, by tameness of $ \left(E,A^{E\prime \prime } \right)$, there exists an open subset $ U_{0 } \subset  X$ with $ X \backslash  U_{0 } $ compact such that the above extension is constant on a neighbourhood of $ \overline{S} \times U_{0 } $.
 Since constant objects can always be extended, we can further extend $ \left(E,A^{E\prime \prime } \right)$ to a larger neighbourhood, which contains a set of type $\widetilde{ \widetilde{S}} \times U_{0 } $.
 Taking $ \widetilde{\widetilde{S}}  $ sufficiently small, this larger neighbourhood contains $ \widetilde{\widetilde{S}}  \times X$.
 The proof of our proposition is complete. 
\end{proof}
In the sequel, we will take the above extension for $ \left(E,A^{E\prime \prime } \right)$ and we will assume $ \widetilde{\widetilde{S} } = \widetilde{S} $.
When no confusion, we will say interchangeably objects defined on $ \overline{M} $ or $ \widetilde{M} $ as well as on $ \overline{S} $ or $ \widetilde{S} $. 
Here, $ \widetilde{M} $ and $ \widetilde{S} $ depend on our objects.

The vector bundle $ E$ is equipped with a filtration \eqref{eq:uucw},  where $ X$ is now replaced by $ \widetilde{M} $.
Then,  
\begin{align}\label{eq:aaet}
 D= {E}/{\overline{T^{*}  \widetilde{M} } \cdot E}.
\end{align}

In current situation, $ E$ is also equipped with a sub-filtration, 
\begin{align}\label{eq:rx1j}
 E \supset \pi^{*} \left(\overline{T^{*}  \widetilde{S} }\right)\cdot E  \supset \cdots \supset \pi^{*} \left( \Lambda^{m} \left(\overline{T^{*}  \widetilde{S} }\right)\right)\cdot E \supset 0 . 
\end{align}
Put 
\begin{align}\label{eq:dexb}
 D^{V}  = E/\pi^{*} \left(\overline{T^{*}  \widetilde{S} }\right)\cdot E.
\end{align}
Then, $ D^{V} $ is a smooth $ \mathbf{Z} $-graded vector bundle on $ \widetilde{M} $, so that 
\begin{align}\label{eq:opai}
 D^{V} /p^{*}  \left(\overline{T^{*} X} \right)\cdot D^{V} = D.
\end{align}
Put
\begin{align}\label{eq:fwkt}
  \mathcal{D}^{V} = \mathscr{O}^{\infty} \left(D^{V} \right).
\end{align}


The direct images $\pi_{*} \mathcal{D}^{V}, \pi_{*} \mathcal{E} $ are sheaves defined on $ \overline{S} $.\footnote{We use the extension of $ \left(E,A^{E\prime \prime } \right) $ as in Proposition \ref{prop:njjy}.}
If $ \widetilde{U} $ is an open subset of $ \widetilde{S} $, we have 
\begin{align}\label{eq:jkkk}
 &\pi_{*} \mathcal{D}^{V} \left(\widetilde{U} \right)= C^{\infty} \left(\widetilde{U} \times X,D^{V}_{|\widetilde{U} \times X}  \right),&\pi_{*} \mathcal{E} \left(\widetilde{U} \right) = C^{\infty} \left(\widetilde{U} \times X,E_{|\widetilde{U} \times X} \right).
\end{align}
Then, $ \pi_{*} \mathcal{D}^{V}$ is an $ \mathscr{O}^{\infty} $-module on $ \overline{S} $, and $ \pi_{*} \mathcal{E} $ is a $ \mathscr{D} $-module on $ \overline{S} $.

Moreover, $ A^{E\prime \prime }_{|\widetilde{U} \times X} $ is a differential operator acting on $ \pi_{*} \mathcal{E} \left(\widetilde{U} \right) $.
This way, we get an operator $ A^{\pi_{*}  \mathcal{E}\prime \prime } $ acting on $ \pi_{*} \mathcal{E}$, so that $ \left(\pi_{*} \mathcal{E}, A^{\pi_{*}  \mathcal{E}\prime \prime }\right) $ is an $ \mathscr{O} $-complex on $\overline{S}$, and defines an object in $ \widetilde{{\rm B}}^{{\rm t}} \left(\overline{S} \right) $.\footnote{Since $ \overline{S} $ is compact, the tameness condition is empty.
However, we still keep the superscript $ {\rm t}$ in the notation.}
More generally, $ \pi_{*} $ is a functor from $ {{\rm B}}^{{\rm b,t}} \left(\overline{M} \right) $ to $ \widetilde{{\rm B}}^{{\rm t}} \left(\overline{S} \right) $.
By \eqref{eq:dexb}, $ A^{\pi_{*} \mathcal{E} \prime\prime } $ induces an operator $ v_{0 }^{\pi_{*} \mathcal{D}^{V} } $ on $ \pi_{*}\mathcal{D}^{V} $, so that $ \left(\pi_{*}\mathcal{D}^{V},v_{0 }^{\pi_{*} \mathcal{D}^{V} } \right)$ is an $ \mathscr{O}^{\infty} $-complex on $ \overline{S} $, which is the diagonal of $ \left(\pi_{*} \mathcal{E},A^{\pi_{*} \mathcal{E} \prime \prime } \right) $ (cf. \eqref{eq:aiww}).



The main results of this section are the following.

\begin{thm}\label{thm:et3l}
  If $ \left(E,A^{E\prime\prime } \right)$ is an object in $ {\rm B^{b,t}_{c}}\left(\overline{M} \right)$, then the direct image $ \pi_{*}  \mathcal{E} $ has coherent cohomology on $ \overline{S}$.
\end{thm}
\begin{proof}
  The proof of our theorem will be given in Section \ref{s:p2}. 
\end{proof}

\begin{prop}\label{prop:wvkt}
 If $ \left(E,A^{E\prime\prime } \right)$ is an object in $ {\rm B^{b,t}_{c}}\left(\overline{M} \right)$, if $ \left(E_{\overline{S} }  ,A^{E_{\overline{S} }  \prime \prime } \right)$ is an object in $ {{\rm B}}^{\rm b,t} \left(\overline{S} \right)$ with associated $ \mathscr{O} $-complex $ \mathcal{E}_{\overline{S} }$, and if $ \phi :  \mathcal{E}_{\overline{S} }  \to \pi_{*} \mathcal{E}$ is a morphism in $\widetilde{ {\rm B}} ^{{\rm t}} \left(\overline{S} \right) $,  then $ \phi$ is a quasi-isomorphism if and only if it induces a quasi-isomorphism on the diagonals,
 \begin{align}\label{eq:vlaj}
  \phi_{0 } : \mathcal{D}_{\overline{S} } \to \pi_{*} \mathcal{D}^{V}.
 \end{align}
\end{prop}
\begin{proof}
 The proof of our proposition will be given in Section \ref{s:p10}. 
\end{proof}

\begin{thm}\label{thm:ete1}
  The following statements hold.
  \begin{enumerate}[\indent a)]
    \item\label{thm:etea} If $ \left(E,A^{E\prime\prime } \right)$ is an object in $ {\rm B^{b,t}_{c}}\left(\overline{M} \right)$, then there exist an object $ \left(E_{\overline{S} }  ,A^{E_{\overline{S} }  \prime \prime } \right)$ in $ {{\rm B}}^{\rm b,t} \left(\overline{S} \right)$ and a quasi-isomorphism in $\widetilde{ {\rm B}} ^{{\rm t}} \left(\overline{S} \right) $,  
    \begin{align}\label{eq:wtkt}
     \phi: \mathcal{E}_{\overline{S} }  \to \pi_{*} \mathcal{E}.
    \end{align}
    \item\label{thm:eteb} If $ \left(\left(\underline{E} ,A^{\underline{E} \prime \prime } \right),\left(\underline{E}_{\overline{S} }, A^{\underline{E}_{\overline{S} } \prime \prime } \right), \underline{\phi} \right)$ is another data as $ \left(\left(E,A^{E\prime \prime } \right),\left(E_{\overline{S} }, A^{E_{\overline{S} } \prime \prime } \right), \phi\right)$ in \ref{thm:etea}),  if $ \psi: \left(E,A^{E\prime \prime } \right) \to \left(\underline{E}, A^{\underline{E} \prime \prime } \right) $ is a morphism\footnote{We need take an extension satisfying analogue properties as Proposition \ref{prop:njjy}.} in $ {\rm B^{b,t}} \left(\overline{M} \right)$, then there exists a unique morphism $ \widetilde{\psi } $ in $ {\rm \underline{B} ^{b,t}}\left(\overline{S} \right) $ such that the  diagram
    \begin{equation}\label{dia:1qca2}
\begin{tikzsd}
          \mathcal{E}_{\overline{S} }   & \underline{{\mathcal{E} }}_{\overline{S} }    \\
          \pi_{*} \mathcal{E}   & \pi_{*}  \underline{{\mathcal{E} }}
          \arrow["\phi ", from=1-1, to=2-1]
          \arrow["\widetilde{\psi }", from=1-1, to=1-2,dotted]
          \arrow["\underline{\phi} ", from=1-2, to=2-2]
          \arrow["\pi_{*} \psi ", from=2-1, to=2-2]
        \end{tikzsd}
    \end{equation}
    in $ \widetilde{\underline{{\rm B}}}^{{\rm t}}  \left(\overline{S} \right)$ commutes.
  \end{enumerate}
\end{thm} 
\begin{proof}
  The proof of Statement \ref{thm:etea}) will be given in Section \ref{s:p1}. 
  We establish Statement \ref{thm:eteb}) by the same method used in Propositions \ref{prop:2gpa} 
  and \ref{prop:y1hu} (cf.~Footnote \ref{footnoteDFA}), while replacing  Proposition \ref{prop:ceuw} by Proposition \ref{prop:wvkt}.
\end{proof}

\subsection{Subcomplexes of $ \pi_{*} \protect\mathcal{D}^{V}  $ and $ \pi_{*} \protect\mathcal{E} $}\label{s:2sub}
\begin{defin}\label{def:zbcx} 
  Let $ \left(\pi_{*} \mathcal{D}^{V} \right)_{c} $ (resp.~$ \left(\pi_{*}  \mathcal{E}\right)_{c} $) be a subsheaf of $ \pi_{*} \mathcal{D}^{V}  $ (resp.~$ \pi_{*} \mathcal{E} $) on $ \overline{S} $, consisting of all the smooth sections $ u$ over an open subset $\widetilde{U}  \subset  \widetilde{S} $ such that for every compact subset $K \subset \widetilde{U}$, the intersection 
  \begin{align}\label{eq:ebsl}
   {\rm Supp}\left(u\right) \cap \pi^{-1} \left(K \right)
  \end{align}
  is compact.
\end{defin}
 
Then, $ \left(\pi_{*} \mathcal{D}^{V} \right)_{c} $ is an $ \mathscr{O}^{\infty} $-module on $ \overline{S} $, and $ \left(\pi_{*}  \mathcal{E}\right)_{c} $ is a $ \mathscr{D} $-module on $ \overline{S} $.
Moreover, $ v_{0 }^{\pi_{*} \mathcal{D}^{V}  } $ restricts on $ \left( \pi_{*} \mathcal{D}^{V} \right)_{c} $ and $ A^{\pi_{*} \mathcal{E}\prime \prime  } $ restricts on $ \left( \pi_{*} \mathcal{E} \right)_{c}$, 
so that $\left( \pi_{*} \mathcal{E} \right)_{c} $ is an object in $ \widetilde{{\rm B}}^{{\rm t}} \left(\overline{S} \right)$ with the diagonal $ \left(\pi_{*} \mathcal{D}^{V} \right)_{c} $  (cf. \eqref{eq:aiww}).

The inclusion defines a morphism of $ \mathscr{O}^{\infty} $-complex on $ \overline{S} $, 
\begin{align}\label{eq:g2el}
 \left( \pi_{*} \mathcal{D}^{V}  \right)_{c} \to \pi_{*} \mathcal{D}^{V},
\end{align}
as well as a morphism in $ \widetilde{{\rm B}}^{{\rm t}} \left(\overline{S} \right)$,
\begin{align}\label{eq:lrex}
\left( \pi_{*} \mathcal{E} \right)_{c} \to \pi_{*} \mathcal{E}.
\end{align}

\begin{prop}\label{prop:q3oq}
 If $ \left(E,A^{E\prime\prime } \right)$ is an object in $ {\rm B^{b,t}_{c}}\left(\overline{M} \right) $, then \eqref{eq:g2el} is a homotopy equivalence of $ \mathscr{O}^{\infty} $-complexes on $ \overline{S} $, and \eqref{eq:lrex} is a homotopy equivalence in $ \widetilde{{\rm B}}^{{\rm t}} \left(\overline{S} \right)$, i.e., an isomorphism in $ \widetilde{{\rm \underline{B}}}^{{\rm t}} \left(\overline{S} \right)$.
\end{prop}
\begin{proof}
  By Proposition \ref{prop:njjy}, $ \left(D,v_{0 } \right)$ is well-defined and constant $ \widetilde{S} \times U_{0 } $.
  Since $ \left(E,A^{E\prime\prime } \right)$ has compact support, $ \left(D,v_{0 } \right)$ is exact on $ \widetilde{S} \times U_{0 } $.
  Let $ v_{0 }^{\dagger}  $ be a constant section of $ \End^{-1} \left(D\right)$ on $ \widetilde{S}  \times U_{0 }  $, such that 
  \begin{align}\label{eq:ulkv}
   \left[v_{0 }, v_{0 }^{\dagger}\right]= 1.
  \end{align}
  By \eqref{eq:ulkv}, we have an identity on $ \widetilde{S} \times U_{0 } $, 
  \begin{align}\label{eq:slok}
   \left[A^{E\prime\prime},v_{0 }^{\dagger} \right]= 1.
  \end{align}
  We extend $ v_{0 }^{\dagger}$ to sections $ v_{0 }^{D,\dagger}, v_{0 }^{E,\dagger}$ of $ \End^{-1} \left(D\right)$, $ \End^{-1} \left(E\right)$ (cf. Footnote \ref{Footenote:Hom}) on $ \widetilde{S} \times X $.
    
  Put 
  \begin{align}\label{eq:1vnf}
  &\phi_{0} = 1- \left[v_{0},v_{0 }^{D,\dagger} \right],& \phi= 1- \left[A^{E\prime\prime },v_{0 }^{E,\dagger} \right].
  \end{align}
  By \eqref{eq:ulkv}-\eqref{eq:1vnf}, we see that $\phi_{0 },  \phi$ induce morphisms of complexes, 
  \begin{align}\label{eq:blus} 
    &\pi_{*} \mathcal{D}^{V}   \to \left(\pi_{*} \mathcal{D}^{V}  \right)_{c},&\pi_{*} \mathcal{E}  \to \left(\pi_{*} \mathcal{E} \right)_{c},
  \end{align}
  which are our desired homotopy inverses of \eqref{eq:g2el} and \eqref{eq:lrex}.
  The proof of our proposition is complete.
\end{proof}

Let $ T_{\mathbf{R} } X$ be the real tangent bundle of $ X$.
Let $ J^{T_{\mathbf{R} } X} $ be the complex structure on $ T_{\mathbf{R} } X$.
Let $ g^{T_{\mathbf{R} } X} $ be a $ J^{T_{\mathbf{R} } X} $-invariant complete Riemannian metric\footnote{ By \cite{NomizuOzeki61}, a complete Riemannian metric $ g^{T_{\mathbf{R} } X}_{0} $ exists always. Then, $ g^{T_{\mathbf{R} } X} \left(\cdot,\cdot\right)=g^{T_{\mathbf{R} } X}_{0 }  \left(\cdot,\cdot\right) + g^{T_{\mathbf{R} } X}_{0 }  \left(J^{T_{\mathbf{R} } X} \cdot,J^{T_{\mathbf{R} } X} \cdot\right)$ is $ J^{T_{\mathbf{R} } X} $-invariant and complete.
} on $ X$.
It induces a Hermitian metric $ g^{TX} $ on $ TX$.
Let $ g^{T \widetilde{S}} $ be a Hermitian metric on $ \widetilde{S}$.
Let\footnote{Most of the notation of this section is defined on $ U_{0 } $ or $ \widetilde{S} \times U_{0 } $ and will be extended non-canonically to $ X$ or $ \widetilde{S} \times X$ later.
To make the notation less cumbersome, sometimes we will drop the subscript $|U_{0}$ or $|\widetilde{S} \times U_{0}$.
} $ g^{D} $ be a constant metric on $ D_{|\widetilde{S} \times U_{0 } } $.

Let $ \nabla^{TX}, \nabla^{T \widetilde{S} }$ be the associated Chern connections on $ TX,T \widetilde{S}$.
Recall that $ \nabla^{D} $ is the constant connection on $ D_{ |\widetilde{S} \times U_{0 }} $. 
Then, $ \nabla^{TX}, \nabla^{T \widetilde{S} },\nabla^{D} $ induce a connection $ \nabla^{E} $ on $E_{|  \widetilde{S} \times U_{0 }} $.

If $ s\in \widetilde{S} $, and if  
\begin{align}\label{eq:gonl}
 U_{0, s} = \left\{s\right\}\times U_{0 },
\end{align}
we have canonical identifications, 
\begin{align}\begin{aligned}\label{eq:ijrr}
 C^{\infty} \left(U_{0, s},D^{V} _{|U_{0, s} }  \right)&= \Omega^{0,\bullet} \left( U_{0,s} ,D_{| U_{0,s} } \right),\\
 C^{\infty} \left(U_{0, s},E_{|U_{0, s} }  \right)&= \Lambda \left(\overline{T^{*}_{s}  \widetilde{S} } \right)\widehat{\otimes } \Omega^{0,\bullet} \left( U_{0,s} ,D_{| U_{0,s} } \right).
\end{aligned}\end{align}

Since $ D_{|U_{0,s} } $ is trivial, $ \Omega^{0,\bullet} \left( U_{0,s} ,D_{| U_{0,s} } \right)$ is independent of $ s$.
The metrics $ g^{TX}$, $ g^{D} $ induce an $ L^{2} $-norm on  $ C^{\infty} \left(U_{0, s},D^{V}_{|U_{0, s} }  \right)$, which is $ s$-independent.
Together with $ g^{T \widetilde{S} } $, using the second identity in \eqref{eq:ijrr}, we get an $ L^{2} $-norm on $ C^{\infty} \left(U_{0, s},E_{|U_{0, s} }  \right)$.
Both these $ L^{2} $-norms are denoted by $ \left\| \cdot \right\|_{L^{2} } $.

Moreover, the Dolbeault differential  $\overline{\partial}^{X} : \Omega^{0,\bullet} \left( U_{0,s} ,D_{| U_{0,s} } \right) \to  \Omega^{0,\bullet+1} \left( U_{0,s} ,D_{| U_{0,s} } \right)$
and its formal $ L^{2} $-adjoint $ \overline{\partial}^{X}_{*} $
do not depend on $ s$.

Put 
\begin{align}\label{eq:qcqh}
 &A_{0 } = 1,&B_{0 } =0.
\end{align}
If $ k\in \mathbf{N}^{*} $, let $ A_{k}, B_{k} $ be the differential operators on $ U_{0 } $ of order $ k$ defined by the following alternating $ k$-products of $ \overline{\partial}^{X}, \overline{\partial}^{X}_{*} $, 
\begin{align}\label{eq:awdl}
 &A_{k} = \overline{\partial}^{X}\overline{\partial}^{X}_{*}\overline{\partial}^{X} \cdots, &B_{k} = \overline{\partial}^{X}_{*} \overline{\partial}^{X} \overline{\partial}^{X}_{*}\cdots.
\end{align}
By \eqref{eq:qcqh} and \eqref{eq:awdl}, for $ k\in \mathbf{N} $, we have  
\begin{align}\label{eq:vcli}
 \left(\overline{\partial}^{X} +\overline{\partial}^{X}_{*}  \right)^{k}= A_{k} +B_{k}. 
\end{align}

For $ u\in \pi_{*} \mathcal{D}^{V}  \left(\widetilde{U} \right)$ or $ \pi_{*} \mathcal{E}  \left(\widetilde{U} \right)$, if $ s\in \widetilde{U} $, we have
\begin{align}\label{eq:x1pg}
 u_{|U_{0,s} } \in C^{\infty} \left(U_{0,s},D^{V} _{|U_{0,s} } \right) \text{ or }  C^{\infty} \left(U_{0,s},E_{|U_{0,s} } \right).
\end{align}
If $ \ell\in \mathbf{N} $, put
\begin{align}\label{eq:supp}
 \left\|  u _{|U_{0,s} } \right\|_{\ell}  = \sqrt{\sum_{k= 0 }^{\ell}   \left(\left\| A_{k}   u _{|U_{0,s} } \right\|_{L^{2} }^{2} +\left\| B_{k}   u _{|U_{0,s} } \right\|_{L^{2}}^{2}\right)} \in \left[0, \infty\right]. 
\end{align}
More generally, if  $ V$ is a smooth vector field on $ \widetilde{S}  $, then $ \nabla^{E}_{V} u$ is well-defined on $ \widetilde{U} \times U_{0} $, so that $ \left\|  \left(\nabla^{ E}_{V} u\right) _{|U_{0,s} } \right\|_{\ell}\in \left[0, \infty\right]$ is also well-defined.



\begin{defin}\label{def:2cjz}
  Let $ \left(\pi_{*}  \mathcal{D}^{V} \right)_{\left(2\right)}$ (resp.~$ \left(\pi_{*}  \mathcal{E}\right)_{\left(2\right)}$) be a subsheaf of $ \pi_{*} \mathcal{D}^{V}  $ (resp.~$ \pi_{*} \mathcal{E} $) on $ \overline{S} $, consisting of all the smooth sections $ u$ over an open subset $\widetilde{U}  \subset  \widetilde{S} $ such that for any compact subset $ K\subset \widetilde{U}$, any $ \ell,\ell^{\prime } \in \mathbf{N} $, and smooth vector fields $ V_{1},\ldots, V_{\ell^{\prime } } $ on $ \widetilde{S} $, 
  \begin{align}\label{eq:cjys}
 \sup_{s\in K}   \left\| \left(\nabla^{E} _{V_{1} }  \cdots \nabla^{E} _{V_{\ell^{\prime }  } }  u\right) _{|U_{0, s} }  \right\|_{\ell } < + \infty.
  \end{align}
\end{defin}

\begin{re}\label{re:f1ah}
 Note that $ \left(\pi_{*} \mathcal{D}^{V} \right)_{\left(2\right)},\left(\pi_{*}  \mathcal{E}\right)_{\left(2\right)}$ depend only on $ g^{TX}_{|U_{0 } }  $, but neither on $  g^{T \widetilde{S} }$ nor $g^{D}$.
\end{re}

As before, $ \left(\pi_{*} \mathcal{D}^{V} \right)_{\left(2\right)}$ is an $ \mathscr{O}^{\infty}  $-module on $ \overline{S} $, and $ \left(\pi_{*}  \mathcal{E}\right)_{\left(2\right)}$ is a $ \mathscr{D} $-module on $ \overline{S} $.
Moreover, $ v_{0 }^{\pi_{*} \mathcal{D}^{V} } $ restricts on $ \left( \pi_{*} \mathcal{D}^{V}  \right)_{\left(2\right)}$ and $ A^{\pi _{*} \mathcal{E} \prime \prime } $ restricts on $ \left( \pi_{*} \mathcal{E} \right)_{\left(2\right)}$, so that $\left( \pi_{*} \mathcal{E} \right)_{\left(2\right)} $ is an object in $ \widetilde{{\rm B}}^{{\rm t}} \left(\overline{S} \right)$ with diagonal $ \left(\pi_{*} \mathcal{D}^{V} \right)_{\left(2\right)}$ (cf. \eqref{eq:aiww}).

The inclusions define morphisms of $ \mathscr{O}^{\infty}  $-complexes on $ \overline{S} $, 
\begin{align}\label{eq:aoup}
 &  \left(\pi_{*} \mathcal{D}^{V}  \right)_{c} \to \left(\pi_{*} \mathcal{D}^{V}  \right)_{\left(2\right)}, &\left(\pi_{*} \mathcal{D}^{V}  \right)_{\left(2\right)}\to \pi_{*} \mathcal{D}^{V} ,
\end{align}
as well as morphisms in $ \widetilde{{\rm B}}^{{\rm t}} \left(\overline{S} \right)$, 
\begin{align}\label{eq:acxm}
&  \left(\pi_{*} \mathcal{E} \right)_{c} \to \left(\pi_{*} \mathcal{E} \right)_{\left(2\right)}, &\left(\pi_{*} \mathcal{E} \right)_{\left(2\right)}\to \pi_{*} \mathcal{E}.
\end{align}

\begin{prop}\label{prop:iu2w}
  If $ \left(E,A^{E\prime\prime } \right)$ is an object in ${\rm {B} ^{b,t}_{c}}\left(\overline{M} \right) $, then 
  \eqref{eq:aoup} are homotopy equivalences of $ \mathscr{O}^{\infty} $-complexes on $ \overline{S} $, and \eqref{eq:acxm} are homotopy equivalences in $ \widetilde{{\rm B}}^{{\rm t}} \left(\overline{S} \right)$, i.e., isomorphisms in $ \widetilde{{\rm \underline{B}}}^{{\rm t}} \left(\overline{S} \right)$.
\end{prop}
\begin{proof}  
 As in the proof of Proposition \ref{prop:q3oq}, the morphisms \eqref{eq:1vnf} induce the homotopic inverses of \eqref{eq:aoup} and \eqref{eq:acxm}. 
\end{proof}

\subsection{A tame splitting}\label{s:split}
In the rest part of this section, we fix a tame\footnote{Outside a compact subset of $ \overline{M} $, the splitting \eqref{eq:nipb} is the canonical identification provided by the tame structure.} splitting on $ \widetilde{M} = \widetilde{S} \times X $, 
\begin{align}\label{eq:nipb}
 E= \Lambda \left(\overline{T^{*} \widetilde{M} } \right)\widehat{\otimes} D.
\end{align}
By \eqref{eq:opai} and \eqref{eq:nipb}, we have 
\begin{align}\label{eq:3w2j}
& D^{V} = p^{*} \Lambda \left(\overline{T^{*} X } \right)\widehat{\otimes} D,&E =  \pi^{*} \Lambda \left(\overline{T^{*} \widetilde{S} } \right) \widehat{\otimes}    p^{*} \Lambda \left(\overline{T^{*} X } \right)\widehat{\otimes} D.
\end{align}

By \eqref{eq:nipb} and \eqref{eq:3w2j}, we have  
\begin{align}\label{eq:rong}
&\pi_{*} \mathcal{D}^{V} \left(\widetilde{S} \right)= C^{\infty} \left(\widetilde{S} \times X,p^{*} \Lambda \left(\overline{T^{*} X } \right)\widehat{\otimes} D  \right), &\pi_{*} \mathcal{E} \left(\widetilde{S} \right)= \Omega^{0,\bullet } \left(\widetilde{S} \times X, D\right).
\end{align}
Moreover, if $ s\in \widetilde{S} $, then 
\begin{align}\label{eq:bjw2}
&C^{\infty} \left(X_{s}, D^{V} _{|X_{s} } \right)=  \Omega^{0,\bullet } \left(X_{s} ,D_{|X_{s}  } \right),& C^{\infty} \left(X_{s}, E_{|X_{s} } \right)= \Lambda \left(\overline{T^{*}_{s}   \widetilde{S} } \right)\widehat{\otimes } \Omega^{0,\bullet } \left(X_{s} ,D_{|X_{s}  } \right).
\end{align}
Also, we have identifications of $ \mathscr{D} $-modules on $ \overline{S} $, 
\begin{align}\label{eq:aiww}
 &\pi_{*} \mathcal{E} = \mathscr{D} \widehat{\otimes} _{\mathscr{O}^{\infty} } \pi_{*} \mathcal{D}^{V},
 &\left(\pi_{*} \mathcal{E}\right)_{c}  = \mathscr{D} \widehat{\otimes} _{\mathscr{O}^{\infty} } \left(\pi_{*} \mathcal{D}^{V}\right)_{c}, &
 &\left(\pi_{*} \mathcal{E}\right)_{\left(2\right)}  = \mathscr{D} \widehat{\otimes} _{\mathscr{O}^{\infty} }\left( \pi_{*} \mathcal{D}^{V}\right)_{\left(2\right)}. 
\end{align}

Using \eqref{eq:nipb}, we can write $ A^{E\prime \prime } $ as in \eqref{eq:ktja}.
If  $ B^{E\prime \prime } = \sum_{i\ge 2}^{} v_{i}$,  we have 
\begin{align}\label{eq:u2hd}
 A^{E\prime \prime } = v_{0 } + \nabla^{D\prime \prime } + B^{E\prime \prime }. 
\end{align}
The operators $ \nabla^{D\prime \prime },B^{E\prime \prime }$ split into 
\begin{align}\label{eq:vij2} 
&\nabla^{D\prime \prime } =  \overline{\partial}^{X} + \nabla^{D,H\prime \prime },&B^{E\prime \prime } = B^{E,X\prime \prime } + B^{E,H\prime \prime },
\end{align}
where $ \overline{\partial}^{X},\nabla^{D,H\prime \prime }$ are respectively the vertical and horizontal partial connections along $ X$ and $ \widetilde{S}$, and $ B^{E,X\prime \prime }$ is the restriction of $ B^{E\prime \prime }$ to $ X$.
Note that $ B^{E,H\prime \prime }$ induces an operator on $ \pi_{*} \mathcal{E} $, which increases the degree of the filtration \eqref{eq:rx1j} by $ 1$.\footnote{The operator $ B^{E,H\prime \prime } $ increases the degree of the total filtration by $ 2$, while the subfiltration by $ 1$, because of terms with mixed nature.
}

\begin{prop}\label{prop:n32d}
 The following identities hold, 
 \begin{align}\label{eq:mab2}
&v^{\pi_{*}  \mathcal{D}^{V}  }_{0} = v_{0 } + \overline{\partial}^{X} + B^{E,X\prime \prime }, &A^{\pi_{*} \mathcal{E} \prime \prime } =  v_{0 }^{\pi_{*}  \mathcal{D}^{V}  }  + \nabla^{D,H\prime \prime }  + B^{E,H\prime \prime }. 
\end{align}
Restricting to $\widetilde{S}  \times  U_{0 } $, we have 
\begin{align}\label{eq:jeja}
 &v_{0 }^{\pi_{*}  \mathcal{D}^{V}  } = v_{0 }  + \overline{\partial}^{X},& A^{\pi_{*} \mathcal{E} \prime \prime } =  v_{0 } + \overline{\partial}^{X} + \nabla^{D,H\prime \prime }.
\end{align}
\end{prop} 
\begin{proof}
  Our proposition is a consequence of  \eqref{eq:u2hd} and \eqref{eq:vij2}.
\end{proof}

We extend the constant metric $ g^{D} $ on $ D_{|\widetilde{S} \times U_{0} }$ to a Hermitian metric on $ D$ over $ \widetilde{S} \times X$.
By \eqref{eq:bjw2}, $ C^{\infty} \left(X_{s}, D^{V}_{|X_ s} \right)$ and $ C^{\infty} \left(X_{s}, E_{|X_ s} \right)$ are equipped with the $ L^{2} $-norms induced by $ \left( g^{TX}, g^{D}\right)$ and by $ \left(g^{T \widetilde{ S} },  g^{TX}, g^{D}\right)$ respectively.
Both these norms will still be denoted by $ \left\| \cdot \right\|_{L^{2} } $.

By the first equation in \eqref{eq:mab2}, for each $ s\in \widetilde{S} $, $ v^{\pi_{*}  \mathcal{D}^{V}  }_{0}$ restricts to a first order differential operator $  v_{0,s}^{\pi_{*} \mathcal{D}^{V} } $ acting on $ C^{\infty} \left(X_{s}, D^{V}_{|X_ s} \right)$.

\begin{defin}\label{def:fas3}
  For $ s\in \widetilde{S}$,  let $ v_{0,*,s} ^{\pi_{*} \mathcal{D}^{V}} $ be the formal $ L^{2} $-adjoint of $ v_{0,s } ^{\pi_{*} \mathcal{D}^{V} } $.
\end{defin}
   
  Then, $ v_{0,*,s } ^{\pi_{*} \mathcal{D}^{V}}$ is a first order differential operator acting on $ C^{\infty}  \left(X_{s}, D^{V} _{|X_{s} } \right)$.
  If $ u\in \left(\pi_{*} \mathcal{D}^{V}  \right)\left(\widetilde{U} \right) $, then $ \left(s,x\right) \mapsto v_{0,*,s } ^{\pi_{*} \mathcal{D}^{V}} u\left(s,x\right)$  defines an element of $ \left(\pi_{*} \mathcal{D}^{V}  \right)\left(\widetilde{U} \right)$.
  This way we obtain an operator $ v_{0, *}^{\pi_{*} \mathcal{D}^{V}}$ acting on $ \pi_{*}  \mathcal{D}^{V}$.
  Denote by $ v_{0, *} $ and $ \overline{\partial}^{X}_{*} $ the operators defined in the same way.

Put 
\begin{align}\label{eq:wjav}
 L^{\pi_{*} \mathcal{D}^{V} }  = \left[v_{0 } ^{\pi_{*} \mathcal{D}^{V} },v_{0,*} ^{\pi_{*} \mathcal{D}^{V} }\right].
\end{align}

  \begin{prop}\label{prop:yemy}
    The following statements hold.
    \begin{enumerate}[\indent a)]
      \item\label{prop:yemya} The operators $ v_{0, *}^{\pi_{*} \mathcal{D}^{V}}, L^{\pi_{*} \mathcal{D}^{V} }$ are $ \mathscr{O}^{\infty}$-linear on $ \overline{S} $, such that 
   \begin{align}\label{eq:bwoe}
      \left[v_{0 }^{\pi_{*} \mathcal{D} ^{V} }, L^{\pi_{*} \mathcal{D}^{V} } \right]= 0.
    \end{align} 
      \item\label{prop:yemyb} Acting on $ \pi_{*} \mathcal{D}^{V}\left(\widetilde{S} \right)= C^{\infty} \left(\widetilde{S} \times X,D^{V} \right)$, $ v_{0, *}^{\pi_{*} \mathcal{D}^{V}}$ (resp.~$  L^{\pi_{*} \mathcal{D}^{V} }$) is a first (resp.~second) order differential operator with smooth coefficients on $ \widetilde{S} \times X$, whose differentiation only in the direction $ X$.   
      \item\label{prop:yemyc} Restricting to $ \widetilde{S}  \times U_{0 } $, we have 
      \begin{align}\label{eq:cecs}
&        v_{0,* }^{\pi_{*}  \mathcal{D}^{V}  }= v_{0, *} + \overline{\partial}^{X}_{*},&L^{\pi_{*} \mathcal{D}^{V} }  =  \left[\overline{\partial}^{X} , \overline{\partial}^{X}_{*}  \right] + \left[v_{0 }, v_{0, *} \right].
      \end{align}
      \item\label{prop:yemyd} The operators $ v_{0, *}^{\pi_{*} \mathcal{D}^{V}}, L^{\pi_{*} \mathcal{D}^{V} } $ preserve $ \left(\pi_{*} \mathcal{D}^{V} \right)_{\left(2\right)} $.  
    \end{enumerate}
  \end{prop}
  \begin{proof}
   Statements \ref{prop:yemya})-\ref{prop:yemyc}) follow immediately from the definition.
   Statement \ref{prop:yemyd}) is a consequence of \ref{prop:yemyc}).
  \end{proof}

  By \eqref{eq:aiww}, the operators $ v_{0, *}^{\pi_{*} \mathcal{D}^{V}}, v_{0, *}, \overline{\partial}^{X}_{*},  L^{\pi_{*} \mathcal{D}^{V} }  $ extend $ \mathscr{D} $-linearly to $ \pi_{*} \mathcal{E}$.
  Proposition \ref{prop:yemy} still holds for these extensions.
  Put 
  \begin{align}\label{eq:eqek}
    L^{\pi _{*} \mathcal{E} } = \left[A^{\pi_{*} \mathcal{E} \prime \prime }, v_{0,*} ^{\pi_{*} \mathcal{D}^{V} }\right].
  \end{align}

 \begin{prop}\label{prop:ywts}
  The following statements hold.
  \begin{enumerate}[\indent a)]
    \item\label{prop:ywtsa} The operator $ L^{\pi_{*} \mathcal{E} }$ is $ \mathscr{D}  $-linear on $ \overline{S} $, such that 
    \begin{align}\label{eq:bwoe1}
      \left[A^{\pi_{*} \mathcal{E} \prime \prime } , L^{\pi_{*} \mathcal{E} } \right]= 0.
    \end{align} 
    \item\label{prop:ywtsb} Acting on $ \pi_{*} \mathcal{E}\left(\widetilde{S} \right)= C^{\infty} \left(\widetilde{S} \times X,E \right)$, $ L^{\pi _{*} \mathcal{E} }$ is a second order differential operator with smooth coefficients on $ \widetilde{S} \times X$, whose differentiation only in the direction $ X$.
    \item\label{prop:ywtsc} Restricting to $ \widetilde{S}  \times U_{0 } $, we have  
    \begin{align}\label{eq:vhlp}
       L^{\pi_{*} \mathcal{E} } =  \left[\overline{\partial}^{X} , \overline{\partial}^{X}_{*}  \right] + \left[v_{0 }, v_{0, *} \right].
    \end{align}
    \item\label{prop:ywtsd}  The operator $ L^{\pi_{*} \mathcal{E} }$  preserves $ \left(\pi_{*} \mathcal{E} \right)_{\left(2\right)}$.
  \end{enumerate}
\end{prop}
\begin{proof}
  The proof of our proposition is the same as the one of Proposition \ref{prop:yemy}.
\end{proof}

By Propositions \ref{prop:yemy} \ref{prop:yemyb}) and \ref{prop:ywts} \ref{prop:ywtsb}), if $ s\in \widetilde{S} $, $ L^{\pi_{*} \mathcal{D}^{V}  } $ and $ L^{\pi_{*} \mathcal{E} } $ restrict to $ X_{s} $.

\begin{defin}\label{def:2ghl}
  Denote by $ L^{\pi_{*} \mathcal{D}^{V}  }_{s}  $ and $ L^{\pi_{*} \mathcal{E} }_{s}  $ the induced operators acting respectively on $C^{\infty} \left(X_{s}, D^{V} _{|X_{s} } \right)$ and $C^{\infty} \left(X_{s}, E_{|X_{s} } \right)$.
\end{defin}

\subsection{A trivialization over $ \widetilde{S} $}\label{s:triv}
In the sequel, we assume that $ \widetilde{S}$ is a small open ball of centre $ 0 $ in $ \mathbf{C}^{m} $.
We take $ g^{T \widetilde{S} } $ to be the restriction of the standard metric of $ \mathbf{C}^{m} $.

Then, 
\begin{align}\label{eq:nsyw}
 X_{0 } = \left\{0 \right\} \times X.
\end{align}
By \eqref{eq:3w2j}, we have 
\begin{align}\label{eq:rgt1}
&D_{|X_{0 }  }^{V} = \Lambda \left(\overline{T^{*} X_{0 } } \right)\widehat{\otimes } D_{|X_{0 } },& E_{|X_{0 }  } = \Lambda \left(\overline{\mathbf{C}^{m*} } \right)\widehat{\otimes } \Lambda \left(\overline{T^{*} X_{0 } } \right)\widehat{\otimes } D_{|X_{0 } }.
\end{align}

We extend the constant connection $ \nabla^{D} $ on $D_{| \widetilde{S} \times U_{0 }} $ to a Hermitian connection on $ D$ over $ \widetilde{S} \times X $ with respect to $ g^{D} $.
If $ s\in \widetilde{S} $, we identify $ D_{|X_{s} } $ with $ D_{|X_{0 }  } $ using the parallel transport along the line $ \left(ts,x\right)|_{0 \le t\le 1}$ with respect to $ \nabla^{D} $.
If $ q : \widetilde{S}   \times X \to X_{0 } $ is the obvious projection, then 
\begin{align}\label{eq:opwq}
& D=q^{*} D_{|X_{0 }  },&g^{D} = q^{*} g^{D} _{|X_{0 } } . 
\end{align}

By \eqref{eq:rgt1} and \eqref{eq:opwq}, we have 
\begin{align}\begin{aligned}\label{eq:h22p}
  \pi_{*} \mathcal{D}^{V} \left(\widetilde{S} \right)&= C^{\infty}  \left(\widetilde{S} \times X, q^{*}  \left(\Lambda \left(\overline{T^{*} X_{0 }}  \right)\otimes D_{|X_{0 }  } \right)\right),\\
 \pi_{*} \mathcal{E} \left(\widetilde{S} \right)&= \Lambda \left(\overline{\mathbf{C}^{m*} } \right)\widehat{\otimes}  C^{\infty}  \left(\widetilde{S} \times X, q^{*}  \left(\Lambda \left(\overline{T^{*} X_{0 }}  \right)\otimes D_{|X_{0 } } \right)\right).
 \end{aligned}\end{align}
As \eqref{eq:bjw2}, if $ s\in \widetilde{S}$, we have 
\begin{align}\label{eq:tpfa}
  &C^{\infty}\left(X_{s},D_{|X_{s} }^{V}  \right) = \Omega^{0,\bullet } \left(X_{0 },D_{|X_{0 } } \right),
  &C^{\infty} \left(X_{s} , E_{|X_{s} } \right)= \Lambda  \left(\overline{\mathbf{C}^{m*}  } \right)\widehat{\otimes } \Omega^{0,\bullet }  \left(X_{0 } ,  D_{|X_{0 } } \right).
\end{align}

\begin{re}\label{re:whjg}
  The two above spaces together with their $ L^{2} $-norm are independent of $ s\in \widetilde{S} $.
  In particular, $L^{\pi_{*} \mathcal{D}^{V}  }_{s},L^{\pi_{*} \mathcal{E} }_{s}$ act on $ s$-independent spaces. 
  Note also that $ \overline{\partial}^{X} $ depends on $ s\in \widetilde{S} $ thought $ \nabla^{D\prime \prime } $.
\end{re}

\subsection{Resolvents as a morphism of sheaves}\label{s:trunc}
Since $ g^{T_{\mathbf{R} } X } $ is complete, by \cite[Theorem 2.2]{Chernoff73}, the operator $ \overline{\partial}^{X_{0 } } +\overline{\partial}^{X_{0 } }_{*} $ as well as all its higher powers are essentially self-adjoint.
For $ \ell\in \mathbf{N} $, denote by $ {\rm Dom}\left(\overline{\partial}^{X_{0 } } +\overline{\partial}^{X_{0 } }_{*}\right)^{\ell}$ the domain of the unique self-adjoint extension.

\begin{defin}\label{def:h31f}
  If $ \ell\in \mathbf{N} $, put 
   \begin{align}\label{eq:x2xk}
  &H^{\ell} \left(X_{0 }, D^{V} _{|X_{0 } } \right)= {\rm Dom}\left(\overline{\partial}^{X_{0 } } +\overline{\partial}^{X_{0 } }_{*}\right)^{\ell},&H^{\ell} \left(X_{0 }, E _{|X_{0 } } \right)= \Lambda \left(\overline{\mathbf{C}^{m*} } \right)\widehat{\otimes} H^{\ell} \left(X_{0 }, D^{V} _{|X_{0 } } \right).
 \end{align}
 Moreover, 
 \begin{align}\label{eq:mllm}
&H^{\infty} \left(X_{0 }, D^{V} _{|X_{0 } } \right)= \bigcap _{\ell\in   \mathbf{N} }H^{\ell} \left(X_{0 }, D^{V} _{|X_{0 } } \right) ,&H^{\infty} \left(X_{0 }, E _{|X_{0 } } \right)= \bigcap _{\ell\in   \mathbf{N} }H^{\ell} \left(X_{0 }, E _{|X_{0 } } \right). 
\end{align} 
 If $ u\in L^{2}  \left(X_{0 }, D^{V} _{|X_{0 } } \right)$ or $ L^{2}  \left(X_{0 }, E _{|X_{0 } } \right)$, put 
  \begin{align}\label{eq:xgjd}
  \left\| u \right\|_{H^{\ell} } = \sqrt{\sum_{k=0 } ^{\ell} \left\| \left(\overline{\partial}^{X_{0 } } +\overline{\partial}^{X_{0 } }_{*}\right)^{k}  u \right\|_{L^{2} }^{2}}\in [0, +\infty].
 \end{align}
\end{defin}

\begin{prop}\label{prop:qooz}
 If $ \ell\in \mathbf{N} $, for $ u\in C^{\infty} \left(X_{0 }, D^{V} _{|X_{0 } } \right)$ or $ C^{\infty} \left(X_{0 } ,E_{|X_{0 } } \right)$ such that $ {\rm Supp}\, u \subset U_{0 } $, we have the identity in $ \left[0,+ \infty\right]$, 
 \begin{align}\label{eq:2yud}
  \left\| u \right\|_{\ell} = \left\| u \right\|_{H^{\ell} }.
  \end{align}
\end{prop}
\begin{proof}
 Consider the vector bundle $ D^{V} _{|X_{0 } } $.
 The other case is the same. 

 We claim that 
  \begin{align}\label{eq:nfzq}
  \left\{v\in L^{2}\left(X_{0 }, D^{V} _{|X_{0 } } \right) : \left(\overline{\partial}^{X_{0 } } +\overline{\partial}^{X_{0 } }_{*}\right)v\in L^{2} \right\}= \left\{v\in L^{2}\left(X_{0 }, D^{V} _{|X_{0 } } \right) : \overline{\partial}^{X_{0 } }v, \overline{\partial}^{X_{0 } }_{*}v\in L^{2} \right\}.
 \end{align}
Moreover, if $ v^{\prime } $ is in the above set such that $ {\rm Supp}\, v^{\prime }  \subset  U_{0 } $, we have 
\begin{align}\label{eq:3kam}
 \left\| \left(\overline{\partial}^{X_{0 } } +\overline{\partial}^{X_{0 } }_{*}\right)v^{\prime }  \right\|^{2}_{L^{2} } = \left\| \overline{\partial}^{X_{0 } }v^{\prime }  \right\|^{2}_{L^{2} } + \left\| \overline{\partial}^{X_{0 } }_{*}v^{\prime }  \right\|^{2}_{L^{2} }.
\end{align} 

Indeed, the inverse inclusion $ \supset $ in \eqref{eq:nfzq} is trivial. 
Assume that $ v$ is contained in the left-hand side of \eqref{eq:nfzq}.
Then $ v$ is in the domain of the maximal extension of $ \overline{\partial}^{X_{0 } } +\overline{\partial}^{X_{0 } }_{*}$.
Since $ \overline{\partial}^{X_{0 } } +\overline{\partial}^{X_{0 } }_{*} $ is essentially self-adjoint, $ v$  is also in the domain of the minimal extension of $ \overline{\partial}^{X_{0 } } +\overline{\partial}^{X_{0 } }_{*}$. 
There exists a sequence $ \left(v_{k}\right)_{k\in \mathbf{N} } $ in $ C^{\infty}_{c} \left(X_{0 }, D^{V} _{|X_{0 } } \right)$ such that as $ k \to +\infty$, we have the $ L^{2} $-convergences,  
\begin{align}\label{eq:vngf}
 &v_{k} \to v, &\left(\overline{\partial}^{X_{0 } } +\overline{\partial}^{X_{0 } }_{*}\right)v_{k} \to \left(\overline{\partial}^{X_{0 } } +\overline{\partial}^{X_{0 } }_{*}\right)v.
\end{align}
Recall that $ \phi_{U_{0 } } $ is defined before \eqref{eq:klag}.
By \eqref{eq:vngf}, as $ k \to +\infty$, we have the $ L^{2} $-convergence, 
\begin{align}\label{eq:mpi2}
 \left(\overline{\partial}^{X_{0 } } +\overline{\partial}^{X_{0 } }_{*}\right)\left(\phi_{U_{0 } } v_{k} \right) \to \left(\overline{\partial}^{X_{0 } } +\overline{\partial}^{X_{0 } }_{*}\right)\left(\phi_{U_{0 } } v\right).
\end{align}
Observe that acting on smooth sections with compact support in $ U_{0 } $, the image of $ \overline{\partial}^{X_{0 } } $  and $ \overline{\partial}^{X_{0 } }_{*}$ is orthogonal.
By \eqref{eq:mpi2}, $ \overline{\partial}^{X_{0 } }\left(\phi_{U_{0 } } v_{k}\right) $ and $ \overline{\partial}^{X_{0 } }_{*}\left(\phi_{U_{0 } } v_{k} \right)$ are $ L^{2} $-Cauchy sequences.
Their distributional limits $ \overline{\partial}^{X_{0 } }\left(\phi_{U_{0 } } v\right)$ and $ \overline{\partial}^{X_{0 } }_{*}\left(\phi_{U_{0 } } v\right)$ are therefore $ L^{2} $, so that as $ k \to +\infty$, we have the $ L^{2} $-convergences,
\begin{align}\label{eq:dxad}
 &\overline{\partial}^{X_{0 } }\left(\phi_{U_{0 } } v_{k}\right)  \to \overline{\partial}^{X_{0 } }\left(\phi_{U_{0 } } v\right),&\overline{\partial}^{X_{0 } }_{*} \left(\phi_{U_{0 } } v_{k}\right)  \to \overline{\partial}_{*} ^{X_{0 } }\left(\phi_{U_{0 } } v\right).
\end{align}
In particular, $ \phi_{U_{0 } } v$ is contained in the right-hand side of \eqref{eq:nfzq}. 
Classical elliptic estimate implies that $ \left(1-\phi_{U_{0 } } \right)v$ is also contained in the right-hand side of \eqref{eq:nfzq}.
This gives us the inclusion $ \subset $ in \eqref{eq:nfzq}.

Since $ v^{\prime } $ is in the set \eqref{eq:nfzq}, we have 
\begin{align}\label{eq:jemh}
 \left\| \left(\overline{\partial}^{X_{0 } } +\overline{\partial}^{X_{0 } }_{*}\right)v^{\prime }  \right\|^{2}_{L^{2} } = \left\| \overline{\partial}^{X_{0 } }v^{\prime }  \right\|^{2}_{L^{2} } + \left\| \overline{\partial}^{X_{0 } }_{*}v^{\prime }  \right\|^{2}_{L^{2} }+ 2 {\rm Re}\left\langle \overline{\partial}^{X_{0 } }v^{\prime } , \overline{\partial}^{X_{0 } }_{*}v^{\prime } \right\rangle.
\end{align}
Since ${\rm Supp}\, v^{\prime } \subset U_{0} $, modifying the function $ \phi_{U_{0 } } $ if necessary, we have $ v^{\prime } = \phi_{U_{0 } } v^{\prime } $.
Using \eqref{eq:dxad} while replacing $ v$ by $ v^{\prime } $, we get
\begin{align}\label{eq:hqgw}
 \left\langle \overline{\partial}^{X_{0 } } v^{\prime } , \overline{\partial}^{X_{0 } }_{*}v^{\prime } \right\rangle = 0. 
\end{align}
By \eqref{eq:jemh} and \eqref{eq:hqgw}, we obtain \eqref{eq:3kam}.

Let us prove our proposition.
When $ \ell = 0$, \eqref{eq:2yud} holds trivially.
When $ \ell= 1$, by \eqref{eq:nfzq}, if one side of \eqref{eq:2yud} is finite, then the other side is also finite. 
In this case, \eqref{eq:2yud} is a consequence of \eqref{eq:3kam}.
When $ \ell \ge 2 $, applying repeatedly \eqref{eq:nfzq} and \eqref{eq:3kam}, we obtain \eqref{eq:2yud} and finish the proof of our proposition.
\end{proof}

\begin{re}\label{re:i2b3}
 By Proposition \ref{prop:qooz}, in defining $ \left(\pi_{*} \mathcal{D}^{V}\right)_{(2)} $ or $ \left(\pi_{*} \mathcal{E}\right)_{(2)} $, one may replace $ \left\| \cdot \right\|_{\ell} $  with $ \left\| \cdot \right\| _ {H^{\ell}}$.
 This definition has the advantage that $ \left\| \cdot \right\| _ {H^{\ell}}$ is taken over the whole $X$ rather than $U_{0}$, yielding better functional-analytic properties.
 However, it relies on the non-canonical splitting \eqref{eq:nipb}.
\end{re}

On $ \widetilde{S}  \times U_{0 } $, $ \left[v_{0 }, v_{0, *} \right]$ is a  positive-definite self-adjoint constant matrix.
Let $ a_{0 }> 0  $ be the smallest eigenvalue of this matrix.
Recall that $ \phi_{U_{0} } $ is defined before \eqref{eq:klag}.

\begin{defin}\label{def:hrci}
  Set 
  \begin{align}\label{eq:tvcd}
   M_{0 } = \left(\overline{\partial}^{X_{0 } } +  \overline{\partial}^{X_{0 } }_{*}  \right)^{2}  + \phi_{U_{0 } } \left[v_{0 }, v_{0,*} \right] + \left(1-\phi_{U_{0 } } \right)a_{0 }. 
  \end{align}   
\end{defin}

As $ \left(\overline{\partial}^{X_{0 } } +  \overline{\partial}^{X_{0 } }_{*}  \right)^{2} $, the operator $ M_{0 } $ is essentially self-adjoint.
The domain of its unique self-adjoint extension is $ H^{2} \left(X_{0 }, D^{V} _{|X_{0 } }\right)$.
Also, 
\begin{align}\label{eq:wpzw} 
 {\rm Sp}\left(M_{0 } \right) \subset \left [a_{0 }, \infty\right ).
\end{align}

\begin{prop}\label{prop:ul1a}
 The following statements hold.
 \begin{enumerate}[\indent a)]
   \item\label{prop:ul1aa} For $ k,\ell\in \mathbf{N} $, if $ R_{k} $ is a differential operator of order $ k$ whose coefficients have compact support, then $ R_{k} $ is bounded from $ H^{\ell+k} \left(X_{0 }, D^{V}_{|X_{0 } } \right)$ to $ H^{\ell} \left(X_{0 }, D^{V}_{|X_{0 } } \right)$.
   \item\label{prop:ul1ab} If $ \lambda \notin {\rm Sp}\left(M_{0 } \right)$, then $ \left(M_{0 } -\lambda \right)^{-1} $ is bounded from $ H^{\ell} \left(X_{0 }, D^{V}_{|X_{0 } } \right)$ to $ H^{\ell+2} \left(X_{0 }, D^{V}_{|X_{0 } } \right)$.
   \item\label{prop:ul1ac} The operator $ R_{1} \left(M_{0 } -\lambda \right)^{-1}  $ on $ L^{2} \left(X_{0 }, D^{V}_{|X_{0 } }\right)$ is compact.
   \item\label{prop:ul1abd} If $ N > 2n+1$, the operator $ \left(R_{1} \left(M_{0 } -\lambda \right)^{-1}\right)^{N} $ on $ L^{2} \left(X_{0 }, D^{V}_{|X_{0 } }\right)$ is of trace class.
 \end{enumerate} 
\end{prop}
\begin{proof}
 Statement \ref{prop:ul1aa}) is a consequence of the classical (local) elliptic estimates.
 
 Assume $ \lambda \notin {\rm Sp}\left(M_{0 } \right)$. 
 We claim that when acting on $ H^{1} \left(X_{0 }, D^{V}_{|X_{0 } } \right)$, we have
\begin{multline}\label{eq:cko1}
 \left(\overline{\partial} ^{X_{0 }  } + \overline{\partial} ^{X_{0 } }_{*}\right)\left(M_{0 } -\lambda \right)^{-1} = \left(M_{0 } -\lambda \right)^{-1}\left(\overline{\partial} ^{X_{0 }  }  + \overline{\partial} ^{X_{0 } }_{*}\right)\\
 -\left(M_{0 } -\lambda \right)^{-1}\left[\overline{\partial} ^{X_{0 }  } + \overline{\partial} ^{X_{0 } }_{*},M_{0 } \right]\left(M_{0 } -\lambda \right)^{-1}.
\end{multline}
Indeed, when acting on distributions, we have
\begin{align}\label{eq:1zv2}
 \left(M_{0 } -\lambda \right)\left(\overline{\partial} ^{X_{0 }  } + \overline{\partial} ^{X_{0 } }_{*}\right) = \left(\overline{\partial} ^{X_{0 }  }  + \overline{\partial} ^{X_{0 } }_{*}\right)\left(M_{0 } -\lambda \right)
 -\left[\overline{\partial} ^{X_{0 }  } + \overline{\partial} ^{X_{0 } }_{*},M_{0 } \right].
\end{align}
Composing \eqref{eq:1zv2} on the right by $ \left(M_{0 } -\lambda \right)^{-1} $, when acting on $ L^{2}\left(X_{0 }, D^{V}_{|X_{0 } } \right) $, we have
\begin{align}\label{eq:1zv3}
 \left(M_{0 } -\lambda \right)\left(\overline{\partial} ^{X_{0 }  } + \overline{\partial} ^{X_{0 } }_{*}\right)\left(M_{0 } -\lambda \right)^{-1}  = \left(\overline{\partial} ^{X_{0 }  }  + \overline{\partial} ^{X_{0 } }_{*}\right) -\left[\overline{\partial} ^{X_{0 }  } + \overline{\partial} ^{X_{0 } }_{*},M_{0 } \right]\left(M_{0 } -\lambda \right)^{-1}.
\end{align}
By \eqref{eq:tvcd}, $ \left[\overline{\partial} ^{X_{0 }  } + \overline{\partial} ^{X_{0 } }_{*},M_{0 } \right]$ is an operator of type $ R_{1} $.
By Statement \ref{prop:ul1aa}), when acting on $ H^{1} \left(X_{0 }, D^{V}_{|X_{0 } } \right)$, the image of the right-hand side of \eqref{eq:1zv3} is contained in $ L^{2} \left(X_{0 }, D^{V}_{|X_{0 } } \right)$.
By \eqref{eq:1zv3}, when acting on $ H^{1} \left(X_{0 }, D^{V}_{|X_{0 } } \right)$, the image of $ \left(\overline{\partial} ^{X_{0 }  } + \overline{\partial} ^{X_{0 } }_{*}\right)\left(M_{0 } -\lambda \right)^{-1}$ is contained in the domain of $ M_{0 } -\lambda$.
Applying $ \left(M_{0 } -\lambda \right)^{-1} $ on the left of \eqref{eq:1zv3}, we get \eqref{eq:cko1}.

Let us prove Statement \ref{prop:ul1ab}) by induction on $ \ell\in \mathbf{N} $.
If $ \ell = 0 $, it is a consequence of the definition of the domain of $ M_{0 } $.
If $ \ell \ge 1$, it is a consequence of the fact that $ \left[\overline{\partial} ^{X_{0 }  } + \overline{\partial} ^{X_{0 } }_{*},M_{0 } \right]$ is type $ R_{1} $, of Statement \ref{prop:ul1aa}), of \eqref{eq:cko1}, and of induction assumptions.

Statement \ref{prop:ul1ac}) is a consequence of Statements \ref{prop:ul1aa}), \ref{prop:ul1ab}), and of classical (local) Sobolev embedding theorem.

Let $ \varphi,\widetilde{\varphi} \in C^{\infty}_{c} \left(X_{0 } ,\mathbf{R} \right)$ be such that $ \varphi =1$ on the support of the coefficient of $ R_{1} $ and $\widetilde{\varphi} = 1$ on the support of $ \varphi$. 
For $ N \ge 1$, we can write
 \begin{align}\label{eq:okdf}
  \left(R_{1} \left(M_{0 } -\lambda \right)^{-1} \right)^{N}= \varphi  \widetilde{\varphi}\left(R_{1} \left(M_{0 } -\lambda \right)^{-1} \right)^{N-1} \widetilde{\varphi} \varphi R_{1} \left(M_{0 } -\lambda \right)^{-1}.
 \end{align}
The operator $ \widetilde{\varphi}  \left(R_{1} \left(M_{0 } -\lambda \right)^{-1} \right)^{N-1}\widetilde{\varphi} $ can be considered as an operator defined on some compact manifold $ X_{c} $, sending the corresponding $ L^{2} $-space to $ H^{N-1} $-space, while $ \varphi$ can be considered as bounded operator from the $ L^{2} $-spaces of $ X_{0 } $ to the one of $ X_{c} $, or vice versa.
If $ N> 2n+1$, since the embedding from $ H^{N-1} $ to $ L^{2} $ on $ X_{c} $ is of trace class, we see that $ \widetilde{\varphi} \left(R_{1} \left(M_{0 } -\lambda \right)^{-1} \right)^{N-1}\widetilde{\varphi} $ is of trace class.
By \eqref{eq:okdf}, we get Statement \ref{prop:ul1abd}).
The proof of our proposition is complete.
\end{proof}

Operators $L_{s}^{\pi_{*}\mathcal{D}^{V}  },  M_{0 }$ extend to $ C^{\infty} \left(X_{0 },E_{|X_{0 } } \right)$ in an obvious way.
In the sequel, we will not distinguish $L_{s}^{\pi_{*}\mathcal{D}^{V}  },  M_{0 }$ with the corresponding unbounded operators with domain $ H^{2} \left(X_{0 }, E _{|X_{0 } } \right)$.

\begin{prop}\label{prop:sthv}
 If $ s\in \widetilde{S} $, the operators $L_{s}^{\pi_{*} \mathcal{D}^{V}  }, L_{s}^{\pi_{*} \mathcal{E} }$ have the same spectrum, which is contained in $ \mathbf{R}_{+} $ and is discrete in $ \left [0, a_{0 } \right )$.
 If $ a\in \left(0, a_{0 } \right)$, if $ \widetilde{S} _{a} \subset \widetilde{S} $ is the set of points $ s\in \widetilde{S} $ such that $ a \not\in {\rm Sp}\left( L_{s}^{\pi_{*} \mathcal{D}^{V}  } \right)$, then $ \widetilde{S} _{a} $ is an open subset of $ \widetilde{S} $.
 In particular, $ \left(\widetilde{S}_{a}\right)_{a\in (0, a_{0 } )}  $ forms an open cover of $ \widetilde{S} $.
\end{prop}
\begin{proof}
Fix $ s\in \widetilde{S} $.  
If $ \lambda \notin {\rm Sp}\left(L_{s}^{\pi_{*} \mathcal{D}^{V}  } \right)$, we have 
\begin{align}\label{eq:biil}
 \lambda -L_{s}^{\pi_{*} \mathcal{E} } = \left(1-\left(L_{s}^{\pi_{*} \mathcal{E} } - L_{s}^{\pi_{*} \mathcal{D}^{V}  }\right)\left(\lambda - L_{s}^{\pi_{*} \mathcal{D}^{V}  }\right)^{-1} \right)\left(\lambda - L_{s}^{\pi_{*} \mathcal{D}^{V}  }\right). 
\end{align}
By \eqref{eq:cecs}, \eqref{eq:vhlp}, and by Proposition \ref{prop:ul1a} \ref{prop:ul1aa}), $ \left(L_{s}^{\pi_{*} \mathcal{E} } - L_{s}^{\pi_{*} \mathcal{D}^{V}  }\right)\left(\lambda - L_{s}^{\pi_{*} \mathcal{D}^{V}  }\right)^{-1} $ is bounded on the $ L^{2} $-space.
Since it is also nilpotent, we deduce that $ \lambda \not\in {\rm Sp}\left(L_{s}^{\pi_{*} \mathcal{E} } \right)$ and 
\begin{align}\label{eq:j3kq}
  \left(\lambda -L_{s}^{\pi_{*} \mathcal{E} }\right)^{-1} = \left(\lambda - L_{s}^{\pi_{*} \mathcal{D}^{V}  }\right)^{-1} \sum_{i= 0 }^{m} \left\{\left(L_{s}^{\pi_{*} \mathcal{E} } - L_{s}^{\pi_{*} \mathcal{D}^{V}  }\right)\left(\lambda - L_{s}^{\pi_{*} \mathcal{D}^{V}  }\right)^{-1}\right\}^{i}. 
 \end{align}
 In particular, $ {\rm Sp}\left(L_{s}^{\pi_{*} \mathcal{E} } \right) \subset {\rm Sp}\left(L_{s}^{\pi_{*} \mathcal{D}^{V}  } \right). $ 
If we exchange the role of $ L^{\pi_{*} \mathcal{D}^{V}  }_{s}  $ and $ L^{\pi_{*} \mathcal{E} }_{s} $ in the previous argument, we get the reverse inclusion.
Thus, 
\begin{align}\label{eq:gi3q}
  {\rm Sp}\left(L_{s}^{\pi_{*} \mathcal{E} } \right) = {\rm Sp}\left(L_{s}^{\pi_{*} \mathcal{D}^{V}  } \right). 
\end{align}

By \eqref{eq:wjav}, $ L_{s}^{\pi_{*}  \mathcal{D}^{V}  }$ with domain $ H^{2} \left(X_{0} , E _{|X_{0 } } \right)$ is self-adjoint and non-negative.
Thus, 
\begin{align}\label{eq:mfzt}
  {\rm Sp}\left(L_{s}^{\pi_{*}  \mathcal{D}^{V}  }\right) \subset \mathbf{R}_{+}. 
\end{align}

Put 
\begin{align}\label{eq:h2u2}
 R_{s} = L_{s}^{\pi_{*} \mathcal{D}^{V}  }   - M_{0 }. 
\end{align}
Since the principal symbols of $ L_{s}^{\pi_{*} \mathcal{D}^{V}  }   $ and $ M_{0 }$ coincide, $ R_{s} $ is a differential operator of order $ 1 $.
By \eqref{eq:cecs} and \eqref{eq:tvcd}, $ R_{s} $ is of type $ R_{1} $ as in Proposition \ref{prop:ul1a} \ref{prop:ul1aa}).

In the rest part of the proof, we assume $ \lambda \in \mathbf{C} $ and $ {\rm Re}\left(\lambda \right)<a_{0 } $.
Note that we have an isomorphism of Hilbert spaces,  
 \begin{align}\label{eq:qtzu}
  M_{0 }  -\lambda : H^{2}\left(X_{0 },  D^{V} _{|X_{0 }  } \right) \to L^{2}\left(X_{0 }, D^{V} _{|X_{0 }  } \right). 
 \end{align}
Put 
 \begin{align}\label{eq:j2my}
  K_{s}\left(\lambda \right) = R_{s} \left(M_{0 }  -\lambda\right)^{-1}.
 \end{align}
 By Proposition \ref{prop:ul1a} \ref{prop:ul1ac}), $ K_{s}\left(\lambda \right) $ is a compact operator.
 By \eqref{eq:h2u2} and \eqref{eq:j2my}, we have  
 \begin{align}\label{eq:fh1n}
  L_{s}^{\pi_{*}  \mathcal{D}^{V}  } -\lambda = \left(1+K_{s}\left(\lambda \right)  \right)\left(M_{0 } -\lambda\right). 
 \end{align}

 By \eqref{eq:qtzu} and \eqref{eq:fh1n}, we see that $ \lambda\notin {\rm Sp}\left(L_{s}^{\pi_{*}  \mathcal{D}^{V}  }\right)$ if and only if $1+ K_{s} \left(\lambda \right)$ is invertible on the $ L^{2} $-space.
 By \eqref{eq:mfzt}, if $ \lambda = -1$, $1+ K_{s} \left(\lambda \right)$ is invertible.
 On the set $ \left\{\lambda\in \mathbf{C} :{\rm Re}\lambda < a_{0 } \right\}$, $1+ K_{s} \left(\lambda \right)$ is a holomorphic family of Fredholm operators of index $ 0 $.
 By the analytic theory of Fredholm operators \cite[Appendix D.3]{Zworski_SemiClassical}, $ 1+ K_{s} \left(\lambda \right) $ is invertible for all $ \lambda $ except at a discrete subset of $ \left\{\lambda\in \mathbf{C} :{\rm Re}\lambda < a_{0 } \right\}$.
 Combining this with \eqref{eq:mfzt}, we see that $ {\rm Sp}\left(L_{s}^{\pi_{*}  \mathcal{D}^{V}  }\right) \cap [0, a_{0 })$ is discrete in $ [0,a_{0})$.
 Moreover, $ \ker \left(1+ K_{s} \left(\lambda \right)\right)$ has finite dimension.
 The spectrum $ \lambda \in {\rm Sp}\left(L_{s}^{\pi_{*} \mathcal{D}^{V}  }\right) \cap [0, a_{0 })$ has finite multiplicities.
 This shows the spectrum of $ L_{s}^{\pi_{*}  \mathcal{D}^{V}  }$ in $ [0, a_{0 } )$ is discrete.

 When $ s$ varies, as bounded operators on $ L^{2} $-space, $ K_{s} \left(\lambda \right)$ varies continuously.
 Therefore, $ \widetilde{S} _{a} $ is an open subset of $ \widetilde{S} $, and forms an open cover of $ \widetilde{S} $ when $ a$ varies in $ (0, a_{0 } )$.
 The proof of our proposition is complete.
\end{proof}

Let $ a\in (0, a_{0 } ) $ such that $ \widetilde{S}_{a} \neq \varnothing $.

\begin{prop}\label{prop:cflb}
  If $ \lambda\in \mathbf{C} $ with $ \left| \lambda  \right| = a$, if $ \widetilde{U}_{a}  \subset \widetilde{S}_{a} $ is an open subset, the pointwise defined operators $ \left(\lambda -L^{\pi_{*} \mathcal{D}^{V}  }_{s}  \right)^{-1}$ and $ \left(\lambda -L^{\pi_{*} \mathcal{E} }_{s}  \right)^{-1}$, $ s\in \widetilde{U}_{a} $, extend respectively to  $ \left(\pi_{*} \mathcal{D}^{V}  \right)_{\left(2\right)}\left(\widetilde{U}_{a} \right) $ and $ \left(\pi_{*} \mathcal{E} \right)_{\left(2\right)} \left(  \widetilde{U} _{a} \right) $ as $ \lambda $-uniformly bounded operators with respect to the semi-norms \eqref{eq:cjys}.
  In particular, we get the operators $ \left(\lambda -L^{\pi_{*} \mathcal{D}^{V}  } \right)^{-1}$ and $ \left(\lambda -L^{\pi_{*} \mathcal{E} } \right)^{-1}$ acting on $ \left(\pi_{*} \mathcal{D}^{V}  \right)_{\left(2\right)|\widetilde{S}_{a} } $ and $ \left(\pi_{*} \mathcal{E} \right)_{\left(2\right)|\widetilde{S}_{a} } $.
\end{prop}
\begin{proof}
Let us show our proposition for $ L^{\pi_{*} {\mathcal{D}^{V}  }  } $.
The other case is the same.

If $ u\in \left(\pi_{*} \mathcal{D}^{V} \right)_{\left(2\right)}\left( \widetilde{U}_{a} \right)$, and if $ u_{s} = u_{|X_{s} } $, by \eqref{eq:fh1n}, we have  
\begin{align}\label{eq:3rny}
   \left(\lambda -L_{s}^{\pi_{*} \mathcal{D}^{V}  } \right)^{-1} u_{s} = \left(M_{0 } -\lambda\right)^{-1}  \left(1+K_{s}\left(\lambda \right)  \right)^{-1} u_{s}. 
\end{align}
We need to show that $ \left(\lambda -L_{s}^{\pi_{*} \mathcal{D}^{V}  } \right)^{-1} u_{s}$ is jointly smooth on $ \widetilde{U}_{a} \times X_{0 }$, and that for any multi-index $ \alpha\in \mathbf{N}^{m} $, $ \ell\in \mathbf{N} $, and any compact subset $ K \subset \widetilde{U}_{a} $, the quantity
\begin{align}\label{eq:brgj}
  {\rm sup}_{\lambda\in \mathbf{C}, |\lambda |= a} {\rm sup}_{s\in K} \left\|\left(\frac{\partial }{\partial s} \right)^{\alpha } \left(\lambda -L_{s}^{\pi_{*} \mathcal{D}^{V}  } \right)^{-1} u_{s} \right\|_{H^{\ell}} 
\end{align}
is finite and controlled by some constant multiple of some semi-norms of $ u $.

By \eqref{eq:cko1}, for all $ \ell\in \mathbf{N}$, the $ \lambda $-family of operators $ \left(M_{0} - \lambda  \right)^{-1} $, $ \lambda\in \mathbf{C} $ with $ \left| \lambda  \right| = a $, is $ \lambda $-uniformly bounded from $ H^{\ell} \left(X_{0 }, D^{V}_{|X_{0 } } \right) $ to $ H^{\ell+2} \left(X_{0 }, D^{V}_{|X_{0 } } \right) $.
By \eqref{eq:j2my}, the map $ s\in \widetilde{U}_{a} \to K_{s} \left(\lambda \right)$ with values in bounded operators from $ H^{\ell} \left(X_{0 }, D^{V}_{|X_{0 } } \right) $ to $ H^{\ell+1} \left(X_{0 }, D^{V}_{|X_{0 } } \right) $ is smooth, and the norms of itself and of its higher $ s $-derivatives are $ \lambda $-uniformly and $ s $-locally uniformly bounded.


We claim that $ \left(1+K_{s} \left(\lambda\right)\right)^{-1} $ is $ \lambda $-uniformly and $ s $-locally uniformly bounded from $ H^{\ell} \left(X_{0 }, D^{V}_{|X_{0 } } \right) $ to $ H^{\ell} \left(X_{0 }, D^{V}_{|X_{0 } } \right) $. 
Indeed, the case $ \ell= 0$ follows from construction.
For $ \ell_{1},\ell_{2} \in \mathbf{N} $, denote $ \left\| \cdot \right\|_{H^{\ell_{1} } ,H^{\ell_{2} } } $ the operator norm for operators from $ H^{\ell_{1} } \left(X_{0 },  D^{V} _{|X_{0 }  } \right) $ to $ H^{\ell_{2} }\left(X_{0 },  D^{V} _{|X_{0 }  } \right)  $.
Using 
\begin{align}\label{eq:zzwz}
 \left(1+K_{s}\left(\lambda \right)  \right)^{-1} = 1 - K_{s}\left(\lambda \right)\left(1+K_{s}\left(\lambda \right)  \right)^{-1},
\end{align}
if $ \ell\in \mathbf{N}^{*} $, we have 
\begin{align}\label{eq:ayo2}
\left\|  \left(1+K_{s}\left(\lambda \right)  \right)^{-1} \right\|_{H^{\ell}, H^{\ell} } \le  1 + \left\| K_{s}\left(\lambda \right) \right\|_{H^{\ell-1} ,H^{\ell} } \left\| \left(1+K_{s}\left(\lambda \right)  \right)^{-1} \right\|_{H^{\ell-1}, H^{\ell-1} }.
\end{align}
By \eqref{eq:ayo2} and by induction argument on $ \ell\in \mathbf{N} $, we get our claim.

We claim further that for any $ \ell\in \mathbf{N} $, the map $ s\in \widetilde{U}_{a} \to \left(1+K_{s}\left(\lambda \right)  \right)^{-1}$ with values in the bounded operators on $ H^{\ell} \left(X_{0 }, D^{V} _{|X_{0 } } \right) $ is smooth, and the norms of itself as well as all the higher $ s$-derivations are $ \lambda $-uniformly and $ s $-locally uniformly bounded.
Indeed, using the previous claim and resolvent identity, we see that the above map is continuous and has the required boundedness.
By induction argument, and using the corresponding properties on $ K_{s} \left(\lambda \right)$, we get our claim in full generality.


By \eqref{eq:cjys}, for all $ \ell\in \mathbf{N} $,  the map $ s\in \widetilde{U}_{a} \to u_{s}\in H^{\ell} \left(X_{0 },D^{V} _{|X_{0 } } \right)$ is smooth.
By \eqref{eq:3rny} and the above claims, the map $ s\in \widetilde{U}_{a} \to \left(\lambda -L_{s}^{\pi_{*} \mathcal{D}^{V}  } \right)^{-1} u_{s}\in H^{\ell} \left(X_{0 },D^{V} _{|X_{0 } } \right)$ is smooth and satisfies the required boundedness stated in \eqref{eq:brgj}.

Let $ \left(s_{0 }, x_{0 } \right)\in X$.
It remains to show the map $\left(s,x\right)\to \left(\lambda -L_{s}^{\pi_{*} \mathcal{D}^{V}  } \right)^{-1} u_{s}\left(x\right)$ is jointly smooth near $ \left(s_{0 }, x_{0 } \right)$.
Let $ \phi_{x_{0 } }\in C_{c}^{\infty} \left(X_{0 } \right)$ be such that $ \phi_{x_{0} } =1$ near $ x_{0} $.
If $ x $ is a point near $ x_{0 } $, if $ \delta_{x } $ is the Dirac distribution at $ x $, we have 
\begin{align}\label{eq:1yoa}
 \left(\lambda -L_{s}^{\pi_{*} \mathcal{D}^{V}  } \right)^{-1} u_{s}\left(x\right) = \left\langle \phi_{x_{0} }  \left(\lambda -L_{s}^{\pi_{*} \mathcal{D}^{V}  } \right)^{-1} u_{s}, \delta_{x} \right\rangle. 
\end{align}
Similar formulas for the higher $ x$-derivatives hold when we replace $ \delta_{x }$ by its corresponding higher derivatives.
Since the map with values in the classical Sobolev spaces $ s\in \widetilde{U}_{a} \to \phi_{x_{0} } \left(\lambda -L_{s}^{\pi_{*} \mathcal{D}^{V}  } \right)^{-1} u_{s}$  is smooth, we see that all the mixed higher partial derivatives of $ \left(\lambda -L_{s}^{\pi_{*} \mathcal{D}^{V}  } \right)^{-1} u_{s}\left(x\right) $ exist and are jointly continuous near $ \left(s_{0 }, x_{0 } \right)$.
Thus, $ \left(\lambda -L_{s}^{\pi_{*} \mathcal{D}^{V}  } \right)^{-1} u_{s}\left(x\right)$ is jointly smooth near $ \left(s_{0 }, x_{0 } \right)$.
The proof of our proposition is complete.
\end{proof}




\subsection{Spectral truncations}\label{s:trunc1}
Let $ a\in \left(0 ,a_{0 } \right) $ be such that $ \widetilde{S}_{a} \neq \varnothing $.

\begin{defin}\label{def:rzbe}
 If $s\in  \widetilde{S}_{a} $, put 
\begin{align}\label{eq:kkep}
  &P_{a,-,s}^{\pi_{*} \mathcal{D}^{V}  }  = \frac{1}{2i \pi } \int_{|\lambda |= a} \left(\lambda - L_{s} ^{\pi_{*} \mathcal{D}^{V}  } \right)^{-1} d \lambda,& P_{a,-,s}^{\pi_{*} \mathcal{E} }  = \frac{1}{2i \pi } \int_{|\lambda |= a} \left(\lambda - L_{s} ^{\pi_{*} \mathcal{E} }  \right)^{-1} d \lambda.
\end{align}
Set 
\begin{align}\label{eq:xryd}
   &P_{ a,+,s}^{\pi_{*} \mathcal{D}^{V}  }= 1- P_{a,-,s}^{\pi_{*} \mathcal{D}^{V}  },&P_{ a,+,s}^{\pi_{*} \mathcal{E} } = 1-P_{a,-,s}^{\pi_{*} \mathcal{E} }.
\end{align}
\end{defin}

\begin{defin}\label{def:rfyl}
 For $ s\in \widetilde{S} _{a} $, denote  
  \begin{align}\label{eq:vuby}
&H^{\infty} _{a,\pm,s}   \left(X_{0 } ,D^{V} _{|X_{0 } }  \right) =  P_{a,\pm,s}^{\pi_{*} \mathcal{D}^{V} } H^{\infty}  \left(X_{0 }, D^{V} _{|X_{0 } } \right),   & H^{\infty}_{a,\pm,s}  \left(X_{0 }, E_{|X_{0 } } \right) = P_{a,\pm,s}^{\pi_{*} \mathcal{E} } H^{\infty}  \left(X_{0 }, E_{|X_{0 } } \right).
  \end{align}
\end{defin}

By Proposition \ref{prop:sthv}, $ H^{\infty} _{a,-,s}   \left(X_{0 } ,D^{V} _{|X_{0 } }  \right)$,  $ H^{\infty}_{a,-,s}  \left(X_{0 }, E_{|X_{0 } } \right) $ are the finite-dimensional spaces of the direct sums of characteristic spaces of $ L^{\pi_{*} \mathcal{D}^{V}  }_{s},  L^{\pi_{*} \mathcal{E} }_{s}$ associated to the eigenvalues $\lambda <a$. 
We have the decompositions, 
\begin{align}\begin{aligned}\label{eq:wcnc}
H^{\infty}   \left(X_{0 } ,D^{V} _{|X_{0 } }  \right) &=  H^{\infty} _{a,-,s}   \left(X_{0 } ,D^{V} _{|X_{0 } }  \right) \oplus H^{\infty} _{a,+,s}   \left(X_{0 } ,D^{V} _{|X_{0 } }  \right), \\
 H^{\infty}  \left(X_{0 }, E_{|X_{0 } } \right) &= H^{\infty}_{a,-,s}  \left(X_{0 }, E_{|X_{0 } } \right) \oplus H^{\infty}_{a,+,s}  \left(X_{0 }, E_{|X_{0 } } \right).
\end{aligned}\end{align}



\begin{prop}\label{prop:rgpe}
  If $ s_{0 }\in \widetilde{S}_{a} $, there exists an open neighbourhood $ \widetilde{U}_{s_{0 }} \subset \widetilde{S}_{a} $ of $ s_{0 }$ such that for all $ s\in \widetilde{U}_{s_{0 }} $, we have isomorphisms of vector spaces,
  \begin{align}\begin{aligned}\label{eq:iapi}
 P_{a,-,s}^{\pi_{*} \mathcal{D}^{V}  }  &:  H^{\infty}_{a,-,s_{0 } }  \left(X_{0 }, D^{V} _{|X_{0 } }  \right) \to H^{\infty}_{a,-,s}  \left(X_{0 }, D^{V} _{|X_{0 }  } \right),\\
 P_{a,-,s}^{\pi_{*} \mathcal{E} }  &:  H^{\infty}_{a,-,s_{0 } }  \left(X_{0 }, E_{|X_{0 } }  \right) \to H^{\infty}_{a,-,s}  \left(X_{0 }, E_{|X_{0 }  } \right).\\
 \end{aligned}\end{align}
\end{prop}
\begin{proof}
 Let us show the first linear map in \eqref{eq:iapi} is an isomorphism.
 The proof for the second one is the same.

 By Proposition \ref{prop:cflb} and Definition \ref{def:rzbe},  $s\in \widetilde{S}_{a} \to  P_{a,-,s_{0 } }^{\pi_{*} \mathcal{D}^{V}  } P_{a,-,s}^{\pi_{*} \mathcal{D}^{V}  } P_{a,-,s_{0 } }^{\pi_{*} \mathcal{D}^{V}  }$ is a smooth map with values in $\End\left( H^{\infty}_{a,-,s_{0 } }  \left(X_{0 }, D^{V} _{|X_{0 } }  \right)\right)$.
 It is equal to the identity when $ s= s_{0 } $.
 There exists a small open neighbourhood $ \widetilde{U}_{s_{0 } } $ of $ s_{0} $ such that if $ s\in \widetilde{U}_{s_{0 } } $, the linear map $ P_{a,-,s_{0 } }^{\pi_{*} \mathcal{D}^{V}  } P_{a,-,s}^{\pi_{*} \mathcal{D}^{V}  } P_{a,-,s_{0 } }^{\pi_{*} \mathcal{D}^{V}  }$  is bijective, so that the first linear map in \eqref{eq:iapi} is injective. 
 
 We claim that 
 the map 
 \begin{align}\label{eq:gefj}
  s\in \widetilde{S}_{a} \to \dim H^{\infty}_{a,-,s} \left(X_{0 }, D^{V}_{|X_{0 } } \right)\in \mathbf{N} 
 \end{align}
 is locally constant.
 Indeed, since $ \left(\lambda -M_{0 } \right)^{-1} $ is holomorphic on $ \lambda\in \mathbf{C} $ with $ \left| \lambda \right| \le  a $, by \eqref{eq:h2u2} and \eqref{eq:kkep}, we have 
  \begin{align}\begin{aligned}\label{eq:fkey}
   P_{a,-,s}^{\pi_{*} \mathcal{D}^{V}  } & = \frac{1}{2i \pi } \int_{|\lambda |= a} \left(\lambda - L_{s} ^{\pi_{*} \mathcal{D}^{V}  } \right)^{-1} -\left(\lambda - M_{0 }  \right)^{-1} d \lambda\\
&   = \frac{1}{2i \pi } \int_{|\lambda |= a} \left(\lambda - L_{s} ^{\pi_{*} \mathcal{D}^{V}  } \right)^{-1} R_{s} \left(\lambda - M_{0 }  \right)^{-1} d \lambda.
  \end{aligned}\end{align}
  Repeating the above procedure, for $ N\in \mathbf{N} $, we have
  \begin{align}\label{eq:ymjp}
   P_{a,-,s}^{\pi_{*} \mathcal{D}^{V}  }  = \frac{1}{2i \pi } \int_{|\lambda |= a} \left(\lambda - L_{s} ^{\pi_{*} \mathcal{D}^{V}  } \right)^{-1} \left(R_{s} \left(\lambda - M_{0 }  \right)^{-1}\right)^{N}  d \lambda.
  \end{align}
  If $ N= 2n+2$, by Proposition \ref{prop:ul1a} \ref{prop:ul1abd}), $\left(\lambda - L_{s} ^{\pi_{*} \mathcal{D}^{V}  } \right)^{-1} \left(R_{s} \left(\lambda - M_{0 }  \right)^{-1}\right)^{N}$ is of trace class.
 Moreover, the proof given there shows that the family $\left(\lambda - L_{s} ^{\pi_{*} \mathcal{D}^{V}  } \right)^{-1} \left(R_{s} \left(\lambda - M_{0 }  \right)^{-1}\right)^{N}$, $ \left(s,\lambda \right)\in \widetilde{S}_{a} \times \left\{\lambda\in \mathbf{C} : \left| \lambda \right| = a \right\}$, of trace class operators, is jointly continuous.
 Therefore, the trace of $ P_{a,-,s}^{\pi_{*} \mathcal{D}^{V}  }$, which coincides with $ \dim H^{\infty}_{a,-,s} \left(X_{0 }, D^{V}_{|X_{0 } } \right) $, is continuous on $ \widetilde{S}_{a} $, and therefore is locally constant.

 By our claim, the injectivity of the first map in \eqref{eq:iapi} on $ \widetilde{U}_{s_{0 } } $ implies its surjectivity on the same open neighbourhood.
 The proof of our proposition is completed.
\end{proof}

\begin{prop}\label{prop:qvah}
 The families $\left\{H^{\infty} _{a,-,s}   \left(X_{0 } ,D_{|X_{0 }  }^{V} \right) \right\} _{s\in \widetilde{S} _{a} }, \left\{H^{\infty}_{a,-,s}  \left(X_{0 }, E_{0 } \right) \right\}_{s\in \widetilde{S} _{a} }$ of vector spaces together with the trivialisations given by \eqref{eq:iapi} define smooth vector bundles $ H^{\infty} _{a,-} \left(X_{0 }, D^{V} _{|X_{0 } } \right), H^{\infty}_{a,-}  \left(X_{0 }, E_{|X_{0 }  } \right) $ on $ \widetilde{S} _{a} $.
\end{prop}
\begin{proof}
It is enough to consider the case $ D^{V} $. 
Let $ s_{0 }, s_{1}\in \widetilde{S}_{a} $.
Let $ \widetilde{U} _{s_{0 } } $ and $ \widetilde{U}_{s_{1} }  $ be as in Proposition \ref{prop:rgpe}.
Assume that $ \widetilde{U} _{s_{0 } } \cap \widetilde{U} _{s_{1} } \neq \varnothing  $.
Let us show the transition map 
\begin{align}\label{eq:tkpr}
 \left(P^{\pi_{*} \mathcal{D}^{V}  }_{a,-,s }P^{\pi_{*} \mathcal{D}^{V}  }_{a,-,s_{1} }\right)^{-1}   \left(P^{\pi_{*} \mathcal{D}^{V}  }_{a,-,s }P^{\pi_{*} \mathcal{D}^{V}  }_{a,-,s_{0 } }\right)  : H^{\infty}_{a,-,s_{0 } }  \left(X_{0 }, D^{V}_{|X_0 } \right) \to H^{\infty}_{a,-,s_{1} }  \left(X_{0 }, D^{V}_{|X_{0} } \right)
\end{align}
is smooth on  $ s\in \widetilde{U} _{s_{0 } } \cap \widetilde{U} _{s_{1} } $.

Let $ \left\{u_{i}\right\},\left\{v_{j}\right\}$ be bases of $ H^{\infty}_{a,-,s_{0 } }  \left(X_{0 }, D^{V} _{|X_{0}  } \right)$, $ H^{\infty}_{a,-,s_{1 } }  \left(X_{0 }, D^{V} _{|X_{0}  } \right)$.
By Proposition \ref{prop:rgpe}, for each $ i$, there exist a family of functions $ \left\{f_{ij} \right\}_{j} $ on $ \widetilde{U} _{s_{0 } } \cap \widetilde{U} _{s_{1} }$, such that 
\begin{align}\label{eq:mocq}
 P^{\pi_{*} \mathcal{D}^{V}  }_{a,-,s } u_{i} = \sum_{j} f_{ij,s} P^{\pi_{*} \mathcal{D}^{V}  }_{a,-,s } v_{j}.
\end{align}
The transition map \eqref{eq:tkpr} under the bases $ \left\{u_{i}\right\},\left\{v_{j}\right\}$ is now given by the matrix $ \left(f_{ij}\right)$.

By \eqref{eq:mocq}, given $ i,k$, we have 
\begin{align}\label{eq:jcgd}
 \left\langle P^{\pi_{*} \mathcal{D}^{V}  }_{a,-,s } u_{i} , P^{\pi_{*} \mathcal{D}^{V}  }_{a,-,s } v_{k} \right\rangle = \sum_{j} f_{ij,s} \left\langle P^{\pi_{*} \mathcal{D}^{V}  }_{a,-,s } v_{j},P^{\pi_{*} \mathcal{D}^{V}  }_{a,-,s } v_{k} \right\rangle.
\end{align}
By Proposition \ref{prop:cflb} and Definition \ref{def:rzbe}, on $ \widetilde{U}_{s_{0 } } \cap \widetilde{U}_{s_{1} } $, $ \left\langle P^{\pi_{*} \mathcal{D}^{V}  }_{a,-,s } u_{i} , P^{\pi_{*} \mathcal{D}^{V}  }_{a,-,s } v_{k} \right\rangle$ and $ \left\langle P^{\pi_{*} \mathcal{D}^{V}  }_{a,-,s } v_{j},P^{\pi_{*} \mathcal{D}^{V}  }_{a,-,s } v_{k} \right\rangle$ are smooth.
By constructions, the matrix  $ \left(\left\langle P^{\pi_{*} \mathcal{D}^{V}  }_{a,-,s } v_{j},P^{\pi_{*} \mathcal{D}^{V}  }_{a,-,s } v_{k} \right\rangle\right)$ is invertible.
Therefore, all the $ f_{ij} $ are smooth on $ \widetilde{U}_{s_{0 } } \cap \widetilde{U}_{s_{1} } $, so that the transition map \eqref{eq:tkpr} is smooth.

Note that the cocycle condition on the transition map is trivial.
The proof of our proposition is complete.
\end{proof}

Take an open set $ S_{a} \subset S \cap \widetilde{S}_{a} $ such that its closure $ \overline{S} _{a} $ in $ \widetilde{S} $ is contained in $ \widetilde{S}_{a} $.
Then, $ S_{a}, \overline{S}_{a}, \widetilde{S}_{a} $ have properties analogous to those of $ S,\overline{S}, \widetilde{S}$.
Also, 
\begin{align}\label{eq:qrpy}
 \overline{S}_{a} \subset \overline{S} \cap \widetilde{S}_{a}.
\end{align}

By Proposition \ref{prop:cflb} and Definition \ref{def:rzbe}, the operators $ P_{a,\pm,s}^{\pi_{*} \mathcal{D}^{V}  }$ and $ P_{a,\pm,s}^{\pi_{*} \mathcal{E} }$ induce the corresponding operators $ P_{a,\pm}^{\pi_{*} \mathcal{D}^{V}  }$ and $ P_{a,\pm}^{\pi_{*} \mathcal{E} }$ acting respectively on $ \left(\pi_{*} \mathcal{D}^{V}  \right)_{\left(2\right)|\overline{S}_{a}}$ and $ \left(\pi_{*} \mathcal{E} \right)_{\left(2\right)|\overline{S}_{a}}$.
By Propositions \ref{prop:yemy} \ref{prop:yemya}) and \ref{prop:ywts} \ref{prop:ywtsa}), $ P_{a,\pm}^{\pi_{*} \mathcal{D}^{V}  }$ are morphisms of $ \mathscr{O}^{\infty} $-complexes on $ \overline{S}_{a}  $, and $ P_{a,\pm}^{\pi_{*} \mathcal{E} }$ are morphisms in $ {\rm {\widetilde{B}}^{t}} \left(\overline{S}_{a} \right)$.

\begin{defin}\label{def:ewmm}
 Let $ \left(\pi_{*} \mathcal{D}^{V}  \right)_{\left(2\right),a,\pm}$ and $ \left(\pi_{*} \mathcal{E} \right)_{\left(2\right),a,\pm}$ be the subsheaves of $ \left(\pi_{*} \mathcal{D}^{V}  \right)_{\left(2\right)|{\overline{S}}_{a}}$ and $ \left(\pi_{*} \mathcal{E} \right)_{\left(2\right)|{\overline{S} }_{a}}$ defined on open subsets $ \widetilde{U} _{a} \subset {\widetilde{S}} _{a}  $ by 
 \begin{align}\begin{aligned}\label{eq:hntb}
  \left(\pi_{*} \mathcal{D}^{V}  \right)_{\left(2\right),a,\pm }\left( \widetilde{U}_{a} \right)&= P_{ a,\pm}^{\pi_{*} \mathcal{D}^{V}  } \left\{\left(\pi_{*} \mathcal{D}^{V}  \right)_{\left(2\right)}\left(\widetilde{U}_{a} \right)\right\},\\
  \left(\pi_{*} \mathcal{E} \right)_{\left(2\right),a,\pm }\left( \widetilde{U}_{a} \right)&= P_{ a,\pm}^{\pi_{*} \mathcal{E} } \left\{\left(\pi_{*} \mathcal{E} \right)_{\left(2\right)}\left(\widetilde{U} _{a} \right)\right\}.
 \end{aligned}\end{align}
\end{defin}

Then, $ v_{0 }^{\pi_{*} \mathcal{D}^{V} } $ and $ A^{\pi_{*} \mathcal{E} \prime \prime } $ restrict to operators $ v^{\pi_{*} \mathcal{D}^{V} }_{0, a,\pm} $ and $ A^{\pi_{*} \mathcal{E} \prime \prime }_{a,\pm} $ acting respectively on $ \left(\pi_{*} \mathcal{D}^{V} \right)_{\left(2\right),a,\pm} $ and $ \left(\pi_{*} \mathcal{E} \right)_{\left(2\right),a,\pm} $, so that $ \left(\pi_{*} \mathcal{D}^{V} \right)_{\left(2\right),a,\pm} $ is an $ \mathscr{O}^{\infty} $-complex on $ \overline{S}_{a} $, and $ \left(\pi_{*} \mathcal{E} \right)_{\left(2\right),a,\pm} $ is an object in $ {\rm \widetilde{B}^{t}}\left(\overline{S}_{a} \right)$.
We have decompositions of complexes,
\begin{align}\begin{aligned}\label{eq:zps3}
 \left(\pi_{*} \mathcal{D}^{V}  \right)_{\left(2\right)|{\overline{S} }_{a} }&=  \left(\pi_{*} \mathcal{D}^{V}  \right)_{\left(2\right),a,- }\oplus \left(\pi_{*} \mathcal{D}^{V}  \right)_{\left(2\right),a,+ },\\
 \left(\pi_{*} \mathcal{E} \right)_{\left(2\right)|\overline{S}_{a} }&=  \left(\pi_{*} \mathcal{E}  \right)_{\left(2\right),a,- }\oplus \left(\pi_{*} \mathcal{E}  \right)_{\left(2\right),a,+ }.
\end{aligned}\end{align}

\begin{prop}\label{prop:lvqq}
  The morphism of $ \mathscr{D} $-modules on $ \overline{S}_{a}  $, 
  \begin{align}\label{eq:jhdj}
   P_{a,\pm}^{\pi_{*} \mathcal{E} } : \mathscr{D} \widehat{\otimes}_{\mathscr{O}^{\infty} }  \left(\pi_{*} \mathcal{D}^{V} \right)_{\left(2\right),a,\pm} \to \left(\pi_{*} \mathcal{E} \right)_{\left(2\right),a,\pm},
  \end{align}
  is an isomorphism, so that $ \left(\pi_{*}\mathcal{D}^{V} \right)_{\left(2\right),a,\pm} $ is the diagonal of $ \left(\pi_{*}\mathcal{E} \right)_{\left(2\right),a,\pm} $.
\end{prop}
\begin{proof}
 The proof is the same as the one given in \cite[(8.10.19), (11.1.12)]{BismutShenWei23}.
\end{proof}

\begin{thm}\label{thm:mitj}
  The sheaves $ \left(\pi_{*} \mathcal{D}^{V}  \right)_{\left(2\right),a,-}$ and $ \left(\pi_{*} \mathcal{E} \right)_{\left(2\right),a,-}$ are locally free of finite rank.
  More precisely, if $ s_{0 }\in \widetilde{S}_{a} $, there exists a small open neighbourhood $ \widetilde{U}_{s_{0 } } \subset \widetilde{S}_{a}  $ of $ s_{0 } $ such that,  
   \begin{align}\begin{aligned}\label{eq:g1zi}
     \left(\pi_{*} \mathcal{D}^{V}  \right)_{\left(2\right),a,-} \left(\widetilde{U}_{s_{0 } }   \right)&= C^{\infty} \left(\widetilde{U} _{s_{0 } }  ,  H^{\infty} _{a,-}   \left(X_{0 } ,D^{V} _{|X_{0 }  }  \right)_{|\widetilde{U} _{s_{0 } }}  \right),\\
    \left(\pi_{*} \mathcal{E} \right)_{\left(2\right),a,-} \left(\widetilde{U} _{s_{0 } }   \right)&= C^{\infty} \left(\widetilde{U} _{s_{0 } }  ,  H^{\infty}_{a,-}   \left(X_{0 } ,E_{|X_{0 }  }  \right)_{|\widetilde{U} _{s_{0 } }}  \right).
   \end{aligned}\end{align}
   In particular, $ \left(\pi_{*} \mathcal{E} \right)_{\left(2\right),a,-}$ is defined by an object in $ {\rm B^{b,t}}\left( \overline{S}_{a} \right) $.
\end{thm}
\begin{proof} As before, we only need to show our theorem for $ \left(\pi_{*} \mathcal{D}^{V}  \right)_{\left(2\right),a,-}$, or sufficiently if $ \widetilde{U}_{s_{0 } } $ is as in Proposition \ref{prop:rgpe}, 
 \begin{align}\label{eq:phxn}
  P_{a,-}^{\pi_{*} \mathcal{D}^{V}  } \left\{\left(\pi_{*} \mathcal{D}^{V}  \right)_{\left(2\right)} \left(\widetilde{U}_{s_{0 } } \right) \right\}= C^{\infty} \left(\widetilde{U} _{s_{0 } } ,H^{\infty}_{a,-}  \left(X_{0 } ,D^{V} _{|X_{0 } } \right)_{|\widetilde{U} _{s_{0 } }} \right). 
 \end{align}

With the notation in the proof of Proposition \ref{prop:qvah}, we see that any element of $ C^{\infty} \left(\widetilde{U} _{s_{0 } } ,H^{\infty}_{a,-}  \left(X_{0 } ,D^{V} _{|X_{0 } } \right)_{|\widetilde{U} _{s_{0 } }} \right)$ has the form 
\begin{align}\label{eq:2sha}
s\in \widetilde{U}_{s_{0} } \to \sum_{i} f_{i,s} P_{a,-,s}^{\pi_{*} \mathcal{D}^{V} }  u_{i},
\end{align}
where $ s\in \widetilde{U}_{s_{0 } } \to f_{i,s} $ is a smooth function. 
This gives us the inverse inclusion $ \supset $ in \eqref{eq:phxn}.

If $ u\in \left(\pi_{*} \mathcal{D}^{V}  \right)_{\left(2\right)} \left(\widetilde{U} _{s_{0 } } \right)  $, there exists a family of functions $\left\{ f_{i}\right\}_{i}  $ on $ \widetilde{U}_{s_{0 } } $ such that 
\begin{align}\label{eq:qmig}
  P_{a,-}^{\pi_{*} \mathcal{D}^{V}   }u= \sum_{i}^{} f_{i} P_{a,-}^{\pi_{*} \mathcal{D}^{V}  }u_{i}. 
\end{align}
By \eqref{eq:qmig}, we have 
\begin{align}\label{eq:3ums}
 \left\langle  P_{a,-,s}^{\pi_{*} \mathcal{D}^{V}  }u_{|X_s} , P_{a,-,s}^{\pi_{*} \mathcal{D}^{V}  }u_{i} \right\rangle = \sum_{i}^{} f_{i,s} \left\langle P_{a,-,s}^{\pi_{*} \mathcal{D}^{V}  }u_{i}, P_{a,-,s}^{\pi_{*} \mathcal{D}^{V}  }u_{j} \right\rangle. 
\end{align}
Since $ \left\langle  P_{a,-,s}^{\pi_{*} \mathcal{D}^{V}  }u_{|X_s} , P_{a,-,s}^{\pi_{*} \mathcal{D}^{V}  }u_{i} \right\rangle$, $ \left\langle P_{a,-,s}^{\pi_{*} \mathcal{D}^{V}  }u_{i}, P_{a,-,s}^{\pi_{*} \mathcal{D}^{V}  }u_{j} \right\rangle$ are smooth, and since the matrix
$ \left(\left\langle P_{a,-,s}^{\pi_{*} \mathcal{D}^{V}  }u_{i}, P_{a,-,s}^{\pi_{*} \mathcal{D}^{V}  }u_{j} \right\rangle\right)$ is invertible, we see that all the $ f_{i} $ are smooth on $ \widetilde{U}_{s_{0 } } $.
This gives us the inclusion $ \subset $ in \eqref{eq:phxn}.

Therefore, we get \eqref{eq:phxn} and finish the proof of our theorem.
\end{proof}

\begin{prop}\label{prop:sn2o}
 The complexes $ \left(\pi_{*} \mathcal{D}^{V} \right)_{\left(2\right),a,+} $ and $ \left(\pi_{*} \mathcal{E} \right)_{\left(2\right),a,+} $ are acyclic.
\end{prop}
\begin{proof}
As before, we only need consider the complex $ \left(\pi_{*} \mathcal{D}^{V} \right)_{\left(2\right),a,+} $.

If $ s\in \widetilde{S}_{a} $, the meromorphic family of operators $ \lambda^{-1} P^{\pi_{*} \mathcal{D}^{V} }_{a,+,s}\left(L^{\pi_{*} \mathcal{D} ^{V} }_{s}  -\lambda \right)^{-1}P^{\pi_{*} \mathcal{D}^{V} }_{a,+,s} $ on $ \lambda\in \mathbf{C} $ with $ \left| \lambda  \right| \le a$ has only one pole at $ \lambda =0$. 
It is simple and its residue is given by $ P^{\pi_{*} \mathcal{D}^{V} }_{a,+,s}\left(L^{\pi_{*} D^{V} }_{s}  \right)^{-1}P^{\pi_{*} \mathcal{D}^{V} }_{a,+,s}$.
Therefore, 
 \begin{align}\label{eq:pnyt}
 \frac{1}{2i \pi }\int_{\left| \lambda  \right| = a}^{}\lambda^{-1} P^{\pi_{*} \mathcal{D}^{V} }_{a,+,s}\left(L^{\pi_{*} \mathcal{D} ^{V} }_{s}  -\lambda \right)^{-1}P^{\pi_{*} \mathcal{D}^{V} }_{a,+,s}  d \lambda = P^{\pi_{*} \mathcal{D}^{V} }_{a,+,s}\left(L^{\pi_{*} \mathcal{D} ^{V} }_{s}  \right)^{-1}P^{\pi_{*} \mathcal{D}^{V} }_{a,+,s}.  
\end{align}
By Proposition \ref{prop:cflb}, the left-hand side of \eqref{eq:pnyt} induces an operator acting on $ \left(\pi_{*} \mathcal{D}^{V} \right)_{\left(2\right),a,+} $ and will be denoted by $ P^{\pi_{*} \mathcal{D}^{V} }_{a,+}\left(L^{\pi_{*} \mathcal{D} ^{V} }  \right)^{-1}P^{\pi_{*} \mathcal{D}^{V} }_{a,+}$.
 
Put 
\begin{align}\label{eq:z31w}
 v^{\pi_{*} \mathcal{D}^{V} }_{0, *,a,+} = P^{\pi_{*} \mathcal{D}^{V} }_{a,+}v_{0, *}^{\pi_{*} \mathcal{D}^{V} }  P^{\pi_{*} \mathcal{D}^{V} }_{a,+}.
\end{align}
By \eqref{eq:wjav}, since $ v^{\pi_{*} \mathcal{D}^{V} }_{0} $ commutes with $ P^{\pi_{*} \mathcal{D}^{V} }_{a,+}$, we have 
\begin{align}\label{eq:ktuc}
 \left[v^{\pi_{*} \mathcal{D}^{V} }_{0, a,+} , v_{0,*, a,+} ^{\pi_{*} \mathcal{D}^{V} }\right]= P^{\pi_{*} \mathcal{D}^{V} }_{a,+} L^{\pi_{*} \mathcal{D} ^{V} } P^{\pi_{*} \mathcal{D}^{V} }_{a,+}.
\end{align}
Since $ v^{\pi_{*} \mathcal{D}^{V} }_{0, a,+}$ commutes with $ P^{\pi_{*} \mathcal{D}^{V} }_{a,+} \left(L^{\pi_{*} \mathcal{D} ^{V} }\right)^{-1} P^{\pi_{*} \mathcal{D}^{V} }_{a,+}$, by \eqref{eq:ktuc}, we have 
\begin{align}\label{eq:y33v}
  \left[v^{\pi_{*} \mathcal{D}^{V} }_{0, a,+} , v_{0,*, a,+} ^{\pi_{*} \mathcal{D}^{V} } P^{\pi_{*} \mathcal{D}^{V} }_{a,+} \left(L^{\pi_{*} D^{V} }\right)^{-1} P^{\pi_{*} \mathcal{D}^{V} }_{a,+}\right]= 1.
\end{align}
This implies $ \left(\pi_{*} \mathcal{D}^{V} \right)_{\left(2\right),a,+} $ is acyclic, and finishes the proof of our proposition.
\end{proof}

\subsection{Proof of Theorem \ref{thm:et3l}}\label{s:p2}
The statement of Theorem \ref{thm:et3l} is local.
We can assume that $ \widetilde{S}$ is a small open ball in $ \mathbf{C}^{m} $ and there is some $ a\in \left(0,a_{0 } \right)$ such that $\widetilde{S} =  \widetilde{S}_{a} $ and choose $ S_{a} =S$.

By Propositions \ref{prop:iu2w}, \ref{prop:sn2o}, and by \eqref{eq:zps3}, the composition of the following embeddings is a quasi-isomorphism 
\begin{align}\label{eq:3mtm}
 \left(\pi_{*} \mathcal{E}\right)_{\left(2\right),a,-} \to  \left(\pi_{*} \mathcal{E}\right)_{\left(2\right)} \to  \pi_{*} \mathcal{E}.
\end{align}
By Theorem \ref{thm:mitj},  $ \left(\pi_{*} \mathcal{E}\right)_{\left(2\right),a,-}$ is defined by an element of $ {\rm B^{b,t}}\left(\overline{S} \right) $.
By \cite[Theorem 5.3.4]{BismutShenWei23},  $ \left(\pi_{*} \mathcal{E}\right)_{\left(2\right),a,-}$ has coherent cohomology, which gives Theorem \ref{thm:et3l}.\qed

\subsection{Proof of Proposition \ref{prop:wvkt}} \label{s:p10}
As before, the statement of Proposition \ref{prop:wvkt} is local.
We assume that $\widetilde{S}$ is a small open ball in $\mathbf{C}^{m}$, that $\widetilde{S} = \widetilde{S}_{a}$ for some $a \in (0, a_{0})$, and that $S_{a} = S$.

As in \eqref{eq:3mtm}, the vertical arrows defined by the inclusions in the following two diagrams \eqref{dia:aqkd0} and \eqref{dia:aqkd} are quasi-isomorphisms.
By this consideration and by Footnote \ref{footnoteDFA}, the conclusion of Proposition \ref{prop:2gpa} still holds for the diagram \eqref{dia:aqkd0} below.
Then, the morphism $ \phi :  \mathcal{E}_{\overline{S} }  \to \pi_{*} \mathcal{E}$ can be lifted to $ \phi_{a,-} $,  so that the diagram
\begin{align}\label{dia:aqkd0}
    \begin{tikzsd}
       &   \left(\pi_{*} \mathcal{E} \right)_{\left(2\right),a,-} \\
      \mathcal{E}_{\overline{S} }   & \pi_{*} \mathcal{E} .
      \arrow["\phi_{a,-}", from=2-1, to=1-2, dotted]
      \arrow["", from=1-2, to=2-2]
      \arrow["\phi", from=2-1, to=2-2]
    \end{tikzsd}
\end{align}
commutes in $ {\rm \widetilde{\underline{B} }^{t}} \left(\overline{S} \right) $.
Taking the diagonal, by Proposition \ref{prop:lvqq}, we get a diagram
\begin{align}\label{dia:aqkd}
    \begin{tikzsd}
       &   \left(\pi_{*} \mathcal{D}^{V}  \right)_{\left(2\right),a,-} \\
      \mathcal{D}_{\overline{S} }   & \pi_{*} \mathcal{D}^{V} .
      \arrow["\phi_{a,-,0 }", from=2-1, to=1-2]
      \arrow["", from=1-2, to=2-2]
      \arrow["\phi_{0 } ", from=2-1, to=2-2]
    \end{tikzsd}
\end{align}
which commutes in the homotopy category of $ \mathscr{O}^{\infty} $-complexes on $ \overline{S} $.

Then $ \phi $ is a quasi-isomorphism if and only if $ \phi_{a,-} $ is a quasi-isomorphism.
By \cite[Proposition 6.4.1]{BismutShenWei23}, this is equivalent to $ \phi_{a,-,0 } $ is a quasi-isomorphism, which is again equivalent to $ \phi_{0 } $ is a quasi-isomorphism.
The proof of Proposition \ref{prop:wvkt} is complete.\qed

\subsection{Proof of Theorem \ref{thm:ete1} \ref{thm:etea})} \label{s:p1}
Theorem \ref{thm:ete1} \ref{thm:etea}) can be established by the same method used in the proof of Theorem \ref{thm:poev} \ref{thm:poeva}), or \cite[Theorem 6.3.6]{BismutShenWei23} since tameness is an empty condition here.

Since the cohomology of $ \pi_{*} \mathcal{D} ^{V} $ is locally isomorphic to the cohomology of some $\left( \pi_{*} \mathcal{D} ^{V}\right)_{\left(2\right),a,-} $, by \cite[Theorem 5.2.1]{BismutShenWei23}, we know that $ \pi_{*} \mathcal{D} ^{V} $ has the same property as the complex $ \mathcal{G} $ introduced in Section \ref{sEssSuj}.

Theorem \ref{thm:ete1} \ref{thm:etea}) follows from Proposition \ref{prop:cfiv} together with a version of Proposition \ref{prop:axew} (see Footnote \ref{footnoteDFA}).\qed

\subsection{Projection in derived category}\label{s:proj}
Recall that $ \widetilde{S} $ is a complex manifold and $ \widetilde{M} = \widetilde{S} \times X$.  
 
\begin{thm}\label{thm:fbt3}
  If $ {\mathcal{F} } $ is an object of $ {\rm D^{b}_{coh}} \left( \widetilde{M}  \right)$ such that for any compact subset $ K $ of $ \widetilde{S}$,  
  \begin{align}\label{eq:rb3y}
    {\rm Supp}\left({{\mathcal{F}} } \right) \cap \pi^{-1} \left(K \right) 
  \end{align}
  is compact, then $ R \pi_{*} {{\mathcal{F}} } $ is an object of $ {\rm D^{b}_{coh}} \left(\widetilde{S} \right)$.
  Moreover, if $ c_{0 } ,c_{1} \in \mathbf{Z} $ with $ c_{0 }  \le c_{1} $, if $ \mathcal{H} {{\mathcal{F}}  } $ is concentrated in degrees $ \left[c_{0 } ,c_{1} \right]$, then the cohomology of $ R \pi_{*} {{\mathcal{F}} } $ is concentrated in degrees $ \left[c_{0} ,c_{1} +n\right]$.
  \end{thm}
 \begin{proof}
  Since the statement is local on $ \widetilde{S} $, it is enough to show our theorem for any open relatively compact subset $ S $ of $ \widetilde{S} $.

  By Example \ref{exa:nqs3}, $ {{\mathcal{F}} } $ induces an object with compact support in $ {\rm C^{b,t}_{coh}}\left(\overline{M} \right)  $, which will still be denoted by $ \mathcal{F} $. 
  By Theorem \ref{thm:poev} \ref{thm:poeva}), there exist an object $ \left(E,A^{E\prime \prime } \right) $ in  $  {\rm B^{b,t}_{c}}\left(\overline{M} \right)  $ with associated $ \mathscr{O} $-complex $ \mathcal{E} $ and a quasi-isomorphism in  ${\rm \widetilde{B} ^{t}} \left( \overline{M}\right)$, 
  \begin{align}\label{eq:okxg} 
   \mathcal{E} \to \mathscr{D} \mathcal{F}. 
  \end{align}  
  Since $ \mathcal{E} $ is soft, by \cite[Proposition 1.8.3]{KashiwaraSchpira} and \cite[IV (4.14) Proposition]{DemaillyBook}, we have 
\begin{align}\label{eq:mpnc}
 R\pi_{*} \mathcal{E} = \pi_{*} \mathcal{E}. 
\end{align}
By \eqref{eq:gvcd}, \eqref{eq:okxg}, and \eqref{eq:mpnc}, the cohomology of $ \pi_{*} \mathcal{E} $ is isomorphic to the cohomology of $ R \pi_{*} \mathcal{F} $ on $ \overline{S} $.  
By Theorem \ref{thm:et3l}, $  R \pi_{*} \mathcal{F}$ has coherent cohomology on $ \overline{S} $.
  
  By the discussion following \eqref{eq:simg}, the  cohomology of $ R\pi_{*} {\mathcal{F} }  $ is concentrated in degrees $ \ge c_{0 } $.
  Let us show it is also concentrated in degrees $ \le c_{1} +n$.
  By Remark \ref{re:dqvn}, the diagonal $ D $ of $ \left(E,A^{E\prime \prime } \right)$ can be constructed to be concentrated in degrees $ \le c_{1} $.
  Then, $ D^{V} $ is concentrated in degrees $ \le c_{1} +n$.
  If $ a\in (0,a_{0 } )$, $ \left(\pi_{*} \mathcal{D}^{V} \right)_{\left(2\right),a,-}$ is also concentrated in degrees $ \le c_{1} +n$.
  By \cite[(5.3.5)]{BismutShenWei23}, the cohomology of $ \left(\pi_{*} \mathcal{E}\right)_{\left(2\right),a,-} $ is again concentrated in degrees $ \le c_{1} +n$.
  This finishes the proof of our theorem.
\end{proof} 

\subsection{Proof of Theorem \ref{thm:main} \textnormal{\ref{main2a}}) when $ X$ is smooth}\label{S:Gs}
Under the assumption that $ X$ is a smooth complex manifold, let us establish Theorem \ref{thm:main} \textnormal{\ref{main2a}}) for any  proper holomorphic map $ f: X\to Y$, and moreover, 
\begin{align}\label{eq:hewx}
 c\left(f\right)\le \dim X.
\end{align}
By Remark \ref{re:fq3s}, we can also assume that $ Y$ is a smooth complex manifold.

Let $ \mathcal{F} $ be an object in $ {\rm D^{b}_{coh}} \left(X\right)$ such that $ \mathcal{H} \mathcal{F} $ is concentrated in degrees $ \left[c_{0 },c_{1} \right]$.
We need to show that the cohomology of $  Rf_{*} \mathcal{F} $ is coherent and concentrated in degrees $ \left[c_{0 },c_{1} +n\right]$.

Let $ X_{f} \subset X\times Y$ be the graph of $ f$.
Then, $ X_{f} $ is a closed complex submanifold of $ X\times Y$. 
Let $ i:x\in  X\to \left(x,f\left(x\right)\right)\in X\times Y$ be the embedding associated to $ X_{f} $.

Let 
\begin{align}\label{eq:ebfu}
 \pi: X\times Y \to Y
\end{align}
be the projection onto the second factor.
Then, 
\begin{align}\label{eq:gwy1}
f= \pi \circ i.
\end{align}

Since $ i$ is a closed embedding, by \eqref{eq:gpfo}, \eqref{eq:gwy1}, and by functoriality \cite[\href{https://stacks.math.columbia.edu/tag/0D5T}{Tag 0D5T}]{stacks-project}, we have a canonical isomorphism in $ {\rm D}^{+} \left(Y\right)$, 
\begin{align}\label{eq:dxym}
  Rf_{*} \mathcal{F} =  R \pi_{*} Ri_{*} \mathcal{F} = R \pi_{*} \left(i_{*} \mathcal{F}\right). 
\end{align}

Note that  
\begin{align}\label{eq:diar}
 {\rm Supp}\left(i_{*} \mathcal{F} \right) \subset X_{f}. 
\end{align}
Since $ f$ is proper, for any compact subset $ K \subset Y$, we see that 
\begin{align}\label{eq:fedz}
  {\rm Supp}\left(i_{*} \mathcal{F} \right) \cap \pi^{-1} \left(K \right)
\end{align}
is compact.

Applying Theorem \ref{thm:fbt3}, while replacing $ \left(\widetilde{S} , \mathcal{F}\right) $ by $ \left(Y,i_{*} \mathcal{F}\right) $, we see that the cohomology of $ R \pi_{*} \left(i_{*} \mathcal{F}\right) $ is coherent and concentrated in degrees $ \left[c_{0 },c_{1} +n\right]$.
By \eqref{eq:dxym}, we get the desired result.
\qed





  
\section{Proof of Theorem \ref{thm:main} \textnormal{\ref{main2a}}) : general case}\label{S:Gg}

The purpose of this section is to establish Theorem \ref{thm:main} \ref{main2a}) in full generality.

This section is organized as follows.
In Section \ref{s:cass}, we introduce the reduction $X_{\mathrm{red}}$ and the singular locus $X_{\mathrm{sing}}$ associated with a complex analytic space $X$.

In Section \ref{s:more}, given a closed reduced complex analytic subspace $Z \subset X$, we study the derived direct image of coherent sheaves supported on $Z$ as well as the derived direct image of the restriction $Rf_{|Z*}$.
In particular, we show that the proof of Theorem \ref{thm:main} \ref{main2a}) can be reduced to the case where $X$ is reduced.

In Section \ref{s:endp}, we prove Theorem \ref{thm:main} \ref{main2a}) in full generality by induction on $\dim X$.
Using Hironaka’s resolution of singularities \cite{Hironaka64a,Hironaka64b}, together with the results of the previous sections, the proof reduces to the case where the source space is smooth or the lower-dimensional closed subspace $X_{\mathrm{sing}}$.

Finally, in Section \ref{s:cont}, applying Andreotti–Grauert’s theorem \cite{AndreottiGrauert62}, we obtain a sharper control on $c(f)$.

\subsection{Complex analytic subspaces}\label{s:cass}
Let $Z$ be a closed complex analytic subspace of $ X$. 
Let $ \mathscr{I}_{Z} $ be the associated coherent ideal sheaf with zero set $ N\left(\mathscr{I}_{Z} \right)$.
Then, 
\begin{align}\label{eq:v3v1}
  \left(Z,\mathscr{O}_{Z} \right)= \left(N\left(\mathscr{I}_{Z}  \right),\left(\mathscr{O}_{X} /\mathscr{I}_{Z}  \right)_{|N\left(\mathscr{I}_{Z}  \right)} \right).
\end{align}
 

Let $ X_{{\rm red}} $ be the reduction of $ X$.
If $ \mathscr{N}_{X} \subset \mathscr{O}_{X} $ is the nilradical of $ \mathscr{O}_{X} $, by \cite[Lemma 4.2.5, Section 4.3.2]{GrauertBook84}, then $ \mathscr{N}_{X} $ is coherent and $ X_{{\rm red}} $ is the closed complex analytic subspace associated to $ \mathscr{N}_{X} $.
Clearly, the underlying topological spaces of $ X_{{\rm red}} $ and $ X$ coincide.
If $ X= X_{{\rm red}} $ or equivalently $ \mathscr{N}_{X} = 0 $, then $ X$ is called reduced.

Let $ X_{{\rm sing}} \subset X$ be the singular locus of $ X$.
By \cite[p.~117]{GrauertBook84}, $ X_{{\rm sing}} $ is a closed reduced complex analytic subspace of $ X$.
Moreover, if $ X$ is reduced, then 
\begin{align}\label{eq:sdgz}
 \dim X_{{\rm sing}} \le  \dim X-1.
\end{align}

\subsection{More reductions via supports and cone constructions}\label{s:more}
Let $ i_{Z} : Z \to X$ be a closed reduced complex analytic subspace of $ X$.  
Let $ c_{0 } ,c_{1} \in \mathbf{Z} $ with $ c_{0 } \le c_{1} $.

\begin{defin}\label{def:3h3i}
  Let $P\left(X,Z,c_{0},c_{1}\right)$ be the class of objects $\mathcal{F} $ in $ {\rm D^{b}_{coh}}(X)$ such that the cohomology of $\mathcal{F}$ is supported on $Z$ and concentrated in degrees $\left[c_{0},c_{1}\right]$.
  
  Let $P_{-1}\left(X,Z,c_{0},c_{1}\right)$ be the subclass of $P\left(X,Z,c_{0},c_{1}\right)$ consisting of objects $\mathcal{F}$ which themselves are supported on $Z$ and concentrated in degrees $\left[c_{0},c_{1}\right]$.

  Let $P_{0}\left(X,Z,c_{0},c_{1}\right)$ be the subclass of $P\left(X,Z,c_{0},c_{1}\right)$ consisting of objects that are isomorphic in ${\rm D^{b}_ {coh}}(X)$ to elements of $P_{-1}\left(X,Z,c_{0},c_{1}\right)$.

  For $r \geq 1$, let $P_{r}\left(X,Z,c_{0},c_{1}\right)$ be the subclass of $P\left(X,Z,c_{0},c_{1}\right)$ consisting of objects $\mathcal{F}$ for which there exist $\mathcal{F}', \mathcal{F}'' \in P_{r-1}(X,Z,c_{0},c_{1})$ together with a distinguished triangle in ${\rm D^{b}_{coh}}(X)$, 
\begin{equation}\label{dia:pfai}
\begin{tikzsd}
      \mathcal{F} ' &	\mathcal{F}  & \mathcal{F}^{\prime \prime }   & \mathcal{F}^{\prime \bullet +1} .
     \arrow["", from=1-1, to=1-2]
     \arrow["", from=1-2, to=1-3]
     \arrow["", from=1-3, to=1-4]
  \end{tikzsd} 
\end{equation}

If $ Z= X_{{\rm red}}  $, for $ r\in \left\{-1\right\} \cup \mathbf{N} $, put 
\begin{align}\label{eq:dkgo}
 &P_{r} \left(X,c_{0 } ,c_{1}\right)= P_{r} \left(X,X_{{\rm red}} ,c_{0 },c_{1}\right),& P \left(X,c_{0 } ,c_{1}\right)= P \left(X,X_{{\rm red}} ,c_{0 },c_{1}\right).
\end{align}
\end{defin}

Observe that if $ r\in \mathbf{N} $ and if $ \mathcal{F}\in  P_{r} \left(X,Z,c_{0 },c_{1}\right) $, then the isomorphic class of $ \mathcal{F} $ in $ {\rm D^{b}_{coh}} \left(X\right)$ is still contained in $ P_{r} \left(X,Z,c_{0 },c_{1}\right)$.
This property does not hold when $ r= -1$.
The isomorphic class of an element in $  P_{-1} \left(X,Z,c_{0 },c_{1}\right)$ is contained in $ P_{0 } \left(X,Z,c_{0 },c_{1}\right)$.

Clearly, the class of all the objects of $ {\rm D^{b}_{coh}} \left(X\right) $ is given by 
\begin{align}\label{eq:ipsk}
 \bigcup_{c_{0 } ,c_{1} \in \mathbf{Z},c_{0 } \le c_{1}} P \left(X,c_{0 },c_{1}\right).
\end{align}
We have a filtration 
\begin{align}\label{eq:oyqz}
  P_{- 1} \left(X,Z,c_{0 } ,c_{1}\right) \subset P_{0  } \left(X,Z,c_{0 },c_{1}\right) \subset \cdots \subset P \left(X,Z,c_{0 },c_{1}\right).
\end{align}

\begin{prop}\label{prop:yewj}
The filtration \eqref{eq:oyqz} is stationary. 
More precisely, 
 \begin{align}\label{eq:zaan}
   P_{c_{1} -c_{0 } } \left(X,Z,c_{0 },c_{1}\right) = P \left(X,Z,c_{0 },c_{1}\right).
 \end{align}
\end{prop}
\begin{proof}
We can and we may assume $ c_{0 } = 0 $. 
Let $ \mathcal{F}\in P \left(X,Z,0 ,c_{1}\right)$.
We need to show 
\begin{align}\label{eq:uylq}
 \mathcal{F}\in P_{c_{1} } \left(X,Z,0 ,c_{1} \right).
\end{align}

If $ c_{1} = 0 $, then $ \mathcal{F} $ is isomorphic to 
$ \tau_{\le 0 } \tau_{\ge 0 } \mathcal{F} $ in $ {\rm D^{b}_{coh}}\left(X\right) $.
Since $ \tau_{\le 0 } \tau_{\ge 0 } \mathcal{F}= \mathcal{H}^{0} \mathcal{F} $, 
by Definition \ref{def:3h3i}, we have \eqref{eq:uylq}.

Assume now $ c_{1} \ge 1$, and that \eqref{eq:uylq} holds for $c_{1}'\le c_{1}-1$.
The canonical inclusion $ \tau_{\le c_{1} -1} \mathcal{F} \to \mathcal{F} $ induces a distinguished triangle in $ {\rm D^{b}_{coh}}\left(X\right)$, 
\begin{equation}\label{dia:v2ll1}
\begin{tikzsd}
    \tau_{\le c_{1} -1} \mathcal{F}  & \mathcal{F}  &	{\rm cone}\left(\tau_{\le c_{1} -1} \mathcal{F},  \mathcal{F}\right) &  \left( \tau_{\le c_{1} -1} \mathcal{F}\right)^{\bullet +1}. 
     \arrow["", from=1-1, to=1-2]
     \arrow["", from=1-2, to=1-3]
     \arrow["", from=1-3, to=1-4]
  \end{tikzsd} 
\end{equation}

By construction, we have 
\begin{align}\label{eq:gofp}
& \tau_{\le c_{1} -1} \mathcal{F}\in P\left(X,Z,0,c_{1} -1\right),&{\rm cone}\left(\tau_{\le c_{1} -1} \mathcal{F},  \mathcal{F}\right) \in P\left(X,Z,c_{1} ,c_{1}\right).
\end{align}
By \eqref{eq:gofp} and induction assumptions, we have 
\begin{align}\label{eq:gnu1}
 & \tau_{\le c_{1} -1} \mathcal{F}\in P_{c_{1} -1} \left(X,Z,0,c_{1} -1\right),&{\rm cone}\left(\tau_{\le c_{1} -1} \mathcal{F},  \mathcal{F}\right) \in P_{0 } \left(X,Z,c_{1} ,c_{1}\right).
\end{align}
By \eqref{dia:v2ll1} and \eqref{eq:gnu1}, we get \eqref{eq:uylq}.
This finishes the induction argument and completes the proof of our proposition.
\end{proof}

\begin{prop}\label{prop:lp3n}
 Given a distinguished triangle in $ {\rm D}^{+}  \left(X\right)$, 
 \begin{equation}\label{dia:pfai1}
\begin{tikzsd}
      \mathcal{F} ' &	\mathcal{F}  & \mathcal{F}^{\prime \prime }   & \mathcal{F}^{\prime \bullet +1},
     \arrow["", from=1-1, to=1-2]
     \arrow["", from=1-2, to=1-3]
     \arrow["", from=1-3, to=1-4]
  \end{tikzsd} 
\end{equation}
 if $ \mathcal{F} ',\mathcal{F}^{\prime \prime }\in P \left(X,Z,c_{0 },c_{1}\right)$, then  $\mathcal{F}\in  P \left(X,Z,c_{0 },c_{1}\right)$.
\end{prop}
\begin{proof}
 By \cite[(1.5.3), Proposition 1.5.6]{KashiwaraSchpira}
and by \eqref{dia:pfai1}, we get a long exact sequence of $ \mathscr{O}_{X} $-modules,
\begin{equation}\label{dia:bq3s12} 
\begin{tikzsd}
    \cdots & \mathcal{H} ^{p} \mathcal{F} '  &	\mathcal{H} ^{p} \mathcal{F}  & \mathcal{H} ^{p} \mathcal{F}^{\prime \prime }  & \cdots
     \arrow["", from=1-1, to=1-2]
     \arrow["", from=1-2, to=1-3]
     \arrow["", from=1-3, to=1-4]
     \arrow["", from=1-4, to=1-5]
  \end{tikzsd} 
\end{equation}
Since all the $ \mathcal{H} ^{p} \mathcal{F} '$ and $ \mathcal{H} ^{p} \mathcal{F}^{\prime \prime } $ are coherent, by the general properties on coherent sheaves \cite[Annex A.4]{GrauertBook84}, we see that $ \mathcal{H} ^{p}  \mathcal{F}$ is also coherent.
Moreover, by the corresponding properties on $ \mathcal{H} ^{p} \mathcal{F} '$ and $ \mathcal{H} ^{p} \mathcal{F}^{\prime \prime }$, we see that $ \mathcal{H} ^{p}  \mathcal{F}$ is supported on $ Z$ and vanishes if $ p \notin \left[c_{0} ,c_{1} \right]$.
The proof of our proposition is completed.
\end{proof}

Let $ Y$ be another complex analytic space, and let $ f: X\to Y$ be a proper holomorphic map.

\begin{prop}\label{prop:1vuw}
  Assume that there exists $ c\left(f\right)\in \mathbf{N} $ such that for all $ c_{0},c_{1}\in  \mathbf{Z}  $ with $ c_{0} \le c_{1} $,  $ Rf_{*} $ sends $ P_{-1} \left(X,Z,c_{0},c_{1} \right)$  to $ P\left(Y,c_{0},c_{1}+c\left(f\right)\right)$.
  Then, $ Rf_{*} $ sends $ P \left(X,Z,c_{0},c_{1} \right)$  to $ P\left(Y,c_{0},c_{1}+c\left(f\right)\right)$.
\end{prop}
\begin{proof}
 By Proposition \ref{prop:yewj}, we can prove our proposition by induction on $ 0 \le r \le c_{1} -c_{0 } $ that $ Rf_{*} $ maps $ P_{r} \left(X,Z,c_{0},c_{1} \right)$  to $ P\left(Y,c_{0},c_{1}+c\left(f\right)\right)$.

 If $ r= 0 $, the result follows from the fact that all the elements $ P_{0 } \left(X,Z,c_{0},c_{1} \right)$ are isomorphic in $ {\rm D^{b}_{coh}}\left(X\right) $ to some elements of $ P_{-1} \left(X,Z,c_{0},c_{1} \right)$.

 Assume now $ r\ge 1$ and that our result holds for $ r-1$. 
 If $ \mathcal{F}\in P_{r} \left(X,Z,c_{0},c_{1} \right)$ such that \eqref{dia:pfai} holds, then we have a distinguished triangle in $ {\rm D}^{+} \left(Y\right)$,
 \begin{align}\label{eq:nzup}
   \begin{tikzsd}
     Rf_{*}  \mathcal{F} ' &	Rf_{*}\mathcal{F}  & Rf_{*}\mathcal{F}^{\prime \prime }  &    \left(Rf_{*} \mathcal{F}'\right)^{\bullet +1}.
     \arrow["", from=1-1, to=1-2]
     \arrow["", from=1-2, to=1-3]
     \arrow["", from=1-3, to=1-4]
  \end{tikzsd} 
\end{align}
By induction assumptions, $ R f_{*} \mathcal{F} ', R f_{*}\mathcal{F}^{\prime \prime }\in P\left(Y,c_{0} ,c_{1}+c\left(f \right)\right)$.
By Proposition \ref{prop:lp3n} and \eqref{eq:nzup}, we get 
\begin{align}\label{eq:1afa}
R f_{*}\mathcal{F}\in P\left(Y,c_{0} ,c_{1}+c\left(f \right)\right).
\end{align}
This finishes the induction argument and completes the proof of our proposition.
\end{proof}


The holomorphic map $ f: X \to Y $ restricts to 
\begin{align}\label{eq:sm21}
 f_{|Z} :Z\to Y.
\end{align}

\begin{prop}\label{prop:izht}
  Assume that there exists $ c\left(f_{|Z} \right)\in \mathbf{N} $ such that for any $ c_{0} ,c_{1}\in \mathbf{Z}  $ with $ c_{0} \le c_{1} $, $ Rf_{|Z*} $ sends $ P\left(Z,c_{0} ,c_{1}\right)$ to $ P\left(Y,c_{0} ,c_{1}+c\left(f_{|Z} \right)\right)$.
  Then, $ Rf_{*} $ sends $ P \left(X,Z,c_{0} ,c_{1} \right)$  to $ P\left(Y,c_{0} ,c_{1}+c\left(f_{|Z} \right)\right)$.
\end{prop}
\begin{proof}
Let $ \mathcal{F} \in P_{-1} \left(X,Z,c_{0} ,c_{1} \right)$.
By Proposition \ref{prop:1vuw}, it is enough to show 
  \begin{align}\label{eq:tpg3}
   Rf_{*} \mathcal{F}\in P\left(Y,c_{0} ,c_{1}+c\left(f_{|Z}\right)\right).
  \end{align}  

As in \cite[p.~20]{GrauertBook84}, we have an exact sequence of $ \mathscr{O}_{X} $-complexes, 
\begin{equation}\label{eq:oafa} 
\begin{tikzsd}
    0  & \mathscr{I}_{Z}  \mathcal{F}  &	\mathcal{F} & i_{Z*} i^{*}_{Z}  \mathcal{F} & 0.
     \arrow["", from=1-1, to=1-2]
     \arrow["", from=1-2, to=1-3]
     \arrow["", from=1-3, to=1-4]
     \arrow["", from=1-4, to=1-5]
  \end{tikzsd} 
\end{equation}

We claim that 
\begin{align}\label{eq:s3vh}
 &   \mathscr{I}_{Z}  \mathcal{F} \in P_{-1} \left(X,Z,c_{0 } ,c_{1} \right),
&i^{*}_{Z}  \mathcal{F}  \in P_{-1} \left(Z, c_{0 } ,c_{1} \right).
\end{align}
Indeed, since  $ \mathcal{F} $ is coherent and supported on $ Z$, we know that $ \mathscr{I}_{Z}  \mathcal{F} $ is also coherent\footnote{Since $ \mathscr{I}_{Z}  $ and $ \mathcal{F} $ are coherent, they are finite generated so that $ \mathscr{I}_{Z}  \mathcal{F} $ is also finite generated.
Moreover, as a submodule of a finite relation generated modules $ \mathcal{F} $, $ \mathscr{I}_{Z} \mathcal{F} $  is also finite relation generated \cite[Annex A.3.2 Examples 2)]{GrauertBook84}.
Therefore, $ \mathscr{I}_{Z} \mathcal{F} $ is coherent.} and supported on $ Z$, so that the first equation of \eqref{eq:s3vh} holds.
By \cite[p.~19]{GrauertBook84}, $ i^{*}_{Z}  \mathcal{F}$ is still coherent on $ Z$, which gives the second equation of \eqref{eq:s3vh}.

Applying the functor $ Rf_{*} $ to \eqref{eq:oafa} and using $ Rf_{*} i_{Z*}=Rf_{*} Ri_{Z*}= R f_{|Z*}   $, by \cite[\href{https://stacks.math.columbia.edu/tag/0150}{Tag 0150},\href{https://stacks.math.columbia.edu/tag/0159}{Tag 0159}]{stacks-project}, we get a distinguished triangle in $ {\rm D}^{+} \left(Y\right)$,
\begin{equation}\label{dia:bq3s} 
\begin{tikzsd}
  R f_{*} \left(\mathscr{I}_{Z}  \mathcal{F}\right) &	Rf_{*}  \mathcal{F} & R f_{|Z*} \left(i^{*}_{Z} \mathcal{F} \right)  &    \left(R f_{*} \left(\mathscr{I}_{Z}  \mathcal{F}\right)\right)^{\bullet +1}.
     \arrow["", from=1-1, to=1-2]
     \arrow["", from=1-2, to=1-3]
     \arrow["", from=1-3, to=1-4]
  \end{tikzsd} 
\end{equation}

By our assumption on $ Rf_{|Z*} $, and by the second equation of \eqref{eq:s3vh}, we have 
\begin{align}\label{eq:crbo}
R f_{|Z*} \left(i^{*}_{Z}  \mathcal{F} \right)\in P\left(Y,c_{0 } ,c_{1}+c\left(f_{|Z} \right)\right).
\end{align}
By Proposition \ref{prop:lp3n}, and by \eqref{dia:bq3s}, \eqref{eq:crbo}, we know that \eqref{eq:tpg3} is a consequence of
\begin{align}\label{eq:masc1}
   R f_{*} \left(\mathscr{I}_{Z}  \mathcal{F}\right) \in P\left(Y,c_{0 } ,c_{1}+c\left(f_{|Z} \right)\right).
\end{align}

By the first equation of \eqref{eq:s3vh}, $ \mathscr{I}_{Z} \mathcal{F} $ has the same property as $ \mathcal{F} $.
We can repeat the previous argument with $\mathcal{F}$ replaced by $\mathscr{I}_{Z}\mathcal{F}$.
Iterating this procedure, it suffices to show that there exists $k \in \mathbf{N}$ such that
\begin{align}\label{eq:masc}
   R f_{*} \left(\mathscr{I}_{Z}^{k} \mathcal{F}\right) \in P\left(Y,c_{0 } ,c_{1}+c\left(f_{|Z} \right)\right).
\end{align}

By Proposition \ref{prop:dpzs}, we know that \eqref{eq:tpg3} is a local property on $ Y$.
Combining this with the above, we need only to show that for any relatively compact open subset $ V$ of $ Y$, there exists $ k\in \mathbf{N}$ such that
\begin{align}\label{eq:qmae}
  \left(Rf_{*}  \left(\mathscr{I}_{Z}^{k}   \mathcal{F}\right)\right)_{|V}\in P\left(V,c_{0 } ,c_{1}+c\left(f_{|Z} \right)\right). 
\end{align}

Since $ V$ is relatively compact in $ Y$, since $ f$ is proper, $ U= f^{-1} \left(V\right)$ is relatively compact in $ X$.
By Rückert Nullstellensatz \cite[Section 3.2.2]{GrauertBook84}, there exists $ k \in \mathbf{N} $ such that 
\begin{align}\label{eq:1too}
\left( \mathscr{I}^{k}_{Z}  \mathcal{F} \right)_{|U} = 0.
\end{align}
By Proposition \ref{prop:dpzs} and by \eqref{eq:1too}, we get 
\begin{align}\label{eq:q1t2}
 \left(Rf_{*}  \left(\mathscr{I}_{Z}^{k}   \mathcal{F}\right)\right)_{|V}= 0,
\end{align}
which implies \eqref{eq:qmae}.
The proof of our proposition is complete.
\end{proof}

\begin{cor}\label{cor:tnlv}
 If the conclusion of Theorem \ref{thm:main} \ref{main2a}) holds for $f_{|X_{{\rm red}} } : X_{{\rm red}} \to Y$, then the conclusion of Theorem \ref{thm:main} \ref{main2a}) holds for $ f: X\to Y$ and 
 \begin{align}\label{eq:lm1i}
  c\left(f\right)\le c\left(f_{|X_ {\rm red}} \right).
 \end{align}
\end{cor}
\begin{proof}
 It is enough to apply Proposition \ref{prop:izht} to $ Z= X_{{\rm red}} $.
\end{proof}

\subsection{End of the proof for Theorem \ref{thm:main} \ref{main2a})}\label{s:endp}
By Remark \ref{re:fq3s} and Corollary \ref{cor:tnlv}, we can assume that $ X$ is reduced and has second countable topology.

Recall that we assume $ \dim X < + \infty $.
Let us prove by induction on $  \dim X\in \mathbf{N} $ that for any complex analytic space $ Y$ and any proper holomorphic map $ f: X \to Y$,  there exists $ c\left(f\right)\in \mathbf{N} $ such that for all $ c_{0 }, c_{1}\in \mathbf{Z} $ with $ c_{0 } \le c_{1} $, and $ \mathcal{F}\in P \left(X,c_{0 } ,c_{1} \right)$, we have 
\begin{align}\label{eq:zmmd}
  Rf_{*} \mathcal{F} \in P\left(Y,c_{0 } ,c_{1}+c\left(f\right)\right).
\end{align}

If $ \dim X=0$, then $ X$ is discrete.
Equation \eqref{eq:zmmd} holds trivially with $ c\left(f\right)= 0 $.

Assume now  $ \dim X \ge 1$ and that \eqref{eq:zmmd} holds  when $ X$ is replaced by any complex analytic spaces $ X^{\prime } $  such that $ \dim X^{\prime } \le \dim X-1$. 
By Section \ref{S:Gs} and by Proposition \ref{prop:1vuw}, we can assume that $ X$ is not smooth and
\begin{align}\label{eq:rhih}
\mathcal{F}\in P_{-1} \left(X,c_{0 } ,c_{1} \right). 
\end{align}

Since $ X$ has second countable topology, by Hironaka's resolution of singularities \cite{Hironaka64a,Hironaka64b}, \cite[Theorem 5.4.2]{ArocaHironakaVice18}, there exist
a smooth complex manifold $ X_{1} $, a normal crossing divisor $D_{1}  \subset  X_{1} $, and   
a proper holomorphic map,  
\begin{align}\label{eq:kru1}
\pi_{1} : X_{1} \to X,
\end{align}
such that $ D_{1} =   \pi_{1}^{-1} \left(X_{{\rm sing}} \right)$ and $ \pi_{1} $ restricts to a biholomorphic map, 
\begin{align}\label{eq:ulhe}
  X_{1} \backslash D_{1} \to X \backslash X_{{\rm sing}}.
\end{align}

By \eqref{eq:rhih}, as in the second equation of \eqref{eq:s3vh}, we have 
\begin{align}\label{eq:h3di}
  \pi_{1}^{*} \mathcal{F}\in P_{-1} \left(X_{1}, c_{0 }, c_{1} \right).
\end{align}
Since $ X_{1} $ is smooth and $ \dim X_{1} = \dim X$, by Section \ref{S:Gs} and \eqref{eq:h3di}, we have 
\begin{align}\label{eq:2kfo}
  &R\pi_{1*} \left( \pi_{1}^{*} \mathcal{F}\right)\in P\left(X,c_{0 },c_{1}+\dim X\right), & R \left(f \pi_{1} \right)_{*} \left(\pi_{1} ^{*} \mathcal{F} \right)\in P\left(Y,c_{0 },c_{1}+\dim X\right).
\end{align}

Let us construct a morphism in $ {\rm D }^{+}  \left(X\right)$, 
\begin{align}\label{eq:ey3f}
  \mathcal{F} \to R\pi_{1*} \left( \pi_{1}^{*} \mathcal{F}\right). 
\end{align}
Indeed, we have a natural morphism of $ \mathscr{O}_{X} $-complexes,
\begin{align}\label{eq:ghwo}
\mathcal{F} \to \pi_{1*}  \left(\pi_{1}^{*} \mathcal{F}\right). 
\end{align}
Let 
\begin{align}\label{eq:pgfx}
  \pi_{1}^{*} \mathcal{F} \to \mathcal{I} _{1} 
\end{align}
be a right injective resolution of $ \pi_{1}^{*} \mathcal{F} $ as in \eqref{eq:evfv}.
We have a canonical isomorphism in $ {\rm D}^{+}  \left(X\right)$,
\begin{align}\label{eq:iyjj}
 \pi_{1*}  \mathcal{I}_{1} \simeq R\pi_{1*} \left( \pi_{1}^{*} \mathcal{F}\right).
\end{align}
By \eqref{eq:ghwo}-\eqref{eq:iyjj}, we get \eqref{eq:ey3f}.
 
We extend \eqref{eq:ey3f} to a distinguished triangle in $ {\rm D^{+}}\left(X\right)$,
 \begin{align}\label{eq:ffuf}
   \begin{tikzsd}
     \mathcal{F} ' &	\mathcal{F}  & R\pi_{1*} \left(\pi_{1}^{*} \mathcal{F}\right)  &     \mathcal{F} ^{\prime \bullet +1}.
     \arrow["", from=1-1, to=1-2]
     \arrow["", from=1-2, to=1-3]
     \arrow["", from=1-3, to=1-4]
  \end{tikzsd} 
\end{align}

We claim that 
\begin{align}\begin{aligned}\label{eq:qyzr}
  \mathcal{F}'\in& P\left(X,X_{{\rm sing}}, c_{0 }, c_{1}+\dim X +1\right),\\
  R f_{*} \mathcal{F} ' \in& P\left(Y,c_{0 }, c_{1} +c\left(f_{|X_{{\rm sing}} }  \right) + \dim X + 1\right).
\end{aligned}\end{align}
Indeed, by Proposition \ref{prop:lp3n}, \eqref{eq:rhih}, the first equation of \eqref{eq:2kfo}, and \eqref{eq:ffuf}, we have 
\begin{align}\label{eq:iwuo}
\mathcal{F}'\in  P\left(X,c_{0 }, c_{1}+\dim X +1\right).
\end{align}
Since \eqref{eq:ulhe} is biholomorphic, we see that \eqref{eq:ghwo} is isomorphic on $ X\backslash X_{{\rm sing}} $, so that \eqref{eq:ey3f} restricts to an isomorphism in ${\rm D}^{+} \left( X\backslash X_{{\rm sing}}\right)$.
Combining this with \eqref{eq:iwuo}, we get the first equation in \eqref{eq:qyzr}.
By \eqref{eq:sdgz}, we can apply the induction assumption to $ f_{|X_{{\rm sing}}  }: X_{{\rm sing}} \to Y$.
By Proposition \ref{prop:izht} and by the first equation in \eqref{eq:qyzr}, we get the second equation in \eqref{eq:qyzr}.

 
Applying the functor $ Rf_{*} $ to \eqref{eq:ffuf}, we get a distinguished triangle in $ {\rm D}^{+} \left(Y\right)$,
 \begin{align}\label{eq:tclh}
   \begin{tikzsd}
     Rf_{*} \mathcal{F} ' &	Rf_{*}   \mathcal{F} & R\left(f\pi_{1} \right)_{*} \left(\pi_{1}^{*} \mathcal{F}\right)   &   Rf_{*} \mathcal{F}^{\prime \bullet +1}.
     \arrow["", from=1-1, to=1-2]
     \arrow["", from=1-2, to=1-3]
     \arrow["", from=1-3, to=1-4]
  \end{tikzsd} 
\end{align}
By Proposition \ref{prop:lp3n}, by the second equation in \eqref{eq:2kfo} and the second equation in \eqref{eq:qyzr}, we get \eqref{eq:zmmd} with 
\begin{align}\label{eq:pohi}
 c\left(f\right)\le c\left(f_{|X_{{\rm sing}} }  \right) + \dim X + 1.
\end{align}
This finishes the induction argument and completes the proof of our theorem. \qed

\subsection{The constant $ c\left(f\right)$ }\label{s:cont}
The estimate on $ c\left(f\right)$ given in \eqref{eq:pohi} is far from optimal.

\begin{prop}\label{prop:a2v3}
 For any proper holomorphic map $ f: X\to Y$, we have 
 \begin{align}\label{eq:qwiq}
  c\left(f\right) \le  \dim X.
 \end{align} 
\end{prop}
\begin{proof}
Indeed, if $ \mathcal{F} $ is a coherent sheaf on $ X$, by \cite[p.~36]{GrauertBook84}, for all $ p\in \mathbf{N} $, $ R^{p} f_{*}\mathcal{F}  $ is the sheaf associated to the presheaf defined for an open set $ V \subset Y$ by
\begin{align}\label{eq:pszx}
 H^{p} \left(f^{-1} \left(V\right), \mathcal{F}_{|f^{-1} \left(V\right)} \right).
\end{align}
By Andreotti-Grauert’s theorem \cite[IX (4.15) Corollary]{DemaillyBook}, if $ p \ge \dim X+1$, \eqref{eq:pszx} vanishes, so that 
\begin{align}\label{eq:mx3z}
 R^{p} f_{*} \mathcal{F} = 0. 
\end{align}

By \eqref{eq:mx3z}, using \cite[\href{https://stacks.math.columbia.edu/tag/07K7}{Tag 07K7}]{stacks-project} or an induction argument similar to the one given in the proof of Proposition \ref{prop:yewj}, we know that for any $ c_{0 }, c_{1}\in \mathbf{Z} $ with $ c_{0 } \le c_{1} $, and $ \mathcal{F} \in P_{-1} \left(X,c_{0 } ,c_{1} \right)$, if $ i \ge c_{1} +\dim X+1$, we have 
\begin{align}\label{eq:ncjv}
 R^{i} f_{*} \mathcal{F} = 0. 
\end{align}
By Proposition \ref{prop:1vuw}, we know that \eqref{eq:ncjv} holds for all $\mathcal{F} \in P\left(X,c_{0 } ,c_{1} \right)$.
Therefore, we get \eqref{eq:qwiq} and finish the proof of our proposition.
\end{proof}

\bibliographystyle{amsalpha}

\def\cprime{$'$}
\providecommand{\bysame}{\leavevmode\hbox to3em{\hrulefill}\thinspace}
\providecommand{\MR}{\relax\ifhmode\unskip\space\fi MR }
\providecommand{\MRhref}[2]{%
  \href{http://www.ams.org/mathscinet-getitem?mr=#1}{#2}
}
\providecommand{\href}[2]{#2}

\end{document}